\documentclass[12pt]{article}
\usepackage{amssymb, amsfonts, amsmath}
\usepackage{amsthm, thmtools, mathrsfs}
\usepackage{enumitem, textcomp, undertilde, yfonts, yhmath, ushort}
\usepackage[colorlinks=true,linkcolor=blue]{hyperref}
\usepackage[retainorgcmds]{IEEEtrantools}

\declaretheorem[style=definition,qed=$\dashv$,numberwithin=section]
{definition}

\declaretheorem[style=plain,sibling=definition]{theorem}
\declaretheorem[style=plain,sibling=definition]{lemma}

\declaretheorem[style=plain,sibling=definition]{corollary}
\declaretheorem[style=definition,sibling=definition]{remark}
\declaretheorem[style=plain,sibling=definition]{claim}
\declaretheorem[style=plain,sibling=definition]{claim*}
\declaretheorem[style=plain,sibling=definition]{subclaim}

\declaretheorem[style=definition,sibling=definition]{notation}
\declaretheorem[style=plain,sibling=definition]{fact}
\numberwithin{equation}{section}

\newcommand{\rest}{\mathord{\upharpoonright}}

\newcommand{\tu}{\textup}
\newcommand{\tulp}{\tu{(}}
\newcommand{\turp}{\tu{)}}

\newcommand{\univ}[1]{\lpole #1\rpole}
\newcommand{\down}{\mathord{\downarrow}}
\newcommand{\PFA}{\mathsf{PFA}}
\newcommand{\dom}{{\rm dom}}
\newcommand{\cp}{{\rm cp }}

\newcommand{\lh}{{\rm lh}}

\newcommand{\hodlower}{\mathrm{hod}}

\newcommand{\card}{\mathrm{card}}
\newcommand{\forces}{\Vdash}

\renewcommand{\models}{\vDash}
\newcommand{\powerset}{\mathfrak{P}}

\def\P{{\mathcal{P}}}
\def\W{{\mathcal{W}}}
\def\Q{{\mathcal{Q}}}

\def\L{{\rm{L}}}

\def\R{{\mathcal{R}}}
\def\H{{\rm{HOD}}}
\def\M{{\mathcal{M}}}
\def\N{{\mathcal{N}}}

\def\T {{\mathcal{T}}}

\def\Ss{{\mathcal{S}}}

\newcommand{\ins}{\trianglelefteq}
\newcommand{\nins}{\ntrianglelefteq}
\newcommand{\pins}{\triangleleft}
\newcommand{\npins}{\ntriangleleft}
\newcommand{\Hull}{\mathrm{Hull}}
\newcommand{\cHull}{\mathrm{cHull}}
\newcommand{\crit}{\mathrm{crit}}
\newcommand{\RR}{\mathbb{R}}
\newcommand{\rSigma}{{\mathrm{r}\Sigma}}

\newcommand{\J}{\mathcal J}
\newcommand{\un}{\cup}
\newcommand{\sub}{\subseteq}
\newcommand{\om}{\omega}
\newcommand{\ZF}{\mathsf{ZF}}
\newcommand{\AD}{\mathsf{AD}}
\newcommand{\ZFC}{\mathsf{ZFC}}
\newcommand{\CC}{\mathbb{C}}
\newcommand{\es}{\mathbb{E}}
\newcommand{\Tt}{\mathcal{T}}

\newcommand{\Ttdot}{\dot{\Tt}}
\newcommand{\PP}{\mathbb{P}}
\newcommand{\Vv}{\mathcal{V}}
\newcommand{\Uu}{\mathcal{U}}
\newcommand{\Uudot}{{\dot{\Uu}}}
\newcommand{\conc}{\ \widehat{\ }\ }
\newcommand{\st}{\ \|\ }
\newcommand{\Lp}{\mathrm{Lp}}
\newcommand{\sats}{\vDash}
\newcommand{\cut}{\backslash}
\newcommand{\Ult}{\mathrm{Ult}}
\newcommand{\core}{\mathfrak{C}}
\newcommand{\pow}{\powerset}
\newcommand{\inter}{\cap}
\newcommand{\xibar}{\bar{\xi}}
\newcommand{\lex}{\mathrm{lex}}
\newcommand{\car}{\mathrm{card}}
\newcommand{\BB}{\mathbb{B}}
\newcommand{\Coll}{\mathrm{Col}}
\newcommand{\rg}{\mathrm{rg}}
\newcommand{\lpole}{\left\lfloor}
\newcommand{\rpole}{\right\rfloor}
\newcommand{\eps}{\epsilon}

\renewcommand{\l}{\mathit{l}}
\newcommand{\Hh}{\mathcal{H}}

\newcommand{\OR}{\mathrm{o}}
\newcommand{\nth}{\mathrm{th}}
\newcommand{\Ord}{\mathrm{Ord}}
\newcommand{\Th}{\mathrm{Th}}
\newcommand{\com}{\circ}
\newcommand{\HC}{\mathrm{HC}}
\newcommand{\cof}{\mathrm{cof}}
\newcommand{\BBB}{\mathfrak{B}}
\newcommand{\Ll}{\mathcal{L}}
\newcommand{\all}{\forall}
\newcommand{\ex}{\exists}
\newcommand{\Ttvec}{\vec{\Tt}}
\newcommand{\Uuvec}{\vec{\Uu}}
\newcommand{\id}{\mathrm{id}}
\newcommand{\cross}{\times}
\newcommand{\Mbar}{\bar{\M}}
\newcommand{\gammabar}{\bar{\gamma}}
\newcommand{\Ttbar}{\bar{\Tt}}

\newcommand{\hod}{\mathrm{HOD}}
\newcommand{\elem}{\preccurlyeq}

\newcommand{\pred}{\text{pred}}
\newcommand{\MC}{\mathsf{MC}}
\newcommand{\OD}{\mathrm{OD}}
\newcommand{\alphabar}{\bar{\alpha}}
\newcommand{\pP}{\mathfrak{P}}
\newcommand{\C}{\mathcal{C}}
\newcommand{\opbk}{\mathscr{B}}
\newcommand{\her}{\mathscr{H}}
\newcommand{\plI}{\mathrm{I}}
\newcommand{\plII}{\mathrm{II}}
\newcommand{\ZFmin}{\mathsf{ZF}^-}
\newcommand{\rank}{\mathrm{rank}}
\newcommand{\psub}{\subsetneq}
\newcommand{\restpred}{\wr}
\newcommand{\MFsharp}{\mathfrak{M}}

\newcommand{\Fop}{\mathcal{F}}

\newcommand{\unique}{\iota}
\newcommand{\g}{\mathrm{g}}
\newcommand{\gTheta}{\mathsf{G}}
\newcommand{\sq}{\mathrm{sq}}

\newcommand{\trancl}{\mathrm{trancl}}

\newcommand{\xvec}{\vec{x}}

\newcommand{\G}{\mathcal{G}}
\newcommand{\gdot}{\dot{g}}

\newcommand{\Ii}{\mathcal{I}}
\newcommand{\maps}{\colon}

\newcommand{\Gammabar}{{\bar{\Gamma}}}

\newcommand{\Aa}{\mathcal{A}}
\newcommand{\betabar}{\bar{\beta}}
\newcommand{\bfrSigma}{\utilde{\rSigma}}

\newcommand{\bfSigma}{\utilde{\Sigma}}
\newcommand{\rPi}{{\mathrm{r}\Pi}}
\newcommand{\bfrPi}{\utilde{\rPi}}

\newcommand{\Ww}{\mathcal{W}}
\newcommand{\Avec}{\vec{A}}
\newcommand{\kappabar}{\bar{\kappa}}

\newcommand{\LpGbgF}{\overline{\Lp}}
\newcommand{\HCc}{\mathrm{cd}}
\mathchardef\mhyphen="2D
\newcommand{\KP}{\mathsf{KP}}
\newcommand{\tame}{\mathrm{t}}
\newcommand{\AC}{\mathsf{AC}}
\newcommand{\DC}{\mathsf{DC}}

\newcommand{\partialto}{\dashrightarrow}

\newcommand{\QQ}{\mathbb{Q}}

\newcommand{\pvec}{\vec{p}}

\newcommand{\eqdef}{=_{\mathrm{def}}}
\newcommand{\onto}{\stackrel{\mathrm{onto}}{\to}}
\newcommand{\cofto}{\stackrel{\mathrm{cof}}{\to}}

\renewcommand{\L}{\mathcal{L}}

\newcommand{\witri}{\widetriangle}
\newcommand{\hmE}{E}
\newcommand{\hmP}{P}

\newcommand{\hmb}{cb}
\newcommand{\hmp}{cp}

\newcommand{\hmmu}{\mu}
\newcommand{\hme}{e}
\newcommand{\hmEdot}{\dot{\hmE}}
\newcommand{\hmPdot}{\dot{\hmP}}

\newcommand{\hmbdot}{\dot{\hmb}}
\newcommand{\hmpdot}{\dot{\hmp}}
\newcommand{\hmmudot}{\dot{\hmmu}}
\newcommand{\hmedot}{\dot{\hme}}
\newcommand{\actext}{F}
\newcommand{\hmPvec}{{\vec{\hmP}}}

\newcommand{\Pvec}{\vec{P}}

\newcommand{\Lim}{\mathrm{Lim}}
\newcommand{\DD}{\mathbb{D}}
\newcommand{\Wwvec}{\vec{\Ww}}
\newcommand{\Sigmavec}{\vec{\Sigma}}

\newcommand{\Vvvec}{{\vec{\Vv}}}
\newcommand{\sigmavec}{\vec{\sigma}}
\newcommand{\Xx}{\mathcal{X}}
\newcommand{\Yy}{\mathcal{Y}}
\newcommand{\Zz}{\mathcal{Z}}
\newcommand{\NNN}{\mathfrak{N}}
\newcommand{\MMM}{\mathfrak{M}}
\newcommand{\gOmega}{{{^\g}\Omega}}
\newcommand{\rmG}{\mathrm{G}}
\newcommand{\GOmega}{{{^\rmG}\Omega}}
\newcommand{\wfp}{\mathrm{wfp}}
\newcommand{\wvec}{\vec{w}}
\newcommand{\IND}{\mathrm{IND}}
\newcommand{\Psishort}{\Psi^{\mathrm{sh}}}
\newcommand{\Ff}{\mathcal{F}}
\newcommand{\hmPsi}{\Psi}
\newcommand{\hmPsidot}{\dot{\Psi}}

\newcommand{\etavec}{\vec{\eta}}
\newcommand{\hpm}{\mathrm{hpm}}
\newcommand{\dfnemph}{\textbf}

\begin{document}
\title{Scales in hybrid mice over $\RR$\renewcommand{\thefootnote}{\fnsymbol{footnote}} 
\renewcommand{\thefootnote}{\arabic{footnote}}}
\author{Farmer Schlutzenberg\footnote{farmer.schlutzenberg@gmail.com, Universit\"at M\"unster, 
Germany}\\
Nam Trang\footnote{ntrang@math.uci.edu, UC Irvine, California, USA}}
\maketitle

\begin{abstract}
We analyze scales in $\Lp^{\GOmega}(\RR,\Omega\rest\HC)$, the stack of sound, projecting, 
$\Theta$-g-organized $\Omega$-mice over $\Omega\rest\HC$, where $\Omega$ is either an 
iteration strategy or an operator, $\Omega$ has appropriate condensation properties, and 
$\Omega\rest\HC$ is self-scaled. This builds on Steel's analysis of scales in $L(\RR)$ and 
$\Lp(\RR)$ (also denoted $K(\RR)$). As in Steel's analysis, we work from optimal determinacy 
hypotheses. One of the main applications of the work is in
the core model induction.
\end{abstract}
\section{Introduction}\label{sec:intro}
There has been significant progress made in the core model induction in recent
years. Pioneered by W. H. Woodin and further developed
by J. R. Steel, R. D. Schindler and others, it is a powerful method for 
obtaining
lower-bound consistency strength for a large class of theories. One of the key
ingredients is the scales analysis in $L(\RR)$, and further, in $\Lp(\RR)$ (also denoted $K(\RR)$);
see Steel's \cite{Scales}, \cite{ScalesK(R)} and
\cite{Scalesweakgap}. Applications include Woodin's proof of
$\AD^{L(\RR)}$ from an $\omega_1$-dense ideal on $\omega_1$ and
Steel's proof that
$\PFA$ implies $\AD^{L(\RR)}$, amongst many others.

To use the core model induction for stronger results (for example, to 
construct models of ``$\AD^+ + \Theta > \Theta_0$") one would like to have a scales analysis
for \emph{hybrid} mice over $\RR$ -- structures beyond $\Lp(\RR)$. In this paper we present such an 
analysis. There have been recent works that make use of
methods and results from this paper, for example \cite{trang2013strength},
\cite{sargsyan2013covering}, and \cite{sargsyan2012non}.

This paper owes a strong debt to Steel's scale constructions in \cite{Scales}, 
\cite{ScalesK(R)} and
\cite{Scalesweakgap},
and to Sargsyan's notion of \emph{reorganized hod premouse}, \cite[\S 3.7]{hod_mice}.
Indeed, these are the two main components, and the main work here is in putting them together.

For the purposes mentioned above, one would particularly like to have a scales analysis 
for something like
$\Lp^{\Sigma}(\RR)$, the stack of ``projecting 
$\Sigma$-mice over $\RR$'', where $\Sigma$ is an iteration strategy with hull condensation. 
Unfortunately, 
the usual definition\footnote{Roughly, that is: Given $\Sigma$-premice 
$\N\ins\M$, with $\N$ reasonably closed, and letting $\Tt$ be the $<_\N$-least 
iteration tree for which $\N$ lacks instruction regarding the 
branch $b=\Sigma(\Tt)$, then $b$ is the next piece of information fed in to 
$\M$ after $\N$.} of 
``$\Sigma$-premouse over $\RR$'' doesn't make sense, because $\RR$ 
is not wellordered. One might try to get around this particular issue by 
arranging $\Sigma$-premice by simultaneously 
feeding in multiple branches instead of feeding them in one by one. But 
it seems difficult 
to define an amenable predicate achieving this, as discussed in 
\ref{rem:bad_non-wo_spm_definition}. 
Even if one could arrange this amenably, the scale constructions in 
\cite{ScalesK(R)} and \cite{Scalesweakgap} do not appear to generalize 
well with such an approach, because of their dependence on the 
close relationship between a mouse 
over $\RR$ and its local $\hod$.

We deal with these problems here by using the hierarchy 
of 
\emph{$\Theta$-g-organized $\Sigma$-premice}, a kind of strategy premouse.
The definition is a simple variant of \emph{g-organization}, which is essentially due to
Sargsyan;
its main content is just that of the \emph{reorganization} of hod premice.
We similarly define \emph{\tu{(}$\Theta$-\tu{)}$\g$-organized $\Fop$-premice} for operators $\Fop$,
where \emph{operators} are defined in \cite{operators}.\footnote{Many readers will probably be 
comfortable reading the present paper without knowledge of \cite{operators}, as the 
particulars of \cite{operators} are not strongly related to our purposes here. In fact, one 
could completely ignore the role of operators and focus entirely on strategy mice, without 
losing any of the main ideas. There is significant overlap between 
\cite{operators} and \S2 of the present paper. For better readability,
the common themes are generally presented in both papers. A few things are omitted in one, but
can be seen in the other.}
Given either $\Omega=\Sigma$ or $\Omega=\Fop$ as above, we only define ($\Theta$-)$\g$-organization 
assuming that $(\Omega,X)$ is \emph{nice} for some $X\in\HC$; this demands 
both a degree of \emph{condensation} and of \emph{generic determination} of $\Omega$; see 
\ref{dfn:determines}.

Given a 
nice $(\Omega,X)$ and \emph{self-scaled} $\Upsilon\sub\HC$
(see \ref{dfn:self-scaled}; this 
holds for $\Upsilon=\emptyset$) we define $\Lp^{\GOmega}(\RR,\Upsilon)$
as the stack of all sound, countably iterable $\Theta$-g-organized $\Omega$-premice built over 
$(\HC,\Upsilon)$, projecting to $\RR$. We 
will 
analyze scales in this structure. If $\Upsilon=\Omega\rest\HC$, the 
analysis can be done from optimal determinacy assumptions.
We remark that when $\Lp^\Omega(\RR,\Upsilon)$ is well-defined (such as when 
$\Omega$ is a mouse operator), we usually have
$\Lp^\Omega(\RR,\Upsilon)\neq\Lp^{\GOmega}(\RR,\Upsilon)$, but
if $\Omega$ \emph{relativizes well} 
(or something similar to this; see \cite[Definition 1.3.21(?)]{cmi}), the two hierarchies agree on 
their $\powerset(\RR)$, and actually have identical extender sequences (see 
\ref{rem:S-construction}).

The scale constructions themselves are mostly a fairly straightforward generalization of 
Steel's work in \cite{Scales}, \cite{ScalesK(R)}, \cite{Scalesweakgap}; we assume that the 
reader is
familiar with these.\footnote{One needs familiarity with said papers for 
\S\S\ref{sec:HOD_J_analysis},\ref{sec:scales} of this paper. If 
the reader has familiarity with just \cite{Scales}, one might read the present paper, referring to 
\cite{ScalesK(R)} and \cite{Scalesweakgap} as needed to fill in details we omit here.} Let 
$(\Omega,X)$ be nice and $\Upsilon$ self-scaled, and 
let $\M$ end a weak gap of 
$\Lp^{\GOmega}(\RR,\Upsilon)$. The construction of new scales over such $\M$ breaks into three 
cases, 
covered in Theorems \ref{gapInTheMouse}, \ref{thm:self_anal_coded} and \ref{weak gap from MC}; 
these are analogous to \cite[Theorems 4.16, 4.17]{ScalesK(R)} and \cite[Theorem 0.1]{Scalesweakgap} 
respectively. Thus, for the first we must assume that $\J(\M)\sats\AD$. In the context of our 
primary application (core model induction), this assumption \emph{will} hold \emph{if}
$\Omega\rest\HC\notin\M|\alpha$ and there are \emph{no divergent $\AD$ pointclasses}; 
see \ref{rem:application}. For the latter two the determinacy assumption is just that $\M\sats\AD$,
but there are also other assumptions necessary.
If $\Upsilon=\Omega\rest\HC$ then the latter two theorems cover all weak gaps, 
and so one never 
needs to assume that $\J(\M)\sats\AD$.

We won't reproduce all the details of the proofs in \cite{ScalesK(R)} and \cite{Scalesweakgap}, but 
will focus on the new features, and fill in some omissions. The most significant of 
the new features are as follows.
First, we must generalize the local $\hod$ analysis of a 
level $\M$ of 
$\Lp(\mathbb{R})$ to that of a 
level $\M$ of $\Lp^{\GOmega}(\RR,\Upsilon)$. As in \cite{ScalesK(R)}, we
establish a level-by-level fine-structural correspondence between 
$\mathcal{H}$, the local $\hod$ of $\M$, and $\M$ itself, above $\Theta^\M$. 
The fact that we are 
using $\Theta$-$\g$-organization is very important in establishing this 
correspondence (and as for $\Lp(\RR)$, the correspondence itself is very 
important in the scales analysis).
Second, an issue not dealt with in \cite{Scalesweakgap}, but with which we deal 
here, is that
a short tree $\T$ on a $k$-suitable premouse $\N$ may introduce
Q-structures
with extenders overlapping $\delta(\T)$. Third, a new case arises in the scale constructions --
at the end of a gap $[\alpha,\beta]$ of $\M$ where $\M|\beta$ is \emph{$\hmP$-active};
that is, strategy information is encoded in the predicate of $\M|\beta$. (It seems this case could 
have been avoided, however, if we had arranged our strategy premice slightly differently.)

The paper is organized as follows. In \S\ref{StrategyPremice} we discuss 
\emph{strategy 
premice} (in the sense of \emph{iteration strategy}) in 
detail, give a new presentation of these, and prove some condensation 
properties thereof, assuming that the iteration strategy involved has hull condensation and has a 
simply definable domain. In \S\ref{G-organized} we discuss
$\g$-organized and $\Theta$-$\g$-organized $\Omega$-premice, and prove related condensation facts.
In \S\ref{sec:HOD_J_analysis} we analyse the local $\hod$ of 
$\M\pins\Lp^{\GOmega}(\RR,\Upsilon)$. In
\S\ref{sec:scales} we analyse the pattern of scales in 
$\Lp^{\GOmega}(\RR,\Upsilon)$.\\

\subsection{Conventions and Notation}
We use \dfnemph{boldface} to indicate a term being defined (though when we define symbols,
these are in their normal font).
Citations such as \cite[Theorem 3.1(?)]{iter_for_stacks} are used to indicate a 
referent that may change in time -- that is, at the time of writing, 
\cite{iter_for_stacks} is a 
preprint and its Theorem 3.1 is the intended referent.

We work under $\ZF$ throughout the paper, indicating choice assumptions where we use 
them ($\DC_\RR$ in particular will be assumed for various key facts). We write $\DC_A$ for 
the restriction of $\DC$ to relations on $A$. $\Ord$ 
denotes the class of 
ordinals. Given a transitive set $M$, $\OR(\M)$ denotes $\Ord\inter M$. We write $\card(X)$ for the
cardinality of $X$, $\pow(X)$ for the power set of $X$, and for $\theta\in\Ord$,
$\powerset(<\theta)$ is the set of bounded subsets of $\theta$ and $\her_\theta$ the set 
of sets hereditarily of size $<\theta$. 
We write $f:X\partialto Y$ to denote a partial function.

We identify $\in[\Ord]^{<\om}$ with the strictly decreasing sequences of ordinals,
so given $p,q\in[\Ord]^{<\om}$, $p\rest i$ denotes the upper $i$ elements of $p$,
and $p\ins q$ means that $p=q\rest i$ for some $i$,
and $p\pins q$ iff $p\ins q$ but $p\neq q$.
The default ordering of $[\Ord]^{<\om}$ is lexicographic (largest element first),
with $p<q$ if $p\pins q$.

Let $\M=(X,A_1,\ldots)$ be a first-order structure with universe $X$ and predicates, 
constants, etc, 
$A_1,\ldots$. We write $\univ{\M}$ for $X$.
If $\Ll$ is the first-order language of $\M$, then definability over $\M$ uses $\Ll$,
unless otherwise specified. If $\Ll'\sub\Ll$, then, for example, $\Sigma_1(\Ll')$ denotes the 
$\Sigma_1$ formulas of $\Ll'$, and if $X\sub\M$, then $\Sigma_1^\M(\Ll',X)$ denotes the relations 
which are $\Sigma_1(\Ll')$-definable over $\M$ from parameters in $X$.
A \dfnemph{transitive structure} is a first-order structure with
with transitive universe. We sometimes blur the distinction between the terms \emph{transitive} and 
\emph{transitive structure}. For example, when we refer to a transitive structure as being 
\dfnemph{rud closed}, it means that its universe is rud closed. For $\M$ a
transitive structure, $\OR(\M)=\OR(\univ{\M})$. An arbitrary transitive set $X$ is also considered 
as the transitive structure $(X)$.
We write $\trancl(X)$ for the transitive closure of $X$.

Given a transitive structure $\M$, we write 
$\J_\alpha(\M)$ for the 
$\alpha^\nth$ step in Jensen's $\J$-hierarchy over $\M$ (for example,
$\J_1(\M)$ is the rud closure of $\trancl(\{\M\})$. We similarly use $\Ss$ to 
denote the function giving Jensen's more refined $\Ss$-hierarchy.
And $\J(\M)=\J_1(\M)$.

We take (standard) \dfnemph{premice} as in \cite{steel2010outline}, and our definition and theory 
of 
\emph{strategy premice} is modelled on \cite{steel2010outline},\cite{FSIT}. Throughout, we 
define most of the notation we use, but hopefully any unexplained terminology is either standard 
or as in those papers.
The article also uses a small part of the theory (and notation) of hod 
mice, as covered in the first parts of \cite{hod_mice}.
(However, the main scale calculations are not related particularly to hod mice,
and can be understood without knowing any theory thereof.)
For discussion of generalized solidity witnesses, see \cite{zeman}.

Our notation pertaining to iteration trees is fairly standard, but here are some points. For $\Tt$ 
a putative iteration tree, we write $\leq_\Tt$ for the tree order of $\Tt$ and $\pred^\Tt$ for the 
$\Tt$-predecessor function.
Let $\alpha+1<\lh(\Tt)$ and $\beta=\pred^\Tt(\alpha+1)$.
Then $M^{*\Tt}_{\alpha+1}$ denotes the $\N\ins M^\Tt_\beta$ such that 
$M^\Tt_{\alpha+1}=\Ult_n(\N,E)$, where $n=\deg^\Tt(\alpha+1)$ and $E=E^\Tt_\alpha$,
and $i^{*\Tt}_{\alpha+1}$ denotes $i^{\N}_{E}$, for this $\N,E$. And for 
$\alpha+1\leq_\Tt\gamma$, $i^{*\Tt}_{\alpha+1,\gamma}=i^\Tt_{\alpha+1,\gamma}\com 
i^{*\Tt}_{\alpha+1}$. Also let $M^{*\Tt}_0=M^\Tt_0$ and $i^{*\Tt}_0=\id$.
If $\lh(\Tt)=\gamma+1$ then $M^{\Tt}_\infty=M^\Tt_\gamma$, etc,
and $b^\Tt$ denotes $[0,\gamma]_\Tt$.

A premouse $\P$ is \dfnemph{$\eta$-sound} iff for every $n<\om$, if 
$\eta<\rho_n^\P$ then $\P$ is $n$-sound, and if $\rho^\P_{n+1}\leq\eta$ then letting 
$p=p_{n+1}^\P$, $p\cut\eta$ is $(n+1)$-solid for $\P$, and $\P=\Hull_{n+1}^\P(\eta\un p)$.
The \dfnemph{$\eta$-core} of $\P$ is $\cHull_{n+1}^\P(\eta\un p_{n+1}^\P)$.
Here $\Hull$ and $\cHull$ are as defined in \ref{dfn:fine_structure}.

\section{Strategy premice}
\label{StrategyPremice}

\begin{definition}
Let $Y$ be transitive. Then $\varrho_Y:Y\to\rank(Y)$ denotes the rank function.
And $\hat{Y}$ denotes $\trancl(\{(Y,\om,\varrho_Y)\})$.
For $M$ transitive, we say that $M$ is \dfnemph{rank closed} iff for every $Y\in M$,
we have $\hat{Y}\in M$ and $\hat{Y}^{<\om}\in M$. Note that if $M$ is rud closed and rank closed 
then $\rank(M)=\Ord\inter M$.
\end{definition}

\begin{definition}[$\J$-structure]
Let $\alpha\in\Ord\cut\{0\}$, let $y$ be transitive, $Y=\hat{y}$,
\[ D=\Lim\inter[\OR(Y)+\om,\OR(Y)+\om\alpha) \]
and let $\hmPvec=\left<\hmP_\beta\right>_{\beta\in D}$ be given.

We define $\J_\beta^{\hmPvec}(Y)$ for $\beta\in[1,\alpha]$, if possible, by recursion on 
$\beta$, as follows. We set
$\J_1^{\hmPvec}(Y)=\J(Y)$ and take unions at limit $\beta$.
For $\beta+1\in[2,\alpha]$, let $R=\J_\beta^{\Pvec}(Y)$ and suppose that 
$\hmPvec_{\OR(R)}=(\hmP_{0},\ldots,\hmP_{n-1})$ for some $n<\om$, and that for each $i<n$, 
$\hmP_i\sub R$ and is amenable to $R$. In this case we define
\[ 
\J_{\beta+1}^{\hmPvec}(Y)=\J(R,\hmPvec\rest R,\hmP_{0},\ldots,\hmP_{n-1}).
 \]
Note then that by induction, 
$\hmPvec\rest R\sub R$ and $\hmPvec\rest R$ is amenable to 
$R$.

For $m<\om$ let $\Ll_{\J,m}$ be the language with binary relation 
symbol 
$\dot{\in}$, predicate 
symbols $\dot{\hmPvec}$ and 
$\hmPdot_i$ for $i<m$, and constant symbol $\hmb$.

Let $m<\om$.
An \dfnemph{$m$-$\J$-structure over $Y$} is an amenable
$\Ll_{\J,m}$-structure
\[ \M=(\J_{\alpha}^\hmPvec(Y),\in^\M,\hmPvec,Y;\hmP_0,\ldots,\hmP_{m-1}), \]
where $\alpha\in\Ord\cut\{0\}$ and $\hmPvec=\left<\hmPvec_{\gamma}\right>_{\gamma\in D}$ 
with domain $D$ defined as above, 
the universe $\univ{\M}=\J_\alpha^{\hmPvec}(Y)$ is defined, 
$\dot{\in}^\M={\in\inter\univ{\M}}$, $\lh(\hmPvec_\gamma)=n$ for each $\gamma\in D$, 
$\dot{\hmPvec}^\M=\hmPvec$, 
$\hmPdot_i^\M=\hmP_i$, and $\hmb^\M=Y$.

Let $\M$ be a $m$-$\J$-structure over $Y$, and adopt the notation above. Let 
$\l(\M)$ denote $\alpha$. For 
$\beta\in[1,\alpha]$ and $R=\J_\beta^{\hmPvec}(Y)$ and $\gamma=\OR(R)$, let
\[ 
\M|\beta=(R,{\in}\inter 
R,\hmPvec\rest R,Y;\hmPvec_{\gamma,0},\ldots,\hmPvec_{\gamma,m-1}) \]
where $\hmPvec_{\OR(\M),i}=\hmP_i$.
We write $\N\ins\M$, and say that $\N$ is an \dfnemph{initial segment} of $\M$, iff $\N=\M|\beta$ 
for 
some $\beta$. Clearly if $\N\ins\M$ then 
$\N$ is 
an $m$-$\J$-structure over $Y$.
We write $\N\pins\M$, and say that $\N$ is a \dfnemph{proper segment} of $\M$, iff $\N\ins\M$ but 
$\N\neq\M$.

A \dfnemph{$\J$-structure} is an $m$-$\J$-structure, for some $m$.
\end{definition}

\begin{definition}
 A $\J$-structure $\M$ over $A$ is \dfnemph{acceptable} iff for all $\N\pins\M$ and all 
$\alpha<\OR(\N)$,
 if there is $X\sub A^{<\om}\cross\alpha^{<\om}$ such that $X\in\J(\N)\cut\N$,
 then in $\J(\N)$ there is a map $A^{<\om}\cross\alpha^{<\om}\onto\N$.
\end{definition}

The following lemma is proven just like the corresponding fact for $L$.

\begin{lemma}\label{lem:canonical_surjection}
Let $\M$ be a $\J$-structure over $A$.
Then there is a map, which we denote $h^\M$, such that
\[ h^\M:A^{<\om}\cross\l(\M)^{<\om}\onto\M\]
whose graph is $\Sigma_1^{\M\rest\L_{\J,0}}$, uniformly in $\M$.
Moreover, for $\N\ins\M$, we have $h^\N\sub h^\M$.
\end{lemma}

\begin{definition}
 Let $\M$ be an acceptable $\J$-structure over $A$ and $\rho<\OR(\M)$. Then $\rho$ is an 
\dfnemph{$A$-cardinal} of $\M$ iff $\M$
has no map $A^{<\om}\cross\gamma^{<\om}\onto\rho$ where $\gamma<\rho$.
Let $\Theta^\M$ denote the least $A$-cardinal of $\M$, if such exists.
We say that $\rho$ is \dfnemph{$A$-regular} in $\M$ iff $\M$ has no map 
$A^{<\om}\cross\gamma^{<\om}\cofto\rho$ where $\gamma<\rho$.
We say that $\rho$ is an \dfnemph{ordinal-cardinal} of $\M$ iff $\M$ has 
no map $\gamma^{<\om}\onto\rho$ where $\gamma<\rho$.
\end{definition}

\begin{lemma}
 Let $\M$ be an acceptable $\J$-structure over $A$ and $0<\xi<\l(\M)$.
 Let $\kappa$ be an $A$-cardinal of $\M$ such that $\kappa\leq\OR(\M|\xi)$.
 Then $\rank(A)<\kappa\leq\xi$ and $\kappa=\OR(\M|\kappa)$.
\end{lemma}

\begin{lemma}\label{lem:H_Theta^M_etc}
 There is a $\Sigma_1$ formula $\varphi\in\Ll_{\J,0}$ such that,
 for any acceptable $\J$-structure $\M$ over $A$, we have the following.
 
 Suppose $\Theta=\Theta^\M$ exists. Then:
 \begin{enumerate}
  \item\label{item:power(a)_unbounded} $\Theta$ is the least 
$\alpha$ such that $\pow(A^{<\om})^\M\sub\M|\alpha$.
  \item\label{item:M_Theta_is_H_Theta} $\univ{\M|\Theta}$ is the set of all $x\in\M$ such that 
$\trancl(x)$ is the surjective image of $A^{<\om}$ in $\M$.
  \item\label{item:uniform_surjections} Over $\M|\Theta$, $\varphi(0,\cdot,\cdot)$ defines 
a function $G:\Theta\to\M|\Theta$ such that for all $\alpha<\Theta$, 
we have $G(\alpha):A^{<\om}\onto\M|\alpha$.
 \item\label{item:Theta_A-regular} $\Theta$ is $A$-regular in $\M$.
\end{enumerate}

 Let $\kappa_0<\kappa_1$ be consecutive $A$-cardinals of $\M$. Then:
 \begin{enumerate}[resume*]
\item $\kappa_1$ is the least $\alpha$ such that 
$\pow(A^{<\om}\cross\kappa_0^{<\om})^\M\sub\M|\alpha$.
\item $\univ{\M|\kappa_1}$ is the set of all $x\in\M$ such that $\trancl(x)$ is the surjective 
image of $A^{<\om}\cross\kappa_0^{<\om}$ in $\M$.
\item Over $\M|\kappa_1$,
$\varphi(\kappa_0,\cdot,\cdot)$ defines a map $G:\kappa_1\to\M|\kappa_1$ such that for all 
$\alpha<\kappa_1$,
we have $G(\alpha):A^{<\om}\cross\kappa_0^{<\om}\onto\M|\alpha$.
\item $\kappa_1$ is $A$-regular in $\M$.
\end{enumerate}
\end{lemma}
\begin{proof}
We just prove parts \ref{item:power(a)_unbounded}--\ref{item:Theta_A-regular}; the others are
similar.

Let $\gamma$ be least such that $\pow(A^{<\om})\inter\M\sub\M|\gamma$.
Then $\gamma$ is a limit ordinal. By acceptability, for every $\alpha<\gamma$,
$\M|\gamma$ has a map $A^{<\om}\onto\M|\alpha$.

Now suppose that $\gamma<\Theta$, and let $g:A^{<\om}\onto\gamma^{<\om}$ be in $\M$.
Let $h=h^{\M|\gamma}$. Then because $g,h\in\M$, clearly $\M$ has a map
$f:A^{<\om}\onto\M|\gamma$, so $\M$ has a map 
$A^{<\om}\onto\pow(A^{<\om})^\M$, a contradiction.

So $\gamma=\Theta$, which gives parts \ref{item:power(a)_unbounded},\ref{item:M_Theta_is_H_Theta}.

Now consider part \ref{item:uniform_surjections}. Let $\alpha<\Theta$.
We will define
$g:A^{<\om}\cross A^{<\om}\onto\M|\alpha$,
and the uniformity in the definition will yield the result.
Let $\beta\in[\alpha,\Theta)$ be least such that
\[ \pow(A^{<\om})\inter\M|\beta\not\sub\M|\alpha.\]
Let $h=h^{\M|\beta}$. Let $x\in A^{<\om}$ be such that for some $y$, $f=h(x,y)$ is 
such that $f:A^{<\om}\to\M|\alpha$ is a surjection (such $x$ exists by acceptability). Let $y$ 
be least such,
and $f=h(x,y)$. Then for $z\in A^{<\om}$,
define $g(x,z)=f(z)$. For all other $(x,z)$, $g(x,z)=\emptyset$. This completes the definition
of $g$, and the uniformity is clear.

Part \ref{item:Theta_A-regular} now follows.
\end{proof}

\begin{corollary}
 Let $\M$ be an acceptable $\J$-structure over $A$ and let $\gamma$ be an $A$-cardinal of 
$\M$.
 If $\gamma$ is a limit of $A$-cardinals of $\M$ then 
$\M|\gamma$ satisfies Separation and Power Set.
If $\gamma$ is not a limit of $A$-cardinals of $\M$ 
then $\M|\gamma\sats\ZFmin$. In particular, $\M|\Theta^\M\sats\ZFmin$.
\end{corollary}

\begin{lemma}\label{lem:cardinals_equiv}
Let $\M$ be an acceptable $\J$-structure over $A$ such that $\Theta^\M$ exists.
Let $\kappa\in[\Theta^\M,\OR(\M))$.
Then $\kappa$ is an $A$-cardinal of $\M$ iff $\kappa$ is an ordinal-cardinal of $\M$.
\end{lemma}
\begin{proof}
 Suppose $\kappa>\Theta=\Theta^\M$ and $\kappa$ is an ordinal-cardinal,
 but $\M$ has a map 
\[ f:A^{<\om}\cross\gamma^{<\om}\onto\kappa \] where $\gamma<\kappa$.
 For each $y\in\gamma^{<\om}$, let $f_y:A^{<\om}\to\kappa$ be $f_y(x)=f(x,y)$,
 and let $g_y$ be the norm associated to $f_y$
 (that is, $f_y(x)<f_y(x')$ iff $g_y(x)<g_y(x')$, and $\rg(g_y)$ is an ordinal).
Then $g_y\in\M$ and $\rg(g_y)<\Theta$, because the prewellorder 
on 
$A^{<\om}$ determined by $f_y$ is in $\M|\Theta$ and $\M|\Theta\sats\ZFmin$. Similarly, 
the 
function $y\mapsto(f_y,g_y)$ is in $\M$.
 Let
 \[ h:\Theta\cross\gamma^{<\om}\onto\kappa \] be as follows.
 Let $(\alpha,y)\in\Theta\cross\gamma^{<\om}$. If $\alpha\notin\rg(g_y)$ then $h(\alpha,y)=0$;
 otherwise $h(\alpha,y)=f(x,y)$ where $g_y(x)=\alpha$. Then $h\in\M$, a
contradiction.
\end{proof}

\begin{definition}\label{dfn:cardinal_successors}
 Let $\M$ be an acceptable $\J$-structure over $A$ and let $\kappa<\OR(\M)$.
Then $(\kappa^+)^\M$ denotes either the least ordinal-cardinal $\gamma$ of 
$\M$ such that $\gamma>\kappa$, if there is such, and denotes $\OR(\M)$ otherwise.
By \ref{lem:cardinals_equiv}, if $\Theta^\M\leq\kappa$, then 
$(\kappa^+)^\M$ is the least $A$-cardinal $\gamma$ of $\M$ such that $\gamma>\kappa$,
if there is such, or is $\OR(\M)$ otherwise. This applies when $\hmE\neq\emptyset$ 
in \ref{dfn:potential_hpm} below.
\end{definition}

\begin{definition}\label{dfn:potential_hpm}
Let $\Ll=\Ll_{\J,2}\un\{\hmpdot,\hmPsidot\}$, where $\hmpdot,\hmPsidot$ are constant 
symbols. Let $\hmEdot=\hmPdot_0$ and $\hmPdot=\hmPdot_1$.

Let $a$ be transitive and $A=\hat{a}$. A \dfnemph{potential hybrid premouse \tu{(}hpm\tu{)} over 
$A$} 
is an amenable $\Ll$-structure
 \[ \M=(\J_\alpha^{\hmPvec}(A),\in^\M,\hmPvec,A;\hmE,\hmP;\hmp,\hmPsi) \]
where $\hmEdot^\M=E$, etc, with the following properties:
\begin{enumerate}
 \item $\bar{\M}=\M\rest\Ll_{\J,2}$ is a $2$-$\J$-structure.
 \item Either $\hmP=\emptyset$ or $\hmE=\emptyset$.
 \item If $\hmE\neq\emptyset$ then $\alpha$ is a limit and there is an extender $F$ over $\M$ such 
that:
\begin{enumerate}[label=--]
\item $\rank(A)<\hmmu=\crit(F)$,
\item $F$ is $A^{<\om}\times\gamma^{<\om}$-complete for all $\gamma < \hmmu$,
\item $\hmE$ is the
amenable code for $F$, as in \cite{steel2010outline}, and the premouse axioms \cite[Definition 
2.2.1]{wilson2012contributions} hold for 
$(\univ{\M},\hmPvec,\hmE)$.
\end{enumerate}
(It follows that $\M$ has a largest cardinal 
$\delta$, and $\delta\leq i_F(\hmmu)$, and $\OR(\M)=(\delta^+)^U$ where $U=\Ult(\M,F)$,
and $i_F(\hmPvec\rest(\hmmu^+)^\M)\rest\OR(\M)=\hmPvec$.)
\item For every $\bar{\N}\ins\bar{\M}$, 
$\N=(\bar{\N};\hmp,\hmPsi)$ is a potential hybrid premouse over $A$ (so 
$\hmp,\hmPsi\in\J(A)$).
\end{enumerate}

Let $\M$ be a potential hpm. We write $\N\ins\M$ iff $\N$ is as above.
Likewise $\N\pins\M$. For $\alpha\leq\l(\M)$, $\M|\alpha$ denotes the $\N\ins\M$ such that 
$\l(\N)=\alpha$, and $\M||\alpha$ denotes the potential hpm $\N$ which is the same as $\M|\alpha$,
except that $\hmE^\N=\emptyset$. (So $\hmP^{\M||\alpha}=\hmP^{\M|\alpha}$ always,
which will help ensure that $\hmP^{\M||\alpha}$ is the kind of structure we want to consider.)
\end{definition}

\begin{remark}
Let $\N$ be a potential hpm over $A$. Suppose $\hmE^\N$ codes an extender $F$. Clearly 
$\kappa=\crit(F)>\Theta^\M>\rank(A)$. By \cite[Definition 
2.2.1]{wilson2012contributions}, we have 
$(\kappa^+)^\M<\OR(\M)$; cf. \ref{lem:H_Theta^M_etc}. Note that
\emph{we allow $F$ to be of superstrong type} (see 
\ref{someNotations}) in accordance with 
\cite{wilson2012contributions}, not \cite[Definition 2.4]{steel2010outline}.\footnote{The main 
point of permitting superstrong 
extenders is that it simplifies 
certain things. But it complicates others. If the reader prefers, one could 
instead require that $F$ \emph{not} be superstrong, but 
various 
statements throughout the paper regarding condensation would need to be modified, along the lines 
of \cite[Lemma 3.3]{FSIT}.}
\end{remark}

\begin{remark}
 From now on we will omit ``$\in^\M$'' from the list of predicates for $\J$-structures $\M$.
\end{remark}

\begin{definition}
\label{someNotations}
Let $\M$ be a potential hpm over $A$.
We say that $\M$ is \dfnemph{$\hmE$-active} iff
$\hmE^\M\neq\emptyset$, and \dfnemph{$\hmP$-active} iff
$\hmP^\M\neq\emptyset$. \dfnemph{Active} means either
$\hmE$-active or $\hmP$-active. \dfnemph{$\hmE$-passive} means not
$\hmE$-active. \dfnemph{$\hmP$-passive} means not $\hmP$-active. \dfnemph{Passive}
means not active. \dfnemph{Type 0} means passive. 
\dfnemph{Type 4} means $\hmP$-active. \dfnemph{Type 1, 2} or \dfnemph{3} 
mean $\hmE$-active, with the usual distinctions.

We write $\actext^\M$ for the extender $F$ coded by $\hmE^\M$ (where $F=\emptyset$ if 
$\hmE^\M=\emptyset$).
We write $\es^\M$ for the function with domain $\l(\M)$, sending 
$\alpha\mapsto\actext^{\M|\alpha}$.
Likewise for $\es^\M_+$, but with domain $\l(\M)+1$.

If $F=F^\M\neq\emptyset$, we say $\M$, or $F$, is \dfnemph{superstrong} iff 
$i_F(\crit(F))=\nu(F)$. We say that $\M$ is \dfnemph{super-small} iff $\M$ has no superstrong 
whole segment.
We define $\M^\sq$ as in \cite{FSIT}. (Unless $\M$ is type $3$,
we have $\M^\sq=\M$.)
\end{definition}

\begin{definition}
Let $\Ll^-=\Ll\cut\{\hmEdot,\hmPdot\}$.
Let $\Ll^+=\Ll\un\{\hmmudot,\hmedot\}$, where $\hmmudot,\hmedot$ are constant symbols.

Let $\N$ be a potential hpm over $A$.

If $\N$ is $\hmE$-active then $\hmmu^\N\eqdef\crit(F^\N)$, and otherwise $\hmmu^\N\eqdef\emptyset$.

If $\N$ is $\hmE$-active type $2$ then $\hme^\N$ denotes the trivial completion of the largest 
non-type $Z$ proper segment of $F$; otherwise $\hme^\N\eqdef\emptyset$.\footnote{In \cite{FSIT}, 
the (analogue of) $\hme$ is referred to by its 
code $\gamma^\M$. We use $\hme$ instead because this does not depend 
on having (and selecting) a wellorder of $\M$.}

If $\N$ is not type 3 then $\core_0(\N)=\N^\sq$ denotes the $\L^+$-structure 
$(\N,\hmmu^\N,\hme^\N)$ (with $\hmmudot^\N=\hmmu^\N$ etc).

If $\N$ is type 3 then define the $\Ll^+$-structure $\core_0(\N)=\N^\sq$ essentially as in 
\cite{FSIT}; so letting 
$\hmPvec=\hmPvec^\N$ and $\nu=\nu(F^\N)$,
 \[ 
\N^\sq=(\J_\nu^{\Pvec\rest\nu}(A),\Pvec\rest\nu,A;\hmE',\emptyset;\hmp^\N,\hmPsi^\N,\hmmu^\N,
\emptyset)\]
 where $\hmE'$ is defined as usual. We also let $(\N^\sq)^{\mathrm{unsq}} = \N$.
\end{definition}

\begin{definition}
\dfnemph{$\Ll^+$-Q-formulas} and \dfnemph{$\Ll^+$-P-formulas} are 
defined analogously to in \cite[\S\S2,3]{FSIT}, using the language $\Ll^+$, but with the 
$\rSigma_1$ of \cite{FSIT} replaced by $\Sigma_1$.
\end{definition}

\begin{lemma}\label{lem:pot_hpm_type_Q-formulas} There are $\Ll^+$-Q-formulas 
$\varphi_0,\varphi_1,\varphi_2,\varphi_4$ and an $\Ll^+$-P-formula $\varphi_3$,
such that for all wellfounded $\Ll^+$-structures $\N$ with $\mu^\N\in\Ord(\N)$:
\begin{enumerate}[label=--]
\item For $i\in\{0,1,2,4\}$, $\N\sats\varphi_i$ iff $\N=\core_0(\M)$ for 
some type $i$ potential hpm $\M$.
\item If $\N=\core_0(\M)$ for 
a type 3 potential hpm $\M$ then $\N\sats\varphi_3$, and if $\N\sats\varphi_3$ then $\hmE^\N$ 
codes an 
extender 
$F$ over $\N$ and if $\Ult(\N,F)$ is wellfounded then 
$\N=\core_0(\M)$ for a type 3 potential hpm.
\end{enumerate}
\end{lemma}
\begin{proof}
This is a routine adaptation of the analogues
\cite[Lemma 2.5]{FSIT}, \cite[Lemma 3.3]{FSIT} respectively, with the added point that we can 
drop the clause ``or $\N$ is of superstrong type'' of \cite[Lemma 
3.3]{FSIT}, because we allow extenders of superstrong type.
\end{proof}
\begin{definition}\label{weak0}
Let $\N$ be a potential hpm. Let $\R$ be an $\Ll^+$-structure (possibly illfounded). Let 
$\pi:\R\to\core_0(\N)$.

We say that $\pi$ is a \dfnemph{weak $0$-embedding} iff $\pi$ is 
$\Sigma_0$-elementary (therefore $\R$ is extensional and wellfounded, so assume $\R$ is 
transitive) and there is $X\sub\R$ such that $X$ is $\in$-cofinal in $\R$ and
$\pi$ is $\Sigma_1$-elementary on elements of $X$, and if $\N$ is type 1 or 2, then letting 
$\mu=\hmmu^\R$,
there is $Y\sub\R|{(\mu^+)^\R}\cross\R$ such that $Y$ is $\in\cross\in$-cofinal in 
$\R|{(\mu^+)^\R}\cross\R$ and $\pi$ is $\Sigma_1$-elementary on 
elements of $Y$.

Let $\M,\N$ be type $i$ potential hpms. A \dfnemph{weak $0$-embedding} $\pi$ from $\M$ to $\N$,
denoted $\pi:\M\to\N$, is a weak $0$-embedding $\pi:\core_0(\M)\to\core_0(\N)$.
(So for example, if $i=3$ then $\dom(\pi)\neq\univ{\M})$.)
\end{definition}

\begin{lemma}\label{lem:Qformula_pres_weak0}
Let $\M$ be a potential hpm, let $\R$ be an $\Ll^+$-structure and let $\pi:\R\to\core_0(\M)$ be
a weak $0$-embedding.

Suppose $\M$ is type $i\neq 3$. Then $\R=\core_0(\N)$ for some type $i$ 
potential hpm $\N$. In fact, for any $\Ll^+$-Q-formula $\varphi$, if $\core_0(\M)\sats\varphi$ then 
$\R\sats\varphi$.

Suppose $\M$ is type 3. For any $\Ll^+$-P-formula $\varphi$,
if $\core_0(\M)\sats\varphi$ then $\R\sats\varphi$.
If $\Ult(\M,F^\M)$ is wellfounded then 
$\R=\core_0(\N)$ for some type 3 potential hpm $\N$.
\end{lemma}
The proof is routine, so we omit it.

\begin{definition}
Let $\R$ be an $\Ll^+$-structure.
Let $\Gamma$ be a collection of $\Ll^+$-formulas with ``$x=\dot{c}$'' in $\Gamma$ for each constant 
$\dot{c}\in\Ll^+$.
Let $X\sub\univ{\R}$. Then
\[ \Hull^\R_\Gamma(X)\eqdef(H,{\in}',\hmPvec',\hmb^\R;
\hmE',\hmP';\hmpdot^\R,\hmPsidot^\R,\hmmudot^\R,\hmedot^\R), \] where $H$ is the set of all 
$y\in\univ{\R}$ such that for some 
$\varphi\in\Gamma$ and $\xvec\in 
X^{<\om}$, $y$ is the unique $y'\in\R$ such that $\R\sats\varphi(\xvec,y')$;
and ${\in}'={\in}^\R\inter H^2$ and $\hmPvec'=\hmPvec^\R\inter H$, 
etc.
If $\R$ is transitive, then
$\C=\cHull^\R_\Gamma(X)$
denotes the $\L^+$-structure which is the transitive collapse of $\Hull^\R_\Gamma(X)$. (That is,
$\univ{\C}$ is the transitive collapse of $H$, and letting $\pi:\univ{\C}\to H$ be the uncollapse, 
$\hmE^\C=\pi^{-1}(\hmE^\R)$, etc.)
\end{definition}

\begin{definition}\label{dfn:fine_structure}
Let $\M$ be a potential hpm and $\R=\core_0(\M)$. The \dfnemph{fine structural notions} for $\M$
are just those of $\R$. We sketch the definition of the \dfnemph{fine structural notions} for $\R$.
For extra details refer to \cite{FSIT},\cite{steel2010outline};
we also adopt some simplifications explained in \cite{local_prikry}.\footnote{The simplifications 
involve dropping the parameters $u_n$, and replacing the use of generalized theories with pure 
theories. These changes are not important, and if the reader prefers, one could redefine things 
more analogously to \cite{FSIT},\cite{steel2010outline}.} Let $A=\hmb^\R$.

We say that $\R$ is \dfnemph{$0$-sound} and let $\rho_0^\R=\OR(\R)$ and 
$p_0^\R=\emptyset$ and $\core_0(\R)=\R$ and
$\rSigma_{1}^{\R}=\Sigma_1^{\R}$. (Here and in what follows, definability uses 
$\Ll^+$.) 

Now let $n<\om$ and suppose that $\R$ is $n$-sound (which will imply that 
$\R=\core_n(\R)$) and that $\om<\rho_n^{\R}$.
We write $\pvec_n^\R=(p_1^\R,\ldots,p_n^\R)$.
Then $\rho_{n+1}^\R$ is the least ordinal $\rho\geq\om$ such that for some 
$X\sub A^{<\om}\cross\rho^{<\om}$, $X$ is $\bfrSigma_{n+1}^\R$ but $X\notin\univ{\R}$.
And $p_{n+1}^\R$ is the least tuple $p\in\Ord^{<\om}$ such that some such $X$ is 
\[ \rSigma_{n+1}^\R(A\un\rho_{n+1}^\R\un\{p,\pvec^\R_n\}).\]
For any $X\sub\univ{\R}$, let
\[ \Hull_{n+1}^\R(X)=\Hull_{\rSigma_{n+1}}^\R(X),\]
and $\cHull_{n+1}^\R(X)$ be 
its transitive collapse.
Then we let
\[ \C=\core_{n+1}(\R)=\cHull_{n+1}^\R(A\un\rho_{n+1}^\R\un\pvec_{n+1}^\R), \]
and the uncollapse map $\pi:\C\to\R$ is the associated \dfnemph{core embedding}.
Define $(n+1)$-\dfnemph{solidity} and $(n+1)$-\dfnemph{universality} for $\R$ as usual (putting all 
elements of $A$ into every relevant hull).
We say that $\R$ is $(n+1)$-\dfnemph{sound} iff $\R$ is $(n+1)$-solid and $\C=\R$ and $\pi=\id$.

Now suppose that $\R$ is $(n+1)$-sound and $\rho_{n+1}^\R>\om$ (so $\rho_{n+1}^\R>\rank(A)$). 
Define $T=T_{n+1}^\R\sub\R$ by letting $t\in T$ iff 
for some $q\in\R$ and $\alpha<\rho_{n+1}^\R$,
\[ t=\Th_{\rSigma_{n+1}}^\R(A\un\alpha\un\{q\}). \]
(This denotes the \emph{pure} $\rSigma_{n+1}$ theory, as opposed to the \emph{generalized} 
$\rSigma_{n+1}$ theory of \cite{FSIT}.\footnote{As in \cite[\S2]{FSIT}, it does not matter which we 
use.}) Define $\rSigma_{n+2}^\R$ from $T$ as 
usual.
\end{definition}

\begin{definition}
Let $k\leq\om$ and let $\M,\N$ be a $k$-sound potential hpms.

A \dfnemph{\tu{(}near\tu{)} $k$-embedding} $\pi:\M\to\N$, literally a \dfnemph{\tu{(}near\tu{)} 
$k$-embedding} $\pi:\core_0(\M)\to\core_0(\N)$, is analogous to the corresponding notion in 
\cite{steel2010outline} (but the elementarity is with respect to the language $\Ll^+$ the fine 
structure is 
that of $\core_0(\M)$ and $\core_0(\N)$). If $k\geq 1$, a \dfnemph{weak $k$-embedding} 
$\pi:\M\to\N$ is likewise, but analogous to the corresponding notion in
\cite[Definition 2.1(?)]{copy_con}.\footnote{Note that 
this definition of \emph{weak $k$-embedding} diverges slightly from the definitions given in 
\cite{FSIT} and \cite{steel2010outline}.} Recall that when $k=\om$, each of these notions are 
equivalent with full elementarity.

A \dfnemph{(weakly, nearly) $k$-good} embedding $\pi:\M\to\N$ is a (weak, near) 
$k$-embedding $\pi:\M\to\N$ such that $\hmb^\M=\hmb^\N$ and $\pi\rest\hmb^\M=\id$.
\end{definition}

\begin{definition}
Let $\N$ be an $\om$-sound potential hpm. We say that $\N$ is \dfnemph{${<\om}$-condensing} (or 
satisfies \dfnemph{${<\om}$-condensation})
iff for every $k<\om$, every $(k+1)$-sound potential hpm $\M$, 
every weak $k$-embedding $\pi:\M\to\N$ such that $\rho=\rho_{k+1}^\M\leq\crit(\pi)$ and 
$\rho<\rho_{k+1}^\N$, either $\M\pins\N$ or $\M\pins\Ult(\N|\rho,F^{\N|\rho})$.
\end{definition}

\begin{definition}
 A \dfnemph{hybrid premouse \tu{(}hpm\tu{)}} is a potential hpm $\M$ such that every $\N\pins\M$ is 
$\om$-sound and ${<\om}$-condensing.
\end{definition}

\begin{lemma}
Lemmas \ref{lem:pot_hpm_type_Q-formulas} and \ref{lem:Qformula_pres_weak0}
hold with every instance of \emph{potential hpm} replaced by \emph{hpm}.
\end{lemma}

We now proceed to defining \emph{strategy premice}, or, \emph{$\Sigma$-premice}, for an iteration 
strategy $\Sigma$. We first define the process we use to feed in branches determined by 
$\Sigma$.
For $\gamma\in\Ord$ and $b\sub\Ord$, we write $\gamma+b$ for 
$\{\gamma+\alpha\st\alpha\in b\}$.
Given a structure $\M$, an iteration tree $\Tt\in\M$ of length $\om\lambda$, and a $\Tt$-cofinal 
branch $b$, Woodin noticed that $\M$ can be extended to a structure $\N$ over which $b$ is added 
with an amenable predicate, with $\N=(\J_\lambda(\M),\OR(\M)+b)$. We will use a variant of this:

\begin{definition}[$\BBB$, $b^\M$]\label{dfn:B(M,T,b)}
Let $\Q$ be an hpm over $A$ with $\N=\hmp^\Q$ transitive. Let
$\lambda>0$ and let
$\Tt$ be an iteration tree\footnote{We formally take an \emph{iteration tree} to include
the entire sequence $\left<M^\Tt_\alpha\right>_{\alpha<\lh(\Tt)}$ of models. So $\N$ is 
determined by $\Tt$, and ``$\Tt$ is an iteration tree on
$\N$'' is 
$\Sigma_0(\Tt,\N)$.} on $\N$, with $\lh(\Tt)=\om\lambda$ and $\Tt\rest\beta\in\Q$ for 
all $\beta\leq\lh(\Tt)$. Let $\zeta\in[1,\lambda]$ and $b\sub\om\zeta$ be such that 
$b\inter\beta\in\Q$ for all 
$\beta<\om\zeta$.

Then $\BBB(\Q,\Tt,\om\zeta,b)$ denotes the potential hpm $\Ss$ such that $\Q\pins\Ss$,
$\l(\Ss)=\l(\Q)+\zeta$, $\hmE^\Ss=\emptyset$,
\[ \hmP^\Ss=\{\Tt\}\cross\left(\OR(\Q)+b\right) \]
and for each $\R$ such that $\Q\pins\R\ins\Ss$,
\[ \hmP^\R=\{\Tt\}\cross\left(\OR(\Q)+[0,\gamma)_\Tt\right) \]
where $\OR(\Q)+\gamma=\OR(\R)$. (Note that $\Ss$ is amenable.)
We also write $b^\Ss=b$ and $\Tt^\Ss=\Tt$, and for $\R,\gamma$ as above, we write 
$b^\R=[0,\gamma)_\Tt$ and $\Tt^\R=\Tt$.

If $\zeta=\lambda$ then we write $\BBB(\M,\Tt,b)$ for $\BBB(\M,\Tt,\om\lambda,b)$.
\end{definition}

Our notion of \emph{$\Sigma$-premouse} $\N$ for an
iteration strategy $\Sigma$, proceeds basically as follows. For certain $\M\pins\N$, we will 
identify an iteration tree $\Tt\in\M$, via $\Sigma$, such that 
$\Sigma(\Tt)$ is not encoded into $\M$, but $\Sigma(\Tt\rest\alpha)$ is encoded into $\M$, for all 
limits $\alpha<\lh(\Tt)$. Let $\Ss=\BBB(\M,\Tt,\Sigma(\Tt))$. In a common case, then either 
$\Ss\ins\N$ or $\N\ins\Ss$. (For the kind of $\Sigma$-premouse most central to this paper,
we will actually need a generalization of this, in which there will be some $\R$ such that 
$\M\pins\R\pins\Ss$ and $\R\pins\N$, but $\N$ disagrees with $\Ss$ above $\R$.)
Clearly if $\lh(\Tt)>\om$ then $\Ss$ codes redundant information between $\M$ and $\Ss$ (the 
branches $\Sigma(\Tt\rest\alpha)$ for $\alpha<\lh(\Tt)$) before coding $\Sigma(\Tt)$ itself over 
$\Ss$. The point 
of this redundancy is that it smooths out the theory a little: it seems
to allow one to prove slightly nicer condensation properties, given that $\Sigma$ itself has nice
condensation properties, while keeping the definition of 
\emph{$\Sigma$-premouse} simple.\footnote{Difficulties that arise if one 
codes $\Sigma$ by only feeding $\Sigma(\Tt)$ itself are discussed in \ref{rem:bad_spm_definition}.} 
The key 
facts are given in \ref{StrategyCondensation} and \ref{lem:spm_con_type_3} below.

We now give some terminology relating to iteration strategies we will use in this section.
Typically the domain of an iteration strategy consists of some simply definable class of trees;
we will assume that it is $\Sigma_0$ definable.

\begin{definition}
Let $\P$ be a transitive structure and $\lambda\leq\Ord$. A \dfnemph{putative $\lambda$-iteration 
strategy for $\P$} is a function $\Sigma$ 
such that $\dom(\Sigma)$ 
is a class of iteration 
trees $\Tt$ on $\P$ of limit length $<\lambda$, and for each $\Tt\in\dom(\Sigma)$, $\Sigma(\Tt)$ 
is a $\Tt$-cofinal branch. Given such a $\Sigma$, we say that $\Sigma$ has \dfnemph{recognizable 
domain} iff there is a $\Sigma_0$ formula $\psi$ in the language of set theory such that for all 
trees $\Tt$ on $\P$, we have $\Tt\in\dom(\Sigma)$ iff $\Tt$ is via $\Sigma$ and $\lh(\Tt)<\lambda$ 
and $\psi(\Tt)$.\footnote{Since $\P=M^\Tt_0\in\trancl(\Tt)$, $\psi$ can reference $\P$, and any of 
the models of $\Tt$.}
A \dfnemph{$\lambda$-iteration strategy for $\P$} is a putative strategy $\Sigma$ 
such that every putative tree via $\Sigma$ is in fact an iteration tree. (Note here 
that $\Sigma(\Tt)$ might fail to be defined for some tree $\Tt$ via $\Sigma$.)
A \dfnemph{\tu{(}putative\tu{)} iteration strategy for $\P$} is a (putative) $\lambda$-iteration 
strategy for $\P$, for some $\lambda$.
\end{definition}

\begin{definition}
Let $\M$ be a potential hpm. Then $\J^{\hpm}(\M)$ denotes the unique potential hpm $\N$ such that 
$\M\pins\N$ 
and $\l(\N)=\l(\M)+1$ and $\hmP^\N=\emptyset$. For ordinals $\alpha$, we define 
$\J_\alpha^{\hpm}(\M)$ inductively as follows. 
\begin{enumerate}[label=--]
\item $\J^{\hpm}_0(\M) = \M$ and $\J_1^{\hpm}(\M) = \J^{\hpm}(\M)$.
\item $\J_{\beta+1}^{\hpm}(\M) = \J^{\hpm}(\J^{\hpm}_\beta(\M))$.
\item For $\lambda$ limit, $\J_\lambda^{\hpm}(\M)$ is the unique passive potential hpm $\N$ such 
that $\N=\lim_{\beta<\lambda}\J_\beta^\hpm(\M)$.
\end{enumerate}

Let $a$ be transitive and $A=\hat{a}$ and $\P,\Psi\in\J(A)$. Then $\J^\hpm(A;\P,\Psi)$
denotes the unique passive potential hpm $\N$ over $A$, with $\hmp^\N=\P$, $\hmPsi^\N=\Psi$ and
$\l(\N)=1$. 
For $\alpha\geq 0$, $\J^\hpm_{1+\alpha}(A;\P,\Psi)$ denotes $\J^\hpm_\alpha(\J^\hpm(A;\P,\Psi))$.
\end{definition}

\begin{definition}\label{dfn:aspm}
An \dfnemph{abstract strategy premouse \tu{(}aspm\tu{)}} is an  
hpm $\M$ such that $\hmp^\M$ is a transitive structure and $\Psi^\M$ is a putative strategy 
for $\hmp^\M$ and there
is $\chi\in\Ord$ and sequences 
$\etavec=\left<\eta_\alpha\right>_{\alpha\leq\chi}$
and $\Sigmavec=\left<\Sigma_\alpha\right>_{\alpha\leq\chi}$
such that $\etavec$ is strictly increasing and continuous, 
$\eta=1$, $\eta_\chi=\l(\M)$, $\Sigmavec$ is an increasing (possibly not strictly) and 
continuous sequence of putative strategies for $\cp^\M$, $\Sigma_0=\Psi^\M$, and for each 
$\alpha<\chi$, either:
 \begin{enumerate}[label=--]
  \item $\M|\eta_{\alpha+1}=\J^\hpm(\M|\eta_\alpha)$ and $\Sigma_{\alpha+1}=\Sigma_\alpha$; or
 \item There is $\Tt\in\M|\eta_\alpha$ such that the following holds.
We have that $\Tt$ is an iteration tree 
via $\Sigma_\alpha$, 
but no proper extension of $\Tt$ is via $\Sigma_\alpha$. Let $\N=\M|\eta_{\alpha}$ and 
$\N'=\M|\eta_{\alpha+1}$ and 
$\theta=\lh(\Tt)$. Then there 
is $b\sub\theta$ such that $\Ss\eqdef\BBB(\N,\Tt,b)$ is defined\footnote{That is, 
$b\inter\beta\in\N$ for all $\beta<\theta$. Note that possibly
$b=\emptyset$ and $\N\pins\Ss$ here. So in this case, $\M$ is still considered a 
$\varphi$-indexed spm, even if there is no $\Tt$-cofinal branch.} and either:
  \begin{enumerate}[label=--]
  \item $\N'=\Ss$, $b$ is a $\Tt$-cofinal branch\footnote{Note that $M^\Tt_b$ might be 
illfounded. But in this 
case $\Tt\conc b$ is not an iteration tree, so there is no $\alpha\leq\theta$ such 
that $\Tt'=\Tt^{\M|\eta_\alpha}_\varphi$ is defined and $\Tt$ is properly extended by $\Tt'$.}
and $\Sigma_{\alpha+1}=\Sigma_\alpha\un\{(\Tt,b)\}$, or
  \item $\N'\pins\Ss$ and $\Sigma_{\alpha+1}=\Sigma_\alpha$.
  \end{enumerate}
\end{enumerate}
Given an aspm $\M$, we write $\chi^\M=\chi$, $\etavec^\M=\etavec$, etc, and 
$\Sigma^\M=\Sigma^\M_\chi$.\footnote{No particular 
demand is made on $\dom(\Sigma^\M)$ (though it is closed under 
initial segment).} We say that $\M$ is a 
\dfnemph{successor} iff $\chi$ is a successor. If $\M$ is a 
successor then $\M^-$ denotes $\M|\eta_{\chi-1}$.
\end{definition}

It is easy to see that the sequences $\etavec,\Sigmavec$ above are unique,\footnote{Adopt the 
notation of \ref{dfn:aspm} and let $\alpha<\chi$.
Then $\eta_{\alpha+1}$ is the least $\eta>\eta_\alpha$ such that either $\eta=\l(\M)$ or
$\hmP^{\M|(\eta+1)}=\emptyset$ or
$\hmP^{\M|(\eta+1)}=\{\Uu\}\cross B$ for some $\Uu,B$ such that 
$B\inter\OR(\M|\eta)=\emptyset$. (This is because $0\in b$ whenever $b$ is a branch through an 
iteration tree.)} so the notation 
$\etavec^\M$, etc, is unambiguous.
We select the trees $\Tt$ for which we add $\Sigma(\Tt)$ in a first-order manner:

\begin{definition}\label{dfn:J_premice}
 Let $\varphi\in\Ll^+$, $\M$ be an hpm and $\Tt\in\M$. We
write $\Tt=\Tt^\M_\varphi$ iff $\hmp^\M$ is transitive and $\Tt$ is a limit length iteration tree 
on $\hmp^\M$
and $\Tt$ is the unique 
$x\in\M$ such that $\M\sats\varphi(x)$.\qedhere
\end{definition}

The generality of the indexing device $\varphi$
in the definition below was
probably influenced by Sargsyan's \cite[Definitions 1.1, 1.2]{hod_mice}.

\begin{definition}\label{dfn:spm}
Let $\varphi\in\Ll^+$. A \dfnemph{$\varphi$-indexed strategy premouse \tu{(}$\varphi$-spm\tu{)}} is 
an 
aspm $\M$ such that letting $\etavec=\etavec^\M$, etc, for every $\alpha<\chi$,
letting $\N=\M|\eta_\alpha$ and $\N'=\M|\eta_{\alpha+1}$, we have:
 \begin{enumerate}[label=--]
  \item If $\Tt^{\N}_\varphi$ is undefined then $\hmP^{\N'}=\emptyset$ (so
$\N'=\J^\hpm(\N)$).
 \item Suppose $\Tt\eqdef\Tt^{\N}_\varphi$ is defined.
Then $\hmP^{\N'}\neq\emptyset$ and $\Tt^{\N'}=\Tt$ (so $\Tt$ is the witness to the corresponding 
clause of \ref{dfn:aspm}) and
$\Tt^\R_\varphi=\Tt$ for all $\R$ such that $\N\ins\R\pins\N'$, but if $\N'\pins\M$ then 
$\Tt^{\N'}_\varphi\neq\Tt$.
\end{enumerate}
 
Let $\M$ be a $\varphi$-spm, and let $\etavec$, etc, be as above. We say that $\M$ is 
\dfnemph{$\varphi$-whole} iff, if $\M$ is a successor and $\Tt\eqdef\Tt^{\M^-}_\varphi$ is 
defined,
then either $\M=\BBB(\M^-,\Tt,b)$ for some $b$, or $\Tt^\M_\varphi\neq\Tt$. 

Let $\Sigma$ be a putative iteration strategy for a transitive structure $\P$. Let 
$\varphi\in\Ll^+$.
A \dfnemph{$\varphi$-indexed $\Sigma$-premouse \tu{(}$(\Sigma,\varphi)$-premouse\tu{)}}, is a 
$\varphi$-spm 
$\M$ such that $\hmp^\M=\P$ and $\Sigma^\M\sub\Sigma$.
\end{definition}

Clearly if $\M$ is a $\varphi$-spm then $\Sigma^\M$ is the least putative strategy $\Sigma$ such 
that $\M$ is a
$\varphi$-indexed $\Sigma$-pm.

It seems difficult to express $\varphi$-indexed spm-hood with Q-formulas. So we 
consider the more general notion of \emph{$\varphi$-indexed possible-spm}, which we can express 
with 
Q-formulas, modulo the usual restrictions.

\begin{definition}
 A \dfnemph{$\varphi$-indexed possible spm} is an hpm $\M$ such that there is a $\varphi$-indexed 
spm 
$\N$ such that either $\M=\N$, or $\N$ is a successor, $\N^-\pins\M$,
$\Tt\eqdef\Tt^{\N^-}_\varphi$ is defined, and letting $\OR(\N)=\OR(\N^-)+\zeta$, there is a 
$\Tt\rest\zeta$-cofinal branch $b$ 
such that $\M=\BBB(\N^-,\Tt,\zeta,b)$.

We adapt terminology and notation for spms to possible-spms in the obvious manner.
\end{definition}

So a $\varphi$-indexed possible spm only fails to be a $\varphi$-spm if, with notation as above, we 
have $\zeta<\lh(\Tt)$ but $b\neq[0,\zeta)_\Tt$. The following lemma is straightforward:

\begin{lemma}\label{lem:spm_type_Q-formulas} Let $\varphi\in\Ll^+$. Then Lemma 
\ref{lem:pot_hpm_type_Q-formulas} holds with every instance of \emph{potential hpm} replaced by 
\emph{$\varphi$-indexed possible spm}.\end{lemma}

\begin{definition}\label{dfn:very_condensing}
 Let $\R,\M$ be $\hmE$-passive possible-spms and $\pi:\R\partialto\M$.
Then $\pi$ is a \dfnemph{very weak 
$0$-embedding} iff $\pi$ is $\Sigma_0$-elementary on its domain and there is an ${\in}$-cofinal set
$X\sub\R$ such that
\[X\un\OR(\R)\un\hmp^\R\un\{\hmp^\R,\hmPsi^\R,\hmb^\R\}\sub\dom(\pi), \]
$\pi\rest\hmp^\R=\id$, and $\pi$ is $\Sigma_1$-elementary on parameters in $X$.

Let $C$ be a class of possible-spms. We say that $C$ is \dfnemph{very 
condensing} iff for all $\hmE$-passive $\M\in C$ and all 
$\hmE$-passive possible-spms $\R$, if there is a very weak 
$0$-embedding $\pi:\R\to\M$ then $\R\in C$.
\end{definition}

\begin{lemma}\label{lem:very_weak_pres_Q-formula}
The truth of $\Ll^+$-Q-formulas is preserved downward under very weak 
$0$-embeddings.
\end{lemma}

We next consider preservation of $\Sigma$-pms, for strategies $\Sigma$ with hull condensation 
(see \cite[Definitions 1.29--1.31]{hod_mice}).

\begin{lemma}\label{StrategyCondensation}
Let $\M$ be a
$\varphi$-indexed spm, not of type 3. Let $\R$ be a $\varphi$-indexed possible spm.

\begin{enumerate}[label=\tu{(}\arabic*\tu{)}]
\item\label{item:piRtoM} Let $\Sigma$ be an iteration strategy with hull condensation. Suppose that 
$\M$ 
is a $\Sigma$-pm, $\hmp^\R=\hmp^\M$ and either \tu{(}i\tu{)} $\hmPsi^\R\sub\Sigma$ and there is 
a very weak 
$0$-embedding $\pi:\R\partialto\M$, or \tu{(}ii\tu{)} there is a weak 
$0$-embedding $\pi:\R\to\M$ above $\hmp^\R$.
Then $\R$ is a $\Sigma$-pm.
\item\label{item:piMtoRgen} Suppose there is $\pi\maps\M\to\R$ such that either:
\begin{enumerate}[label=\tu{(}\alph*\tu{)}]
\item\label{item:Sigma_2} $\pi$ is $\Sigma_2$-elementary, or
\item\label{item:cof_Sigma_1_Tt_undef} $\pi$ is cofinal $\Sigma_1$-elementary and either $\M$ is a 
limit or $\Tt^{\M^-}_\varphi$ is undefined, or
\item\label{item:cof_Sigma_1_Tt_def} $\pi$ is cofinal $\Sigma_1$-elementary, $\M$ is a 
successor and $\Tt=\Tt^{\M^-}_\varphi$ is defined and either $b^\M\in\M$ or $\pi$ is 
continuous at $\lh(\Tt)$.\footnote{Cf. \ref{rem:tree_cof_meas}.}
\end{enumerate}
Then $\R$ is a $\varphi$-indexed spm.
\end{enumerate}
\end{lemma}
\begin{proof}
Part \ref{item:piRtoM}: We just consider the case (i).
(So by \ref{dfn:very_condensing}, $\R,\M$ are $\hmE$-passive and $\pi$ is above $\hmp^\R$.)
We may assume that $\R$ is a successor and 
every proper segment of $\R$ is a $\Sigma$-pm, since $\pi$ induces very weak $0$-embeddings (in 
fact, fully elementary on their domains) from the proper segments of $\R$ to proper segments of 
$\M$.
It follows that $\M$ is a successor and $\pi(\R^-)=\M^-$. We may assume that
$\Ttbar=\Tt^{\R^-}_\varphi$ is defined, so $\pi(\Ttbar)=\Tt=\Tt^{\M^-}_\varphi$ is defined.
Let $\OR(\R^-)+\gammabar=\OR(\R)$ and $\OR(\M^-)+\gamma=\OR(\M)$.
Then $\pi$ induces a hull embedding from $(\Ttbar\rest\gammabar)\conc b^\R$ to 
$(\Tt\rest\gamma)\conc b^\M$. Since the latter is via $\Sigma$, as is $\Ttbar$,
hull condensation gives that $b^\R=\Sigma(\Ttbar\rest\gammabar)$, so $\R$ 
is a $\Sigma$-pm.

We leave \ref{item:piMtoRgen}\ref{item:Sigma_2} and 
\ref{item:piMtoRgen}\ref{item:cof_Sigma_1_Tt_undef} to the reader. 
Consider \ref{item:piMtoRgen}\ref{item:cof_Sigma_1_Tt_def}. Note that $\pi(\M^-)=\R^-$,
and since $\Tt=\Tt^{\M^-}_\varphi$ is defined, so is $\pi(\Tt)=\Tt^{\R^-}_\varphi$. Let 
$\OR(\M^-)+\gamma=\OR(\M)$, so 
$\OR(\R^-)+\gamma'=\OR(\R)$, where $\gamma'=\sup\pi``\gamma$.
Then $b^\M$ is $\Tt\rest\gamma$-cofinal, and since $\pi``b^\M\sub b^\R$, $b^\R$ is 
$\pi(\Tt)\rest\gamma'$-cofinal. So we may assume that $\gamma'<\lh(\pi(\Tt))$,
and must see that $b^\R=[0,\gamma')_{\pi(\Tt)}$.

Suppose $b^\M\in\M$. Then because $\pi$ 
is $\Sigma_1$-elementary, $b^\R=\pi(b^\M)\inter\gamma'$. If $\gamma'<\pi(\gamma)$ then since 
$\pi(b^\M)$ is $\pi(\Tt)\rest\pi(\gamma)$-cofinal, we are done. If $\gamma'=\pi(\gamma)$ then 
since $\gamma<\lh(\Tt)$, so $b^\M=[0,\gamma)_\Tt$, so $b^\R=\pi(b^\M)$ and we are done.

Now suppose 
that $b^\M\notin\M$ and $\pi$ is continuous at $\lh(\Tt)$. Then $\gamma=\lh(\Tt)$ and 
$\gamma'=\lh(\pi(\Tt))$, contradiction. \end{proof}

\begin{corollary}\label{cor:Sigma-pm_very_condensing}
For any strategy $\Sigma$ with hull condensation and any $\varphi\in\Ll^+$, the class of 
$\varphi$-indexed $\Sigma$-pms $\M$ such that $\hmPsi^\M=\emptyset$ is very condensing.
\end{corollary}

A type 3 analogue of \ref{StrategyCondensation} follows easily from 
\ref{StrategyCondensation}:

\begin{lemma}\label{lem:spm_con_type_3}
Let $\M$ be a type 3 $\varphi$-indexed spm.
Let $\R$ be an $\Ll^+$-structure with $\hmp^\R=\hmp^\M$.
\begin{enumerate}[label=--]
\item\label{item:piRtoM_tp3} Let $\Sigma$ be an iteration strategy with hull condensation. Let 
$\kappa=\hmmu^\M$ and suppose 
$U^\M=\Ult(\M|(\kappa^+)^\M,F^\M)$ is a $\Sigma$-pm. Let
$\pi:\R\to\core_0(\M)$ be a 
weak $0$-embedding with $\pi\rest\hmp^\R=\id$. Let 
$\mu=\hmmu^\R$ and $U^\R=\Ult(\R|(\mu^+)^\R,F^\R)$. Suppose there is an elementary
$\pi':U^\R\to U^\M$
with $\pi\sub\pi'$.

Then $\R=\Q^\sq$ for some type 3 $\varphi$-indexed $\Sigma$-pm $\Q$,
and $U^\R$ is also a $\varphi$-indexed $\Sigma$-pm.

\item\label{item:piMtoRgen_tp3} Suppose there is $\pi\maps\core_0(\M)\to\R$ such 
that either
\tu{(}i\tu{)} $\pi$ is $\Sigma_2$-elementary, or
\tu{(}ii\tu{)} $\pi$ is cofinal and $\Sigma_1$-elementary.
Let $\mu=\hmmu^\R$ and suppose that $U^\R$ \tu{(}as above\tu{)} is wellfounded.

Then $\R=\Q^\sq$ for some type 3, $\varphi$-indexed spm.
\end{enumerate}
\end{lemma}

We now define $\Sigma$-iterability for $\Sigma$-premice $\M$.
The main point is that the iteration strategy should produce iterates which are $\Sigma$-premice.
One needs to be a little careful, however, because the iterates might contain iteration trees 
outside of the domain of $\Sigma$.

\begin{definition}\label{dfn:F-iterability}
Let $\Sigma$ be an iteration strategy, $\varphi\in\Ll^+$ and $X=(\Sigma,\varphi)$. Let 
$\M$ be a $X$-pm. A \dfnemph{putative $X$-iteration tree $\Tt$ on $\M$} is defined as usual,
with the added requirement that $M^\Tt_\alpha$ is an $X$-pm for each $\alpha+1<\lh(\Tt)$ (and for 
each such $\alpha$, $E^\Tt_\alpha\in\es_+(M^\Tt_\alpha)$). Let $\Tt$ be a putative $X$-tree on $\M$.
We say that $\Tt$ is a \dfnemph{well-putative $X$-iteration tree} iff $\Tt$ is a 
the models of $\Tt$ are
all wellfounded. We say that $\Tt$ is an \dfnemph{$X$-iteration tree} 
iff $M^\Tt_\alpha$ is an $X$-pm for all $\alpha+1\leq\lh(\Tt)$.

Let $k<\om$ and let $\M\in\opbk$ be a $k$-sound $X$-pm.
Let $\theta\in\Ord$.
The iteration game 
$\G^{X,\M}(k,\theta)$ has the rules of $\G^\M(k,\theta)$, except for the following 
differences. Let $\Tt$ be the putative tree being produced.
For $\alpha+1\leq\theta$, if both players meet their requirements at all stages $<\alpha$, then,
in stage $\alpha$, player $\plII$ must first ensure that $\Tt\rest\alpha+1$ is a well-putative 
$X$-tree, and if $\alpha+1<\theta$, that $\Tt\rest\alpha+1$ is an 
$X$-tree. Given this, if $\alpha+1<\theta$, player $\plI$ then selects 
$E^\Tt_\alpha$, but we replace that requirement that $\lh(E^\Tt_\beta)<\lh(E^\Tt_\alpha)$ for all 
$\beta<\alpha$, with the requirement that $\lh(E^\Tt_\beta)\leq\lh(E^\Tt_\alpha)$ for all 
$\beta<\alpha$.\footnote{
Thus, if we reach a putative tree $\Tt$ of length $\theta$, then
$\plII$ wins iff either $\theta$ is a limit or 
$M^\Tt_{\theta-1}$ is wellfounded. If $\theta=\alpha+1$, we cannot in general expect 
$M^\Tt_{\alpha}$ to be an $X$-pm. For example, suppose that $\theta=\om_1+1$ 
and $\Sigma$ is an $(\om_1+1)$-strategy for some $\P\in\HC$. 
Then $M^\Tt_{\om_1}$ could have $\varphi$-whole successor proper segments $\N$ such 
that $\Uu=\Tt^{\N^-}_\varphi$ is defined, but $\lh(\Uu)>\om_1+1$. In this case 
$\Uu\notin\dom(\Sigma)$, so $\N$ is not an $X$-pm.
In applications such as comparison, in this circumstance we only need to know that 
$M^\Tt_{\om_1}$ is wellfounded. So we still decide the game in favour of player $\plII$ in this 
situation.}

Let $\alpha,\theta\in\Ord$.
The iteration game 
$\G^{X,\M}(k,\alpha,\theta)$ is defined just as $\G^\M(k,\alpha,\theta)$,
with the differences that (i) the rounds are runs of $\G^{X,\Q}(q,\theta)$ for some 
$\Q,q$,\footnote{Recall that (considering the rules of $G^\M(k,\alpha,\theta)$) if a round of 
$\G^{X,\M}(k,\alpha,\theta)$ reaches a tree of length $\theta$, then the game finishes at that 
point. So $\Q$ here will certainly be an $X$-pm.}and 
(ii) if $\alpha$ is a limit and neither player breaks 
any rule, and $\Ttvec$ is the sequence of trees played, then player $\plII$ wins iff
$M^{\Ttvec}_\infty$ is defined (that is, the trees eventually do not drop on their main branches, 
etc) and wellfounded.

The game $\G^{X,\M}_{\max}(k,\alpha,\theta)$ is like $\G^{X,\M}(k,\alpha,\theta)$, except that 
player $\plI$ may not drop in model or degree between rounds. (For example, in both games, 
after the first round has produced a successor length $k$-maximal 
tree $\Tt_0$, the second round forms a $q$-maximal tree $\Tt_1$ on $\Q$, for a certain $(\Q,q)$. In 
$\G^{X,\M}_{\max}$, $\Q=M^{\Tt_0}_\infty$ and $q=\deg^{\Tt_0}(\infty)$, whereas in $\G^{X,\M}$, 
player $\plI$
chooses $\Q\ins M^{\Tt_0}_\infty$ and $q\leq\om$, with $q\leq\deg^{\Tt_0}(\infty)$ if 
$\Q=M^{\Tt_0}_\infty$. Likewise at the start of every later round.)

If $\alpha$ is a limit ordinal, the game $\G^{X,\M}(k,{<\alpha},\theta)$ is like 
$\G^{X,\M}(k,\alpha,\theta)$, except that if the game runs through $\alpha$ rounds with no player 
breaking any rules within those rounds, then player $\plII$ wins automatically, irrespective of 
whether 
the direct limit model is defined or wellfounded. Likewise $\G^{X,\M}_{\max}(k,{<\alpha},\theta)$.

Now \dfnemph{$X$-$(k,\theta)$-iteration strategy}, \dfnemph{$X$-$(k,\alpha,\theta)$-maximal 
iterability}, etc, are defined from these games in the obvious manner.

The game $\G^{X,\M}_\hodlower(k,\alpha,\theta)$ is just like $\G^{X,\M}(k,\alpha,\theta)$,
except that if at the end of round $\beta$ a successor length normal tree $\Tt_\beta$ has been 
produced, and both players have met all their obligations up to that point,
and $b^{\Tt_\beta}$ drops in model or degree, then player $\plII$ wins.
\dfnemph{Hod $X$-$(k,\alpha,\theta)$-iteration strategy} and \dfnemph{-iterability} are defined 
using $\G_\hodlower^{X,\M}(k,\alpha,\theta)$.
\end{definition}
\begin{remark}
The requirement, in $\G^{\M}(k,\theta)$, that $\lh(E^\Tt_\beta)\leq\lh(E^\Tt_\alpha)$ for 
$\beta<\alpha$, is weaker than requiring $\lh(E^\Tt_\beta)<\lh(E^\Tt_\alpha)$, 
because of superstrongs. See \cite[Remark 2.44(?)]{operators} regarding this 
and changes to the comparison algorithm that are needed to accommodate superstrongs.
\end{remark}

\begin{remark}\label{rem:tree_cof_meas}
Lemma \ref{StrategyCondensation} left open the possibility that 
$\R$ fails to be a $\varphi$-indexed spm, when $\pi:\M\to\R$ is cofinal and 
$\Sigma_1$-elementary, $\M$ is a successor, $\Tt=\Tt^{\M^-}_\varphi$ is defined, 
$b^\M\notin\M$ and $\pi$ is 
discontinuous at $\lambda=\lh(\Tt)$, so $\M$ is $\varphi$-whole, 
$\lambda'=\sup\pi``\lambda<\lh(\Tt')$ where 
$\Tt'=\pi(\Tt)=\Tt^{\R^-}_\varphi$, and $b^\R\neq[0,\lambda')_{\Tt'}$.
Now let $X=(\Sigma,\varphi)$, where $\Sigma$ has hull condensation and $\varphi\in\Ll^+$, and 
suppose further that $\M$ is a $X$-iterable $X$-pm, as witnessed by some strategy 
$\Lambda$. We describe two standard circumstances below which will then lead to contradiction.

First, suppose that $\pi:\M\to\R$ is via $\Lambda$. Then because $\Lambda$ is a $X$-iteration 
strategy, $b^\R=[0,\lambda)_{\Tt'}$, a contradiction. 

Second, suppose that $\Sigma$ has hull condensation,
$\pi$ is any degree $0$ iteration embedding of $\M$ ($\pi$ need not be via any
iteration strategy). We will show that $b^\M\in\M$, for a contradiction.

Because $\pi$ is a degree $0$ iteration embedding, the discontinuity implies 
that
$\M\sats$``There is $E\in\es$ which is a total measure and $\lh(\Tt^\M)$ has
cofinality $\kappa=\crit(E)$''.
Let $C\in\M$, $C\sub\lh(\Tt)$ be a club of ordertype $\kappa$.
Then
\[ \sigma=i^\M_E:\M\to U=\Ult_0(\M,E) \]
is continuous at all points of $C$. Let $\zeta=\sup
\sigma``\lh(\Tt)$. Then $\sigma``C=\sigma(C)\inter\zeta$ is club in $\zeta$. But
\[ U\sats\text{``}\zeta<\lh(\sigma(\Tt))\text{ and }\cof(\zeta)=\kappa\text{ is
uncountable''.}\]
So $[0,\zeta)_{\sigma(\Tt)}\inter \sigma``C$ is club in $\zeta$,
and $C'\in\M$ where $C'$ is the club
\[ C'=C\inter\sigma^{-1}``[0,\zeta)_{\sigma(\Tt)}.\]
Because $\M$ is $X$-iterable, $\sigma(\Tt)$ is via $\Sigma$. But then by hull condensation, 
$\Sigma(\Tt)$ is the downward $\leq_\Tt$-closure of $C'$, which is in $\M$.
\end{remark}
\begin{definition}
Let $\M$ be an hpm and $\N\ins\M$. We say that $\N$ is a \dfnemph{cutpoint} of $\M$ 
iff for all $\P\ins\M$, if $\N\pins\P$ and $F^\P\neq\emptyset$ then
$\OR(\N)\leq\crit(F^\P)$.
And $\N$ is a \dfnemph{strong cutpoint} of $\M$ iff likewise, but with the conclusion 
replaced with ``$\OR(\N)<\crit(F^\P)$''.
\end{definition}

\begin{definition}[$\Lp^{(\Sigma,\varphi)}$]\label{dfn:Lp^X}
Let $\Sigma$ be a strategy with hull condensation for a transitive structure $\P\in\HC$, 
$\varphi\in\Ll^+$ and $X=(\Sigma,\varphi)$. Let $a$ be transitive and $A=\hat{a}$, with 
$\P\in\J(A)$. Assume $\DC_A$.

Let $n\leq\om$ and let $\M$ be an $n$-sound $X$-pm over $A$ (and 
$\eta\leq\OR(\M)$). We say that 
$\M$ is \dfnemph{countably (above-$\eta$) $X$-$(n,\om_1+1)$-iterable} iff for
every countable hpm $\Mbar$, if $\P=\hmp^{\Mbar}$ and there is an elementary $\pi:\Mbar\to\M$
then $\Mbar$ is (above-$\bar{\eta}$) $X$-$(n,\om_1+1)$-iterable
(where $\bar{\eta}$ is 
the collapse of $\eta$).
 
$\Lp^X(a)$ denotes the stack of all 
countably $X$-$(\om,\om_1+1)$-iterable $X$-premice $\M$ over $A$ such that $\M$ is 
fully sound 
and projects to $\om$.\footnote{$\DC_A$ is enough to prove that this is a stack.
For let $\M,\N$ be such $X$-premice. Because $\M,\N$ are generated by ordinals and 
elements of $A$, by taking elementary substructures which do not collapse $A$,
we may assume that there are maps $A^{<\om}\onto\M$ and $A^{<\om}\onto\N$. But then 
by $\DC_A$, we may assume that $A,\M,\N$ are countable, so we can compare $\M,\N$ 
as usual.}
Assuming $\DC_\RR$, and letting $B\sub\HC$, $\Lp^X(\RR,B)$ 
denotes $\Lp^X((\HC,B))$, and $\Lp^X(\RR)$ denotes $\Lp^X(\HC)$.\footnote{Since $\RR$ 
is not transitive, this is not an abuse of notation.}

Let $\N$ be an $X$-premouse. Then 
$\Lp^X_+(\N)$ 
denotes the stack 
of all $X$-premice $\M$ such that either $\M=\N$, or $\N\pins\M$, $\N$ is a strong cutpoint 
of 
$\M$,
$\M$ is $\OR(\N)$-sound, and there is $n<\om$ such that $\rho_{n+1}^\M\leq\OR(\N)<\rho_n^\M$ 
and $\M$ is countably above-$\OR(\N)$ 
$X$-$(n,\om_1+1)$-iterable. Note that $\Lp_+^X(\N)$ might have a largest element, which 
projects strictly across $\OR(\N)$ and is not $\om$-sound.
\end{definition}

\begin{definition}
 Let $\Sigma$ be an iteration strategy, $\varphi\in\Ll^+$, $X=(\Sigma,\varphi)$ and $\M$ be an 
$X$-pm. Let $k\leq\om$. Then $\M$ is \dfnemph{$X$-$k$-fine} iff 
for each $j\leq k$, we have (i) $\core_j(\M)$ is a $j$-solid $X$-pm, (ii) if $j<k$ then 
$\core_j(\N)$ is $(j+1)$-universal, and (iii) if $k=\om$ then $\core_\om(\N)$ is 
${<\om}$-condensing.
\end{definition}

\begin{lemma}\label{lem:fine_structure}
Let $\Sigma,\P,\varphi,X,a,A$ be as in \ref{dfn:Lp^X} \tu{(}so we assume
$\DC_A$\tu{)}.
Then:
\begin{enumerate}[label=--]
\item For $k<\om$, every $k$-sound, countably
$X$-$(k,\om_1,\om_1+1)$-iterable $X$-pm $\M$ over $A$ is $X$-$(k+1)$-fine.
 \item Every $\om$-sound, countably $X$-$(\om,\om_1,\om_1+1)$-iterable 
$X$-pm over $A$ is ${<\om}$-condensing.
\item Every countably $X$-$(0,\om_1,\om_1+1)$-iterable 
$X$-pseudo-premouse over $A$ is an $X$-pm.
\item There is no countably $X$-$(0,\om_1+1)$-iterable 
$X$-bicephalus over $A$.
\end{enumerate}
\end{lemma}
\begin{proof}
Consider for example the proof that $\M$ is $X$-$(k+1)$-fine. We may assume that $\M$ is countable, 
by $\DC_A$. If $\AC$ holds (recall that our background theory is $\ZF$) then using the condensation 
lemmas \ref{StrategyCondensation} and \ref{lem:spm_con_type_3}, it is 
straightforward to see that the proofs of the copying construction, weak 
Dodd-Jensen\footnote{$\DC_\RR$ seems to be 
used in the construction of an iteration strategy with the weak Dodd Jensen property.} and the 
fundamental fine structural theorems go through.
But we may assume $\ZFC$,
because letting $x\in\RR$ code $\M$ and $\Lambda$ be an iteration strategy for $\M$ as 
hypothesized, then we can pass to $W=L^{\Lambda,\Sigma}[x]$ (where we feed $\Lambda,\Sigma$ into 
$W$ like with strategy mice; we do not care about fine structure for $W$), replacing $\Sigma$ with 
$\Sigma'=\Sigma\inter W$.
\end{proof}

We will build $\Sigma$-mice by background construction:

 \begin{definition}\label{dfn:Fop_con}
Let $a$ be transitive and $A=\hat{a}$. Let $\Sigma$ be an iteration strategy for a transitive 
structure in $\J(A)$, let $\varphi\in\Ll^+$ and let $X=(\Sigma,\varphi)$. An
\dfnemph{$L^X[\es,A]$-construction (of length $\chi$)} is a 
sequence $\CC=\left<\N_\alpha\right>_{\alpha<\chi}$ such that for all $\alpha<\chi$:
\begin{enumerate}[label=--]
\item $\N_\alpha$ is a $X$-pm over $A$ and $\l(\N_0)=1$.
\item If $\alpha$ is a limit then $\N_\alpha=\liminf_{\beta<\alpha}\N_\beta$.
\item If $\alpha+1<\chi$ then letting $\N=\N_{\alpha+1}$, either:
\begin{enumerate}[label=--]
\item $\N$ is $\hmE$-active and
$\N||\OR(\N)=\N_\alpha$ and letting $\kappa=\hmmu^\N$, then $\Ult(\N|(\kappa^+)^\N,F^\N)$ is an 
$X$-pm, or
\item $\N_\alpha$ is $X$-$\om$-fine and 
$\M\eqdef\core_\om(\N_\alpha)\pins\N$ and 
$\l(\N)=\l(\M)+1$.\qedhere
\end{enumerate}
\end{enumerate}
\end{definition}

We will consider fully backgrounded $L^\Sigma[\es,A]$-constructions. Assume $\DC_A$. Then 
given 
$\N_\alpha$ and supposing that $\N_\alpha$ is $X$-$k$-fine,
countable $X$-$(k,\om_1,\om_1+1)$-iterability will be enough to verify that $\N_\alpha$ is 
$X$-$(k+1)$-fine. This iterability will be established (where we can) by the standard arguments, 
using the condensation lemmas.

\begin{remark}\label{rem:bad_spm_definition} Our definition of $\Sigma$-premice (for an iteration 
strategy $\Sigma$) differs a little from the standard one. The standard one is along
the lines of: given $\M|\alpha$, letting $\Tt\in\M|\alpha$ be the
$<_{\M|\alpha}$-least tree for which $\M|\alpha$ does not know $\Sigma(\Tt)$,
and $\om\lambda=\lh(\Tt)$, let $\M|(\alpha+\lambda)=(\J_\lambda(\M|\alpha),B)$,
such that $B$ codes $\Sigma(\Tt)$ amenably.

Whatever one's definition of $\Sigma$-premice, one would probably like to know that an ultrapower 
of a $\Sigma$-premouse is also a $\Sigma$-premouse.
As has been observed by others, this is not true of the hierarchy described
above. For suppose $\M|\alpha$, $\Tt$ and $\lambda$ are as above, and $\lh(\Tt)$
has measurable cofinality $\kappa$ in $\M|(\alpha+\lambda)$, and $E$ is an
extender over $\M=\M|(\alpha+\lambda)$ with $\crit(E)=\kappa$. Then
$U=\Ult_0(\M,E)$ is not in the hierarchy. For $i_E$ is discontinuous at
$\lh(E)$, but $\OR(U)=\sup i_E``\OR(\M)$.

There seem to have been two approaches used to correct this problem (other than the one we use)
used by others. 
One is to feed
in all initial segments of $\Sigma(\Tt)$ (even though they have been fed in
earlier), immediately prior to feeding in $\Sigma(\Tt)$ itself. But this
approach seems flawed. For $(*)$ let $\M'$ be a structure in this hierarchy, 
with
$B^{\M'}\neq\emptyset$, but $B^{\M'}$ coding a non-$\Tt'$-cofinal (for the relevant tree $\Tt'$) 
branch $[0,\om\gamma')_{\Tt'}$ (for some $\om\gamma'<\lh(\Tt')$). Let $\pi:\M\to\M'$ be
fully elementary. Then clearly $B^\M$ codes $[0,\om\gamma)_\Tt$ where
$\pi(\Tt)=\Tt'$ and $\pi(\gamma)=\gamma'$, and $\om\gamma<\lh(\Tt)$. But we
need that $[0,\om\gamma)_\Tt\sub\Sigma(\Tt)$, and this is not clear
(even if $\Sigma$ has hull condensation).

The other correction, which is better, is to simply not feed in $\Sigma(\Tt)$ in
the case that $\lh(\Tt)$ has measurable
cofinality in $\M|(\alpha+\lambda)$ (as witnessed by some measure on $\es^\M$).
For by the argument in \ref{rem:tree_cof_meas}, $\M$ already has $\Sigma(\Tt)$ as an element,
and there is a uniform procedure which $\M$ can use to determine $\Sigma(\Tt)$.

Thus, one must show that the relevant ultrapowers and substructures of models
in the resulting hierarchy are also in the hierachy. It is easy to see that
ultrapowers preserve the relevant first-order properties.

So let $\M'$ be
a $\Sigma$-premouse and let $\pi:\M\to\M'$ be a weak $0$-embedding. We want to
know that $\M$ is a $\Sigma$-premouse, given that $\Sigma$ has hull condensation. We just need to 
verify the first-order
properties.

We need to rule out the possibility that
$B^\M=\emptyset$ (and therefore $B^{\M'}=\emptyset$), but there is some
$B\neq\emptyset$ such that $(\M,B)$ is a $\Sigma$-premouse. Let $\Tt\in\M$ be
the relevant tree (with $B$ coding $\Sigma(\Tt)$). Because $\pi$ is a weak
$0$-embedding, this implies that $\Tt'=\pi(\Tt)$ is the least tree for which
$\M'$ does not know $\Sigma(\Tt')$, and $\pi$ is discontinuous at $\lh(\Tt)$.
Suppose also that $\M=\core_1(\M')$ and $\pi$ is the core map,
and $\M'$ is $(0,\om_1,\om_1+1)$-iterable. Then by the usual proof of
solidity (with a little extra argument to deal with the possibility that $\M$ is
not a $\Sigma$-premouse), $\M$ and $\M'$ are $1$-solid and
$\pi(p_1^\M)=p_1^{\M'}$, and then using the comparison argument in the proof of
universality, and the commutativity of $\pi$ with the resulting iteration
embeddings, one gets that $\lh(\Tt)$ has measurable cofinality in
$\M$, and therefore $\M$ is in fact a $\Sigma$-premouse, contradiction. (For the
higher degree core maps, the present situation cannot arise, just by
elementarity.)

Now suppose that $B^{\M'}\neq\emptyset$. It is easy to see that $B^\M$ codes
some branch $b$ through $\Tt$, and that $B^\M\inter\M$ is cofinal in
$\OR(\M)$ (by the $\Sigma_1$-elementarity of $\pi$ on a set cofinal in $\OR(\M)$).
But $b$ need not be $\Tt$-cofinal. (For example, if $\OR(\M')$ has
uncountable cofinality, it is easy to find $\N\pins\M$ such that letting
$\M=(\N,B^{\M'}\inter\N)$ and $\pi=\id$, then $\pi:\M\to\M'$
is a weak $0$-embedding, and $\Tt=\Tt'$.) If we have that $\pi$ is
$\Sigma_1$-elementary on a set $X\sub\OR(\M)$ which is both cofinal in $\OR(\M)$
and cofinal in $\lh(\Tt)$, then $b$ will be cofinal in $\Tt$.

These arguments give that the models produced in an $L[\es,\Sigma]$-construction
will all be $\Sigma$-mice, as long as iterates of countable
elementary substructures are realizable back into models of the construction,
in the usual manner. But we opted for the hierarchy for $\Sigma$-premice defined
in \S\ref{StrategyPremice} because it has stronger condensation properties, and without
assuming any iterability.
\end{remark}

\section{G-organization}
\label{G-organized}

Let $\Omega$ be either an operator or an iteration strategy. In this section we 
implement some ideas of Grigor Sargsyan, defining \emph{$\g$-organized $\Omega$-premice}.
This will be useful assuming that $\Omega$ has a certain absoluteness property, which we first 
describe.

\begin{definition} Let $a$ be transitive and $A=\hat{a}$. We say that $A$ is 
\dfnemph{self-wellordered \tu{(}swo'd\tu{)}} iff
$a=\trancl(x\un\{x,\prec\})$ for some transitive set $x$, and wellorder $\prec$ of $x$.
For swo'd $A$ and $\prec$ as above, let $\prec_A$ denote the canonical wellorder of $A$ determined 
by $\prec$.
\end{definition}

\begin{definition}
Let $\psi$ be a $\Sigma_0$ formula in the language of set theory.\footnote{$\psi$ will be 
used 
to restrict the class of 
iteration trees being considered; for example, $\psi(x)$ might say that ``$x$ is a normal tree''} 
Then $\varphi_{\psi,\min}(x)$ denotes the formula in the free variable $x$ asserting,
over abstract spms $\M$ with $\hmb^\M$ swo'd: ``Let $\prec$ be the canonical wellorder of 
the universe. Then $x$ 
is the $\prec$-least limit length iteration tree $\Tt$ on $\hmp$ according to $\Sigma^V$ such that 
$\Sigma^V(\Tt)$ is undefined and $\psi(\Tt)$ holds''.

Let $\varphi_{\min}=\varphi_{\text{``true''},\min}$.
\end{definition}

\begin{definition}\label{dfn:M_1^Sigma,sharp(A)}
Either:
\begin{enumerate}[label=--]
\item let $\Sigma$ be an iteration strategy for a transitive structure 
$\P\in\HC$, let $\varphi\in\Ll^+$ and $X=(\Sigma,\varphi)$,
let $a\in\HC$ and $A=\hat{a}$ with $\P\in\J_1(A)$, or
\item let $X=\Fop$ be an operator over $\opbk$, $\kappa\leq\OR(\opbk)$ be an 
uncountable cardinal
and $A\in\witri{C_\Fop}\inter\HC$.
\end{enumerate}

We write $\M_1^{X,\#}(A)$ for the (unique) sound, 
non-$1$-small $X$-pm $\M$ over $A$, such that $\M$ is $X$-$(0,\om_1)$-iterable,
and if $\cof(\om_1)>\om$, $\M$ is $X$-$(0,\om_1+1)$-iterable (given such an $\M$ 
exists).\footnote{$\ZF$ 
proves uniqueness. For let $\M\neq\N$ be such $X$-pms. Let $(\Tt,\Uu)$ be their length 
$\om_1$ comparison if $\cof(\om_1)=\om$, or length $\om_1+1$ comparison otherwise. Let 
$z\in\RR$ code $(\M,\N)$ and
let $W=L[z,\Tt,\Uu]$. Then $\M,\N\in\HC^W$ and $W\sats\AC$, and therefore if 
$\cof(\om_1)=\om$ then $W\sats$``$\gamma$ is a limit cardinal'', where $\gamma=\om_1$.
So working in $W$ we can reach a contradiction as 
usual.} Let $\kappa$ be an uncountable cardinal. We say that $\M_1^{X,\#}(A)$ is 
\dfnemph{$X$-$\kappa$-naturally iterable} iff $\MFsharp=\M_1^{X,\#}(A)$ exists and either:
\begin{enumerate}[label=\tu{(}\alph*\tu{)}]
 \item\label{item:cof(kappa)>om} $\cof(\kappa)>\om$ and $\MFsharp$ is $X$-$(0,\kappa+1)$-iterable, 
or
 \item\label{item:cof(kappa)=om} $\cof(\kappa)=\om$ and $\MFsharp$ is $X$-$(0,\kappa)$-iterable.
\end{enumerate}
When this holds, let $\Lambda^{X,\kappa}_\MFsharp$ denote the unique\footnote{Much as 
before,
$\ZF$ proves uniqueness.} $X$-$(0,\kappa)$-strategy for 
$\MFsharp$ which, if $\cof(\kappa)>\om$, extends to an $X$-$(0,\kappa+1)$-strategy;
also if $\cof(\kappa)>\om$ let $\Lambda^{X,\kappa+1}_\MFsharp$ denote the unique\footnote{Likewise.}
$X$-$(0,\kappa+1)$-strategy for $\MFsharp$.

Let $\Sigma$ be an iteration strategy for a transitive structure $\P\in\HC$,
with recognizable domain, as witnessed by a $\Sigma_0$ formula $\psi$ (in the language of set 
theory), with $\psi$ least such. Then $\varphi_{\min}^{\Sigma}$ denotes
$\varphi_{\psi,\min}$. Let $\varphi=\varphi^\Sigma_{\min}$. Then we abbreviate the pair 
$(\Sigma,\varphi)$ with $\Sigma$. So a $\Sigma$-pm is a 
$(\Sigma,\varphi)$-pm, etc.
\end{definition}

\begin{definition}\label{dfn:suitable}
We say that $(\Omega,\varphi,X,A,\kappa)$ is 
\dfnemph{suitable} iff $\kappa$ is an uncountable cardinal, $A=\hat{a}$ for some transitive 
$a\in\HC$,
$\M_1^{X,\#}(A)$ exists and is 
$X$-$\kappa$-naturally iterable, and either:
\begin{enumerate}[label=--]
 \item $\Omega=\Sigma$ is a $\kappa$-strategy with hull condensation and recognizable domain, for a 
transitive structure $\P\in\HC\inter\J(A)$, $\varphi\in\Ll^+$, and $X=(\Sigma,\varphi)$, or
\item $\Omega=X=\Fop$ is a total operator over $\opbk$, where $\opbk$ is an operator background 
with 
$\kappa=\OR(\opbk)$,  
$C_\Fop$ is the (possibly swo'd) cone of $\opbk$ 
above $a$, and $\Fop$ condenses finely above $a$.
\end{enumerate}

For suitable $t=(\Omega,\varphi,X,A,\kappa)$, let $\Omega_t=\Omega$, etc.

Let $(\Omega,A)$ be given. We say that $(\Omega,A)$ is \dfnemph{suitable} iff
either (i) $\Omega$ is a $\kappa$-strategy $\Sigma$, with recognizable domain, for some transitive 
structure $\P\in\HC$ and uncountable cardinal $\kappa$,  $A=\hat{a}$ for some $a\in\HC$, $A$ is 
swo'd, and 
\[ t_{\Sigma,A}\eqdef(\Sigma,\varphi^{\Sigma}_{\min},(\Sigma,\varphi^{\Sigma}_{\min}),A,\kappa)\] 
is suitable, or (ii) $\Omega$ is an operator $\Fop$ over $\opbk$, and letting $\kappa=\OR(\opbk)$,
\[ t_{\Fop,A}\eqdef(\Fop,0,\Fop,A,\kappa) \]
is suitable.
\end{definition}

\begin{lemma}\label{lem:Sigma_N_condensation}
Let $t=(\Omega,\varphi,X,A,\kappa)$ be suitable, $\MFsharp=\M_1^{X,\#}(A)$ and 
$\eta=\cof(\kappa)$. Then:
\begin{enumerate}
\item\label{item:branch_hull_con} $\Lambda^{X,\kappa}_\MFsharp$ has branch condensation and hull 
condensation.
\item\label{item:stacks_max_iter_cof=om} If $\eta=\om$ then $\MFsharp$ is
$X$-$(0,{<\om},\kappa)$-maximally iterable.
\item\label{item:stacks_max_iter_cof>om} If $\eta>\om$ 
then $\MFsharp$ is $X$-$(0,\eta,\kappa+1)$-maximally 
iterable.\footnote{We also get 
$X$-$(0,\eta,\kappa+1)$-iterability (without \emph{maximal}) for the strategy case,
but for reasons covered in \cite{operators}, we cannot expect the same if 
$X$ is an operator.}
\end{enumerate}
\end{lemma}
\begin{proof}
These facts come from the uniqueness of $\Lambda^{X,\kappa}_\MFsharp$, together with 
the the condensation properties proved in this section 
for strategy mice, and the condensation properties and copying arguments of 
\cite{operators} in the case that $X$ is an operator.
Part \ref{item:branch_hull_con} follows routinely from these items.
Parts \ref{item:stacks_max_iter_cof=om} and \ref{item:stacks_max_iter_cof>om} are essentially by 
\cite[Theorem 3.1(?)]{iter_for_stacks}.
The latter results are literally stated and proved only for standard premice, but the arguments 
there go through using the properties and arguments just mentioned.
\end{proof}

\begin{remark}\label{rem:stacks_strategy}
What is behind the foregoing proof (in terms of the details contained in \cite{iter_for_stacks}), 
is as follows. If $\eta=\om$ let $\Lambda=\Lambda^{X,\kappa}_\MFsharp$ and $\theta=\kappa$.
If $\eta>\om$ let $\Lambda=\Lambda^{X,\kappa+1}_\MFsharp$ and $\theta=\kappa+1$.
An $X$-$(0,{<\eta},\theta)$-maximal strategy $\Psi$ for $\M$ is computed,
extending $\Lambda$ (and therefore $\Lambda^{X,\kappa}_\MFsharp$), and such that the restriction of 
$\Psi$ to an $X$-$(0,{<\eta},\kappa)$-maximal strategy $\Psi'$, lifts to 
$\Lambda^{X,\kappa}_\MFsharp\sub\Lambda$. 
(Stacks via $\Psi$ which have a last tree $\Tt$ of length $\kappa+1$ can lift to a normal tree 
$\Uu$ of length $>\kappa+1$, in which case $\Uu$ cannot literally be via $\Lambda$, but for 
instance,
$\Uu\rest\kappa+1$ is via $\Lambda$.)

If $\Tt$, via $\Psi'$, has successor length, then $M^\Tt_\infty$ is 
$X$-$(\deg^\Tt(\infty),{<\eta},\theta)$-maximally iterable,
via the tail 
$\Psi^*$ 
of $\Psi$.
Moreover, given a normal tree $\Uu$ on $M^\Tt_\infty$ of limit length $<\kappa$, via 
$\Psi^*$, and $c=\Psi^*(\Uu)$, either there is a 
Q-structure for $M(\Tt)$ in $L^{X}_\kappa(M(\Tt))$, which determines $c$ as usual, or else 
neither $b^\Tt$ nor $c$ drop 
in model or degree and $i^\Uu_c\com i^\Tt(\delta^\MFsharp)=\delta(\Uu)$.

If $\eta>\om$ then clearly any
$X$-$(n,{<\eta},\kappa+1)$-maximal strategy extends to an $X$-$(n,\eta,\kappa+1)$-maximal 
strategy. So part 
\ref{item:stacks_max_iter_cof>om} follows readily from the above.
Note also that any strategy witnessing part
\ref{item:stacks_max_iter_cof=om} (part \ref{item:stacks_max_iter_cof>om}) must extend 
$\Lambda^{X,\kappa}_\MFsharp$ (must extend $\Lambda^{X,\kappa+1}_\MFsharp$).
\end{remark}

\begin{definition}
In the preceding context, let $\Lambda^{X,({<\eta},\kappa)}_\MFsharp$ denote $\Psi'$.
\end{definition}

The following absoluteness property ensures that g-organization is useful:

\begin{definition}\label{dfn:determines}
Let $t=(\Omega,\varphi,X,A,\kappa)$ be suitable and $\MFsharp=\M_1^{X,\#}(A)$.
We say that $t$ \dfnemph{determines itself 
on generic extensions} iff there
are formulas $\Phi,\Psi$ in
$\Ll^+$ and some $\gamma > \delta^\MFsharp$ such that 
$\MFsharp|\gamma
\vDash\Phi$ and for any non-dropping
$\Lambda^{X,\kappa}_\MFsharp$-iterate $\N$ of
$\MFsharp$ via a countable tree $\Tt$ based on $\MFsharp|\delta^{\MFsharp}$, any $\N$-cardinal 
$\delta$, 
any
$\gamma\in\Ord$ such that $\N|\gamma\models\Phi\ \&\ $``$\delta$ is
Woodin'', and any $g$ which is set-generic over $\N|\gamma$ (with $g\in V$),
we have that $\R\eqdef(\N|\gamma)[g]$
is closed under $\Omega$, and $\Omega\rest\R$ is defined 
over
$\R$ by
$\Psi$. We say such a pair $(\Phi,\Psi)$ \dfnemph{generically determines 
$t$}.

Let $A\in\HC$ and let $\Omega$ be either an operator or an iteration strategy.
We say that $(\Omega,A)$ is \dfnemph{nice} iff $(\Omega,A)$ is 
suitable and $t_{\Omega,A}$ determines itself on generic extensions. We say that $(\Phi,\Psi)$ 
\dfnemph{generically determines}
$(\Omega,A)$ iff $(\Phi,\Psi)$ generically determines $t_{\Omega,A}$.
\end{definition}

\begin{lemma}\label{lem:determines_extend}
Let $\N,\delta$, etc, be as in \ref{dfn:determines}, except that we allow 
$\Tt$ to have any length ${<\kappa}$, and allow $g$ to be in a set-generic extension of $V$. Then 
$\R$ is closed 
under $\Omega$ 
and $\Omega'\rest\dom(\Omega)=\Omega\rest\R$ where $\Omega'$ is the interpretation of $\Psi$ over 
$\R$.
\end{lemma}
\begin{proof}
We first give the proof assuming that $\Omega=\Sigma$ is a strategy, and then point out the 
differences 
for the other case. Suppose the lemma fails. Let $x\in\R$ be a counterexample to the claimed
agreement between $\Sigma,\Sigma'$. So $\Uu\eqdef x\in\dom(\Sigma)\sub V$. Let $\PP$ be some 
forcing, 
and $H\sub\PP$ be $V$-generic, such that $g\in V[H]$. Because $a\in\HC$, $\N$ is wellorderable,
and so by $\Sigma^1_1$-absoluteness, we may assume $\PP=\Coll(\om,\OR(\N))$.
Moreover, letting $z\in\RR$ code $a,M^\Uu_0,\MFsharp$, we may assume that $g\in W\eqdef 
L[z,\Tt,\Uu,\Sigma(\Uu)]$.

Work in $W$, where $\AC$ holds.
Let $\dot{g}$ be a $\PP$-name for $g$. Let $\dot{\Uu}\in\N|\gamma$ be such that 
$\dot{\Uu}^g=\Uu$. Fix $p\in H$ forcing ``$\dot{g}$ is 
$\check{\N}|\check{\gamma}$-generic 
and $\check{\dot{\Uu}}^{\dot{g}}=\check{\Uu}$''; for simplicity assume that $p=\emptyset$. Let
$\alpha$ be large and let
\[ \pi:M\to L_\alpha[z,\Tt,\Uu,\Sigma(\Uu)] \]
be elementary, with $M$ 
countable and all 
relevant objects in $\rg(\pi)$. Write
$\pi(\Ttbar)=\Tt$, etc.

Now work in $V$. Note that $\bar{\Uu}$ is via $\Sigma$ and $\bar{\Uu}\in\dom(\Sigma)$ because 
$\Sigma$ has hull condensation and recognizable domain. By \ref{lem:Sigma_N_condensation}, $\Ttbar$ 
is 
via $\Lambda^X_\MFsharp$. For 
any $H^*$ which is $\bar{\PP}$-generic over $M$, letting $g^*=\bar{\gdot}^{H^*}$, 
we then have
\[ \overline{\Uudot}^{g^*}=\bar{\Uu}\in\overline{\N|\gamma}[g^*],\]
and letting $\Sigma^*$ be 
the 
interpretation of $\Psi$ over $\overline{\N|\gamma}[g^*]$, 
by \ref{dfn:determines} we have
\begin{equation}\label{eqn:agmt_G} 
\Sigma(\bar{\Uu})=\Sigma^*(\bar{\Uu})\in\overline{\N|\gamma}[g^*].\end{equation}
So $\Uu\in\dom(\Sigma')$ (by the above, this is forced by $\PP$),
and so $\Sigma'(\Uu) \neq\Sigma(\Uu)$, by choice of $\Uu$. 
By hull condensation, $\overline{\Sigma(\Uu)} = \Sigma(\bar{\Uu})$, and so by 
line (\ref{eqn:agmt_G}), $\overline{\Sigma(\Uu)}=\Sigma^*(\bar{\Uu})$ for any $H^*$. So in
$M$, 
$\bar{\PP}$ forces that $\overline{\Sigma(\Uu)} = \Sigma^*(\bar{\Uu})$. Therefore 
$\PP$ forces that $\Sigma(\Uu) = \Sigma'(\Uu)$. Contradiction.

Now consider the case that $\Omega=\Fop$ is an operator. The argument is almost the same. The 
coarse 
condensation (a component of fine condensation) of $\Fop$ above $a$, and the fact that $a\in\HC$, 
replaces the use of hull condensation and the recognizability of $\dom(\Sigma)$. Much as 
before, we can assume that 
$\PP=\Coll(\om,Z)$ for some transitive $Z\in\opbk$. Because $\opbk\sats\DC$ we can form 
an appropriate countable elementary substructure $M$ of some large enough set in $\opbk$. We omit 
further detail.
\end{proof}

We next consider some issues pertaining to hod mice;
see \cite{hod_mice} for background.\footnote{We assume only a basic knowledge of hod mice; more 
than enough is covered in the first sections of \cite{hod_mice}. As mentioned earlier, the actual 
analysis of scales does not depend particularly on the theory of hod mice, and is developed in 
parallel for standard mice.}

\begin{definition}
 A pointclass is \dfnemph{smooth} iff it contains all open sets and is closed under continuous 
preimage, intersections, unions and real quantifiers.
\end{definition}

\begin{remark}\label{rem:Lp^Gamma,Sigma}
Assume that $\om_1$ is regular and let $\Gamma$ be smooth pointclass.
Let $a\in\HC$ be swo'd.
Let $\Sigmavec$ be the join of a sequence of strategies for a sequence $\vec{\P}$ of transitive 
structures in $\J(a)$ (possibly the sequence has length $0$ or $1$).
As in \cite[Definition 2.26]{hod_mice}, $\Lp^{\Gamma,\Sigmavec}(a)$ denotes the stack 
of all sound $\Sigmavec$-premice $\M$ over $a$ which project to $a$,
such that in $\Gamma$ there is a $\Sigmavec$-$(\om,\om_1,\om_1)$-iteration strategy for $\M$ 
which extends to a $\Sigmavec$-$(\om,\om_1,\om_1+1)$-strategy.\footnote{In \cite{hod_mice}, 
the definition is stated in the context of $\AD^+$, so the extension to $\om_1+1$ exists. Here as 
elsewhere, a $\Sigmavec$-$(\om,\om_1,\om_1+1)$-strategy is only required to ensure wellfoundedness 
of the last model of successor length trees of size $\om_1$, not $\Sigmavec$-correctness.}
Here we are demanding a full $\Sigmavec$-$(\om,\om_1,\om_1+1)$-strategy,
not just a hod strategy. This is somewhat at odds with our usual practice in this paper, of dealing 
only with hod strategies for hod premice; it is done for consistency with \cite{hod_mice}.
Fortunately, if each strategy in $\Sigmavec$ has hull condensation, we could have actually
defined $\Lp^{\Gamma,\Sigmavec}$ using hod strategies, or in fact using normal strategies, and 
gotten the same result:
\end{remark}

\begin{lemma}
Suppose $\om_1$ is regular and let $\Gamma,a,\vec{\P},\Sigmavec$ be as in 
\ref{rem:Lp^Gamma,Sigma}. Suppose that every strategy in $\Sigmavec$ has hull condensation.
Then $\Lp^{\Gamma,\Sigmavec}(a)$ is the stack of all sound $\Sigmavec$-premice over 
$a$ which project to $a$ and such that there is a $\Sigmavec$-$(\om,\om_1)$-strategy for $\N$ in 
$\Gamma$ which extends to a $\Sigmavec$-$(\om,\om_1+1)$-strategy.
\end{lemma}
\begin{proof}
This is by the proof of \ref{lem:Sigma_N_condensation}, together with 
\cite[\S3(?)]{iter_for_stacks} and the closure of $\Gamma$ under real quantifiers.
\end{proof}

\begin{definition}
 Let $\P$ be a hod premouse and $\R\pins\Ss\pins\P$
 be such that $\R$ is a cutpoint of $\Ss$
 and $\Ss\pins\P(\alpha)$ where $\alpha$ is least such that $\R\pins\P(\alpha)$.
 Suppose either $\Ss$ projects $\leq\OR(\R)$, or $\OR(\R)$ is the largest cardinal of $\Ss$.
 Then $\Ss^*(\R)$ denotes the \emph{$*$-translation} of $\Ss$ above $\R$
 (much as in \cite[\S7]{DMATM}; so $\Ss^*(\R)$ is approximately a strategy premouse over $\R$, and 
in particular, $\OR(\R)$ is a \emph{strong} cutpoint of 
$\Ss^*(\R)$). If $\OR(\R)$ is the largest cardinal of $\Ss$ then $\Ss^*$ denotes $\Ss^*(\R)$.
\end{definition}

We now 
define \emph{$\Gamma$-fullness$^*$ preserving} much as 
\emph{$\Gamma$-fullness preserving} is defined in \cite[Definition 2.27]{hod_mice}, but with a few 
modifications, the most significant of which is that we make requirements regarding dropping 
iterates, and related to this, the fact that we consider all cutpoints, not just 
strong cutpoints.
(thus, because $\R$ is, by definition, a strong cutpoint of $\Lp^{\Gamma,\Sigma}(\R)$, we must use 
$\Ss^*$ where $\Ss$ is used in \cite{hod_mice}).

\begin{definition}\label{dfn:Gamma-fullness*}
Suppose $\om_1$ is regular and $(\P,\Sigma)$ is a hod pair with $\P\in\HC$ and $\Gamma$ is a 
smooth pointclass. Then $\Sigma$ is \dfnemph{$\Gamma$-fullness$^*$ 
preserving} 
iff the following two conditions hold:
\begin{enumerate}
 \item\label{item:Gamma-fullness*_non-dropping} Let $(\Ttvec,\Q)\in I(\P,\Sigma)\inter\HC$ and 
$\gamma\leq\lambda^\Q$. Then
\begin{enumerate}[label=--]
 \item for all cutpoints (not just strong) $\eta$ of $\Q(0)$,
\[ (\Q|(\eta^+)^{\Q(0)})^*=\Lp^{\Gamma}(\Q|\eta), \]
 \item if $\gamma=\alpha+1$ then for all cutpoints $\eta$ of $\Q(\alpha+1)$ with 
$\eta\geq\OR(\Q(\alpha))$,
\[ (\Q|(\eta^+)^{\Q(\alpha+1)})^*=\Lp^{\Gamma,\Sigma_{\Q(\alpha),\Ttvec}}(\Q|\eta), \]
 \item and if $\gamma$ is a limit then for all cutpoints $\eta$ of $\Q(\gamma)$ with 
$\eta\geq\delta_\gamma^\Q$,
\[ 
(\Q|(\eta^+)^{\Q(\gamma)})^*=\Lp^{\Gamma,\oplus_{\beta<\gamma}\Sigma_{\Q(\beta),\Ttvec}}
(\Q|\eta). \]
\end{enumerate}
\item\label{item:Gamma-fullness*_dropping} Let $(\Ttvec,\Tt)$ be a countable tree via $\Sigma$, 
consisting of a stack $\Ttvec$ followed 
by a normal tree $\Tt$, such that $\Tt$ has successor length and $b^\Tt$ drops.
Let $\Q=M^{\Tt}_\infty$ and $\lambda=\lambda^\Q$.
Let $\gamma$ be least such that $\OR(\Q(\lambda))<\lh(E^\Tt_\gamma)$ and
let $\Uu=\Ttvec\conc(\Tt\rest(\gamma+1))$. (Note $b^\Uu$ does not drop.)
Let $\R,\Ss$ be such that
$\Q(\lambda)\ins\R\pins\Ss\ins\Q$ and $\R$ is a cutpoint of $\Ss$
and $\Ss$ projects $\leq\OR(\R)$ and is $\OR(\R)$-sound (so either $\Ss\pins\Q$ or all 
generators of $\Tt$ are $<\OR(\R)$). Then
\[ \Ss^*(\R) \pins\Lp^{\Gamma,\Sigma_{\Q(\lambda),\Uu}}(\R).\qedhere \]
\end{enumerate}
\end{definition}

\begin{definition}
 Let $(\P,\Sigma)$ be a hod pair with $\P\in\HC$. We say that $\Sigma$ has \dfnemph{weak hull 
condensation} iff for all transitive $W,X$ satisfying $\ZF^-$, with $W\in\HC$, and fully elementary 
$\pi:W\to X$,
if $\P\in\HC^W$ and $\Ttvec\in W$ and $\pi(\Ttvec)$ is a stack on $\P$ via $\Sigma$,
then $\Ttvec$ is via $\Sigma$.
\end{definition}

\begin{definition}
A premouse or hod premouse $\P$ is \dfnemph{reasonable} iff $\P$ is super-small,
all $\N\ins\P$ (including $\N=\P$) satisfy the conclusions of \cite[4.11, 4.12,
4.15]{mim},
and if $\P$ is a premouse then all $\N\pins\P$ are ${<\om}$-condensing,
and if $\P$ is a hod premouse then for all $\N\pins\P$, $\N$ is ${<\om}$-condensing with respect to 
embeddings $\pi:\Hh\to\N$ such that $\crit(\pi)\geq\delta_\alpha^\P$ for all 
$\alpha$ such that $\delta_\alpha^\P\leq\OR(\N)$. Reasonableness is 
preserved by fine structural iteration, as super-smallness is $\rSigma_2$ and the other conditions 
are $\rPi_1$.

A hod pair $(\Sigma,\P)$ is \dfnemph{within scope} iff $\DC_\RR$ holds, $\P\in\HC$ is reasonable,
is below $\AD_\RR$+``$\Theta$ is regular'',
$\Sigma$ is a \emph{hod} $(\om,\kappa,\kappa+1)$-strategy 
for $\P$, where $\kappa$ is some regular uncountable cardinal,
$\Sigma$ is $\Gamma$-fullness$^*$ preserving for some smooth pointclass $\Gamma$,
$\Sigma$ has branch condensation, and if $\kappa>\om_1$ then $\Sigma$ has weak hull condensation.
\end{definition}

\begin{definition}
Let $(\P,\Sigma)$ be a hod pair. We say that $\Sigma$ has \dfnemph{factor hull condensation}
iff whenever:
\begin{enumerate}[label=--]
 \item $\Ttvec,\Uuvec$ are stacks via $\Sigma$ and $i^{\Ttvec},i^{\Uuvec}$ exist;
 let $\M=M^{\Ttvec}_\infty$ and $\N=M^{\Uuvec}_\infty$,
 \item $\pi:\M\to\N$ is elementary and $\pi\com i^{\Ttvec}=i^{\Uuvec}$,
\item $\Wwvec$ is a stack on $\N$ via $\Sigma_{\N,\Uuvec}$, and
\item $\Vvvec$ is a stack on $\M$ and $\pi\Vvvec$ is a hull of $\Wwvec$,
\end{enumerate}
then $\Vvvec$ is via $\Sigma_{\M,\Ttvec}$.
\end{definition}

Factor hull condensation trivially implies hull condensation.
But the following lemma is more interesting; part of its proof uses ideas similar to those in 
Sargsyan's \cite[Proposition 2.41]{hod_mice}.

\begin{lemma}\label{lem:factor-hull_condensation}
 Let $(\P,\Sigma)$ be a hod pair within scope. Then $\Sigma$ has factor hull condensation.
\end{lemma}
\begin{proof}
By weak hull condensation, we may assume that all trees we deal with are countable.
(If $\kappa=\om_1$ then because $\om_1$ is regular, it is easy to see that we may still reduce to 
countable trees, without using weak hull condensation.)

Suppose the lemma fails. Let
$\Ttvec,\Uuvec,\M,\N,\pi$ be a counterexample. A \dfnemph{bad system} is a countable system
\[ 
\left(\left<\Vvvec_i,\Vvvec_i^*,\Wwvec_i,\Wwvec_i^*,\sigmavec_i,\sigmavec_i^*,\beta_i\right>_{i\leq 
n},\left<\alpha_i,\pi_i\right>_{i\leq n+1}\right) \]
where 
\begin{enumerate}
 \item $(\Ttvec,\Vvvec_0,\ldots,\Vvvec_n)$ and $(\Uuvec,\Wwvec_0,\ldots,\Wwvec_n)$ are terminally 
non-dropping stacks on $\P$ via $\Sigma$. Let $\M_0=\M$ and $\N_0=\N$ and 
$\M_{i+1}=M^{\Vvvec_{i}}_\infty$ and $\N_{i+1}=M^{\Wwvec_i}_\infty$.
 \item $\beta_0=\lambda^\M$ and $\alpha_0=\lambda^\M+1$ and $\beta_i<\alpha_i$ and
 $\alpha_{i+1}\leq i^{\Vvvec_i}(\beta_i)$.
 \item $\Vvvec_{i}$ is based on $\M_{i}(\beta_{i})$.
 \item $\Vvvec_i=(\Vvvec_i',\Vv_i)$, where $\Vv_i$ is a normal tree (so $\Vv_i$ 
is terminally non-dropping and $\M_{i+1}=M^{\Vv_i}_\infty$).
 \item $\Vvvec_i^*=(\Vvvec_i',\Vv_i^*)$ is a stack on $\M_i$, based on $\M_i(\beta_i)$,
 where $\Vv_i^*$ is a normal extension of $\Vv_i$, 
$\Vv_i=\Vv_i^*\rest(\gamma_i+1)$, where $\gamma_i$ is the least $\gamma$ such that
$\OR(\M_{i+1}(\alpha))<\lh(E^{\Vv_i^*}_\gamma)$ for all $\alpha<\alpha_{i+1}$, 
$\Vv_i^*\rest[\gamma_i,\lh(\Vv_i^*))$ is based on $\M_{i+1}(\alpha_{i+1})$, and $\Vv_i^*$ has 
successor length and is terminally dropping.
 \item $\Wwvec_i=(\Wwvec_i',\Ww_i)$, where $\Ww_i$ is a normal tree.
 \item $\Wwvec_i^*=(\Wwvec_i',\Ww_i^*)$, where $\Ww_i^*$ is a normal extension of $\Ww_i$.
 \item $(\Uuvec,\Wwvec_0,\ldots,\Wwvec_{i-1},\Wwvec_i^*)$ is via $\Sigma$.
 \item $(\Ttvec,\Vvvec_0,\ldots,\Vvvec_{i-1},\Vvvec_i^*)$ is not via $\Sigma$.
\item $(\Ttvec,\Vvvec_0,\ldots,\Vvvec_{i-1},\Vvvec_i^*\rest\lh((\Vvvec_i^*)-1))$ is via $\Sigma$.
 \item $\pi_0=\pi$ and $\pi_{i}:\M_{i}\to\N_{i}$,
 and $\pi_i\Vvvec_i$ is a hull of $\Wwvec_i$, as witnessed by $\sigmavec_i$, 
with the final node of $\pi_i\Vvvec_i$ corresponding to the final node of $\Wwvec_i$, and
$\pi_{i+1}$ is the composition of the final copy and hull embedding maps.
 \item $\pi_i\Vvvec_i^*$ is a hull of $\Wwvec_i^*$, as witnessed by $\sigmavec_i^*$,
 and $\sigmavec_i\sub\sigmavec_i^*$, and $\Wwvec_i^*$ has no proper segment $\Ww'$ such that 
$\pi_i\Vvvec_i^*$ is a hull of $\Ww'$ as witnessed by $\sigmavec_i^*$.
\end{enumerate}

Because of our choice of $\Ttvec,\Uuvec,\pi$, it is easy to see using branch condensation and weak 
hull condensation (the latter to give countability, and the former to ensure a dropping branch) 
that there is a bad system with $n=0$. Using $\DC_\RR$, it follows that there is a bad system $B$ 
for which no proper extension is also a bad system. Let the notation above be used to describe $B$.

Let $\R=\M_{n+1}(\alpha_{n+1})$ and $\varrho=\pi_{n+1}\rest\R$ and 
$\Ss=\N_{n+1}(\varrho(\alpha_{n+1}))$, so $\varrho:\R\to\Ss$ is elementary.
Let $\eta=\sup_{\alpha<\alpha_{n+1}}\OR(\R(\alpha))$.
Let $\eta^\Ss=\varrho(\eta)$. Let $\Psi^\Ss$ be the above-$\eta^\Ss$,
$(\om,\om_1+1)$-strategy for $\Ss$ given by normally extending $\Wwvec^*_n$, and continuing to use 
$\Sigma$. Let $\Psi$ be the $\varrho$-pullback of $\Psi^\Ss$, for $\R$.
Let $\Xx=(\Ttvec,\Vvvec_0,\ldots,\Vvvec_n)$.

\begin{claim}\label{clm:Psi_is_correct} $\Psi$ is a 
$\oplus_{\alpha<\alpha_{n+1}}\Sigma_{\R(\alpha),\Xx}$-strategy.
\end{claim}
\begin{proof}
If not then, again using weak hull 
and branch condensation, it is easy to produce a bad system properly extending $B$, a contradiction.
\end{proof}

Now let $\Vv$ be the tree on $\R$ which is equivalent to 
$\Vvvec_n^*\rest[\gamma_n,\lh(\Vvvec_n^*)-1)$, let $b,c$ be the $\Vv$-cofinal branches determined 
by $\Vvvec_n^*$ and $\Sigma$ respectively. So $b\neq c$ and $b$ drops. By the claim and using 
$\Sigma$,
we may successfully compare the phalanxes $\Phi(\Vv\conc b)$ and $\Phi(\Vv\conc c)$, producing 
(padded) trees $\Yy,\Zz$ extending $\Vv\conc b$ and $\Vv\conc c$ respectively. Moreover, all models 
of $\Yy,\Zz$ are $\oplus_{\alpha<\alpha_{n+1}}\Sigma_{\R(\alpha),\Xx}$-hod premice. Let 
$\delta=\delta(\Vv)$ and $\lambda=\lh(\Vv)$.

\begin{claim} $c$ does not drop, and therefore $\alpha_{n+1}$ is not a limit.\end{claim}
\begin{proof} This is a standard argument, but we give it as it is not too long,
and we need it elsewhere. Suppose $c$ drops. 
It suffices to see that at for every $\alpha\geq\lambda$, $[0,\alpha]_\Yy$ and $[0,\alpha]_\Zz$ 
drops, since then standard fine structure yields a contradiction.
Suppose this fails. Then there is $\alpha\geq\lambda$ such that either 
$E=E^\Yy_\alpha$, or $E=E^\Zz_\alpha$, has $\crit(E)<\delta$. Let $\alpha$ be least such.
 Then $[0,\alpha']_\Yy$ and $[0,\alpha']_\Zz$ drop for each $\alpha'\in[\lambda,\alpha]$.
Let $\Q_b=\Q(\Vv,b)$ and $\Q_c=\Q(\Vv,c)$.
Then $\Q_b\neq\Q_c$, so $\delta$ is Woodin in $M^\Yy_{\alpha}||\lh(E)$. So if 
there is any $F\in\es^{M^\Yy_\alpha||\lh(E)}$ such that $\crit(F)<\delta<\lh(F)$, we easily get 
that $[0,\beta]_\Yy$ and $[0,\beta]_\Zz$ are dropping for all $\beta>\alpha$ (as Woodins are 
cutpoints of hod premice). So suppose $E$ is the least extender overlapping $\delta$, so 
$\alpha=\lambda$. Let $\kappa=\crit(E)$. Then $\kappa$ is a measurable limit 
of Woodins and strong cutpoints of $M(\Vv)$. Let $\gamma$ be least such that 
$\kappa<\lh(E^\Vv_\gamma)$. Then for all $\beta<\gamma$, $\lh(E^\Vv_\beta)<\kappa$.
Let $\Q=M^{*\Yy}_{\lambda+1}$, or $\Q=M^{*\Zz}_{\lambda+1}$, according to whether $E$ is used in 
$\Yy$ or $\Zz$. Note that $\kappa$ is a cutpoint of $\Q$. But then $\Yy\rest[\gamma,\lh(\Yy))$ is 
and $\Zz\rest[\gamma,\lh(\Zz))$ are equivalent to above-$\kappa$, normal trees on $\Q$.
So if $\Q\pins M^{\Vv}_\gamma$ then we are done, and if $\Q=M^\Vv_\gamma$
then note that $[0,\gamma]_\Vv$ drops (as our hod premice are below $\AD_\RR+$``$\Theta$ is 
measurable''), so again we are done.
\end{proof}

\begin{claim}\label{clm:c_maximal} We have:
\begin{enumerate}[label=--]
 \item $\delta$ is Woodin in $M^\Vv_c$, so $\delta=\delta_{\alpha_{n+1}}^{M^\Vv_c}$, the largest 
Woodin of $M^\Vv_c$.
\item $\delta$ is a strong cutpoint of $\Q_b$.
\item $M^\Vv_c\pins\Q_b$.
\end{enumerate}
\end{claim}
\begin{proof} Neither $\Q_b$, nor $\Q_c$ if it exists, can have overlaps of $\delta$, since 
otherwise $M^\Vv_c$ has a measurable limit of Woodins, which implies $c$ drops, contradiction.
But if $\delta$ is not Woodin in $M^\Vv_c$ then as before, $\Q_b\neq\Q_c$, so comparison gives a 
contradiction.

So the comparison of $\Q_b$ with $M^\Vv_c$ is above $\delta$, and succeeds, and this easily gives 
that $M^\Vv_c\pins\Q_b$.
\end{proof}

Let $\tau:M^{\Vvvec_n^*}_\infty\to M^{\Wwvec_n^*}_\infty$ be the final map given by the hull 
embedding (by the minimality of $\Wwvec_n^*$ with respect to $\sigmavec_n^*$,
the final model of $\Vvvec_n^*$ does indeed correspond to the final model of $\Wwvec_n^*$).
Let $\Q'$ be the lift of $\Q=\Q_b$ under $\tau$, and let $\tau_\Q=\tau\rest\Q$.
Let $\delta'=\tau_\Q(\delta)$.
Let $\Xx'=(\Uuvec,\Wwvec_0,\ldots,\Wwvec_{n-1},\Wwvec_n^*)$.
Let $\alpha=\alpha_{n+1}-1$ (possibly $\alpha=-1$) and $\alpha'=\tau_\Q(\alpha)$.
Using \ref{dfn:Gamma-fullness*}(\ref{item:Gamma-fullness*_dropping}),
let $\Upsilon'$ be an above-$\delta'$, $\Sigma_{\Q'(\alpha'),\Xx'}$-$(\om,\om_1,\om_1+1)$-strategy 
for $\Q$, whose restriction to countable trees is in $\Gamma$.
Let $\Upsilon$ be the $\tau$-pullback of $\Upsilon'$.
Like in Claim \ref{clm:Psi_is_correct}, we then get:

\begin{claim}
 $\Upsilon$ is a $\Sigma_{\Q(\alpha),\Xx}$-strategy, and the restriction of $\Upsilon$ to countable 
trees is in $\Gamma$ \tu{(}where $\Q(-1)=\emptyset)$; note $\Q(\alpha)=\R(\alpha)$\tu{)}.
\end{claim}

But then $\Q_b\pins\Lp^{\Gamma,\Sigma_{\R(\alpha),\Xx}}$, which with Claim 
\ref{clm:c_maximal}, contradicts $\Gamma$-fullness$^*$ preservation for $\Sigma$,
completing the proof.
\end{proof}

The following lemma, related to \cite[\S2]{maxcore}, is due
to Steel. However, the standard proof seems to have a gap (in the proof of
Claim \ref{clm:E_on_seq} below).
A correct proof of what is essentially the lemma appeared in \cite[\S5]{mim}, but that
proof is somewhat buried in another context, so we give a proof here for convenience. We state 
and prove the 
lemma literally only for pure $L[\es]$-constructions, but
it is easy to adapt it to strategy mice and other variants.

\begin{lemma}[Stationarity of
\protect{$L\left[\es\right]$} constructions]\label{lem:back_stat}
Let $\gamma$ be an uncountable cardinal. Let $\P$ be a reasonable
$k$-sound premouse, $\Psi$ a $(k,\gamma+1)$-strategy for 
$\P$ and $\CC=\left<\N_\alpha\right>_{\alpha\leq\gamma}$ be
a fully backgrounded $L[\es]$-construction.

Suppose that
for each active $\N_{\alpha+1}=(\N_\alpha,E)$ there is an extender $E^*$ such
that \tu{(}a\tu{)} $\card(\P)<\crit(E^*)$, \tu{(}b\tu{)}
$F\rest\nu(E)\sub E^*$, \tu{(}c\tu{)} if $\P$ is
non-tame then $i_{E^*}(\Psi)\rest
V_\eta\sub\Psi$ where $\eta$ is the sup of all $\delta+1$ such that $\delta$
is Woodin in $\N_\alpha$.

Then there is $\xi\leq\gamma+1$ such that:
\begin{enumerate}[label=\tulp\arabic*\turp]
 \item\label{item:<xi_below_P} for each $\alpha<\xi$, we have
$\N_\alpha\ins\P'$ for some $\Psi$-iterate $\P'$ of $\P$, and
\item\label{item:N_xi_reaches_P}
if $\xi\leq\gamma$ then there is a tree $\Tt$ via $\Psi$, of successor length,
$\N_\xi=M^\Tt_\infty$ and $b^\Tt$ does not drop in
model.\end{enumerate}
\end{lemma}
\begin{proof}
It suffices to prove that if \ref{item:<xi_below_P} holds at $\xi$, but
\ref{item:N_xi_reaches_P} does not, then \ref{item:<xi_below_P} holds at $\xi+1$. This is easy in 
all cases except when
$\xi=\alpha+1$ and $\N_{\alpha+1}=(\N_\alpha,E)$
for some $E$, so suppose this is the case. Let $E^*$ be a background extender
for $E$ and let $j=i_{E^*}$.
Let $\Tt$ be the tree witnessing the lemma's conclusion for $\alpha$. We assume
that $\Tt$ has minimal possible length. We must show that $E$ is used in $\Tt$.
Let $\nu=\nu(E)$ and $\kappa=\crit(E)$. The main point is the following claim:

\begin{claim} There is $\beta<\lh(j(\Tt))$ such that $\nu\leq\nu(E^\Tt_\beta)$
and $E\rest\nu\sub E^\Tt_\beta$.\end{claim}
\begin{proof}
As in the proof
that comparison of premice terminates, we have
$M^{j(\Tt)}_\kappa=M^\Tt_\kappa$ and $\kappa<_{j(\Tt)}j(\kappa)$ and
$i^{j(\Tt)}_{\kappa,j(\kappa)}$ exists and
\begin{equation}\label{eqn:it_emb_j} i^\Tt_{\kappa,j(\kappa)}\rest
M^{\Tt}_\kappa=j\rest M^\Tt_\kappa. \end{equation}
So let $\beta+1<_\Tt j(\kappa)$ be such that $\pred^\Tt(\beta+1)=\kappa$. We
claim that $\beta$ works. For let
\[ k:\Ult(\N_\alpha,E)\to j(\N_\alpha) \]
 be the factor embedding. Then $\crit(k)\geq\nu(E)$, and if $E$ is type 2 then
$\crit(k)\geq\lh(E)$. So 
$\N_\alpha$,
$M^\Tt_\kappa$, $M^\Tt_\beta$ and
$M^{j(\Tt)}_{j(\kappa)}$ agree
below $(\kappa^+)^{N_\alpha}$. So $E^\Tt_\beta$ measures all sets measured by
$E$ and by line (\ref{eqn:it_emb_j}) we have that $E\rest\nu'\sub
E^\Tt_\beta\rest\nu'$, where $\nu'=\min(\nu,\nu(E^\Tt_\beta))$. Now if
$(\kappa^+)^{\N_\alpha}<(\kappa^+)^{M^\Tt_\kappa}$ then
$\crit(k)=(\kappa^+)^{\N_\alpha}$, so $E$ is type 1 and
$\nu=(\kappa^+)^{\N_\alpha}$, so we are done. So assume
$(\kappa^+)^{\N_\alpha}=(\kappa^+)^{M^\Tt_\kappa}$, and assume $\nu'<\nu$. Since
also $(\kappa^+)^{M^\Tt_\kappa}\leq\nu'$, the ISC applies to $E\rest\nu'$. So
$E\rest\nu'\in\N_\alpha$, although $E\rest\nu'\notin j(\N_\alpha)$. So $E$ is not
type 2. So $E$ is type 3, but then $\lh(E^\Tt_\beta)<\nu$, contradicting the
fact that $\N_\alpha||\nu=j(\N_\alpha)||\nu$.\end{proof}

\begin{claim}\label{clm:E_on_seq}
Either:
\begin{enumerate}[label=--]
 \item $E\in\es_+(M^{j(\Tt)}_\beta)$, or
 \item $\R\eqdef M^{j(\Tt)}_\beta|\nu$ is active with extender
$F$ and $E\in\es_+^{\Ult(\R,F)}$.
\end{enumerate}
\end{claim}
\begin{proof}
If $(\kappa^+)^{N_\alpha}=(\kappa^+)^{M^{j(\Tt)}_\beta}$ this is just by the ISC.
So suppose $(\kappa^+)^{N_\alpha}<(\kappa^+)^{M^{j(\Tt)}_\beta}$. Then $E$ is type
1, the normal measure derived from $E$ is a submeasure of the normal 
measure derived from $E^{j(\Tt)}_\beta$, and
$M^{j(\Tt)}_\beta||\nu=\N_\alpha||\nu$. Thus, we can use \cite[4.11, 4.12,
4.15]{mim} (as $\P$ is reasonable) given that if $\R$ is
active with a type 3 extender $F$ then
\begin{equation}\label{eqn:ult_by_F} \Ult(\R,F)||\lh(E) =
\N_\alpha. \end{equation}
So suppose $F=F^\R\neq\emptyset$. We have 
$\Tt\rest(\kappa+1)=j(\Tt)\rest\kappa+1$, and note that
$\Tt$ uses no extenders with index in the interval $(\kappa,\nu)$, as $E$ is type 1, and $j(\Tt)$
uses no extender with index in the interval
$(\kappa,(\kappa^+)^{M^\Tt_\kappa})$. So
$M^\Tt_\kappa|\nu=\R$, and since
$N_\alpha|\nu$ is passive, therefore $E^\Tt_\kappa=F$. But then
$\Tt$ uses no extender with index in the interval
$(\nu,\lh(E))$, and line (\ref{eqn:ult_by_F}) is true.
\end{proof}

Now let $\lambda$ be least such that $\lh(E^{j(\Tt)}_\lambda)\geq\lh(E)$, and
let $\xi$ be the largest limit ordinal such that $\xi\leq\lambda$.
By the following claim, we clearly have that $j(\Tt)\rest\lambda+1$ is via
$\Psi$, which completes the proof.

\begin{claim}
$j(\Tt)\rest\xi+1=\Tt\rest\xi+1$.
\end{claim}
\begin{proof}
We have $\N_\alpha\ins M^\Tt_\infty$ and $j(\N_\alpha)\ins M^{j(\Tt)}_\infty$.
Let $\chi$ be the largest cardinal of $\N_\alpha$ and $\eps$ be the
largest limit cardinal of $j(\N_\alpha)||\lh(E)$. Then
$\eps\leq\chi$ and
$\N_\alpha||(\eps^+)^{\N_\alpha}=j(\N_\alpha)||(\eps^+)^{\N_\alpha}$ (possibly 
$(\eps^+)^{\N_\alpha}<(\eps^+)^{j(\N_\alpha)}$) and
$\univ{\N_\alpha}\sub j(\N_\alpha)$. These things follow from ${<\om}$-condensation,
considering the factor embedding $k$. Now let $\delta=\delta(j(\Tt)\rest\xi)$;
it follows that $\delta\leq\eps$. So
$\N_\alpha|\delta=j(\N_\alpha)|\delta$, and it suffices to see that for each
$\xi'\leq\xi$, we have $[0,\xi')_{j(\Tt)}=[0,\xi')_\Tt$. We prove this by
induction on $\xi'$. So assume $\Tt\rest\xi'=j(\Tt)\rest\xi'$. We may assume
$\xi'\geq\kappa$, so $\delta'=\delta(\Tt\rest\xi')\geq\kappa$ also. Now if
$\N_\alpha\sats$``$\delta'$ is not Woodin'' then let $\Q\pins M^\Tt_{\xi'}$ be
the Q-structure for $\delta'$. Then $\Q\pins \N_\alpha$, so $\Q\pins j(\N_\alpha)$,
so $\Q\pins M^{j(\Tt)}_{\xi'}$. Therefore
$[0,\xi')_\Tt=[0,\xi')_{j(\Tt)}$, as required. So suppose
$\N_\alpha\sats$``$\delta'$ is Woodin''. Since $\kappa\leq\delta'<\lh(E)$, and
so by Claim \ref{clm:E_on_seq}, $\P$ is non-tame.
So by our hypothesis, $j(\Psi)\rest V_{\delta'+1}\sub\Psi$, so
$[0,\xi')_{j(\Tt)}=[0,\xi')_\Tt$.\end{proof}\let\qed\relax
\end{proof}

\begin{definition}\label{dfn:P,Sigma-bounded_hod_pair_con}
Let $(\P,\Sigma)$ be a hod pair, within scope, and $\kappa$ be such that $\Sigma$ is a hod
$(\om,\kappa,\kappa+1)$-strategy. Let 
$a\in\HC$ be such that 
$\P\in\J_1(\hat{a})$ and $\hat{a}$ is swo'd. Suppose that $\MFsharp=\M_1^{\Sigma,\#}(\hat{a})$ 
exists and is $\Sigma$-$\kappa$-naturally iterable. Let $\NNN$ be any non-dropping 
$\Lambda^{\Sigma,\kappa}_\MFsharp$-iterate of $\MFsharp$. Let $\delta=\delta^\NNN$ and 
$\Sigma_\P^{\NNN}=\Sigma\rest\NNN$. Let $\chi\leq\delta+1$.

A \dfnemph{$(\P,\Sigma)$-bounded hod pair construction of $\NNN$,
of length $\chi$}, is a sequence
\[ 
\DD=\left<(\CC_\beta,\Tt_\beta,\alpha_\beta,\Q_\beta,\R_\beta,\M_\beta,\Sigma_\beta)\right>_{
\beta<\chi} \]
with the following properties holding inside $\NNN$ for all $\beta<\chi$:
\begin{enumerate}[label=--]
 \item $\Tt_\beta$ is a terminally non-dropping, successor length, normal tree on $\P$ via 
$\Sigma^\NNN_\P$, and $\Q_\beta=M^{\Tt_\beta}_\infty$ and $\R_\beta=\Q_\beta(\beta)$.
\item $\Tt_\alpha\psub\Tt_\beta$ for $\alpha<\beta$.
\item $\Tt_0$ is based on $\P(0)|\delta_0^\P$.
\item If $\beta+1<\chi$ then $\Tt_{\beta+1}$ is based 
on $\Q_\beta(\beta+1)|\delta_{\beta+1}^{\Q_\beta}$,
and is above $\delta_\beta^{\Q_\beta}$.
\item If $\beta$ is a limit then 
$\Tt_\beta=\Tt^*_\beta\conc\Sigma(\Tt^*_\beta)$ where $\Tt^*_\beta=\lim_{\alpha<\beta}\Tt_\alpha$.
\item $\Sigma_\beta$ is the strategy for $\R_\beta$ which is the tail of $\Sigma^\NNN_\P$.
\item $\CC_0$ is the maximal $L[\es]$-construction\footnote{Here and below, all 
background extenders are required to come from $\es^\NNN$.} of $\NNN|\delta$.
\item If $\beta+1<\chi$ then $\CC_{\beta+1}$ is the maximal 
$L^{\Sigma_\beta}[\es](\R_\beta)$-construction of $\NNN|\delta$.
\item If $\beta$ is a limit, $\CC_\beta$ is the maximal 
$L^{\Sigma^*_\beta}[\es](\R^*_\beta)$-construction of $\NNN|\delta$, where 
$\Sigma^*_\beta=\oplus_{\alpha<\beta}\Sigma_\alpha$ and 
$\R^*_\beta=\oplus_{\alpha<\beta}\R_\alpha$.
\item $\alpha_\beta<\delta$ and $\R_\beta=\NNN_{\alpha_\beta}^{\CC_\beta}$.
\item For all $\alpha<\alpha_0$ there is a successor length normal tree $\Tt$ on 
$\P$, via $\Sigma^\NNN_\P$, based on $\P(0)$, such that $\N_\alpha^{\CC_0}\ins 
M^\Tt_\infty$, and either $b^\Tt$ drops or $\N_\alpha^{\CC_0}\pins M^\Tt_\infty(0)$.
\item If $\beta+1<\chi$ then for all $\alpha<\alpha_{\beta+1}$ there is a successor length 
normal 
tree $\Tt$ on $\Q_{\beta}$, based on 
$\Q_{\beta}(\beta+1)$, above $\R_\beta=\Q_\beta(\beta)$, with
$\Tt_\beta\conc\Tt$ via $\Sigma^\NNN_\P$, and such that 
$\N_\alpha^{\CC_{\beta+1}}\ins M^\Tt_\infty$,
and either $b^\Tt$ drops or $\N_\alpha^{\CC_{\beta+1}}\pins M^\Tt_\infty(\beta+1)$.
\item If $\beta$ is a limit then for all $\alpha<\alpha_\beta$, either 
$\N_\alpha^{\CC_\beta}\pins\R_\beta$, or there is a successor length normal tree $\Tt$ on 
$\R_\beta$, above $\delta_\beta^{\R_\beta}$, with $\Tt_\beta\conc\Tt$ via $\Sigma^\NNN_\P$,
such that $b^\Tt$ drops and $\N_\alpha^{\CC_\beta}\ins M^\Tt_\infty$.
\item $\MMM_\beta$ is the least $\MMM\pins\NNN$ such that $\OR(\MMM)$ is a successor cardinal and 
$\beta,\alpha_\gamma<\OR(\MMM)$ for all $\gamma\leq\beta$. Let 
$\Lambda_{\MMM_\beta}$
be the $(\om,\Ord,\Ord)$-maximal strategy for $\MMM_\beta$, guided by Q-structures computed from 
ordinals and $\Sigma^\NNN_\P$. Then $\Sigma_\beta$ is exactly the strategy for $\R_\beta$ 
induced by lifting to $\Lambda_{\MMM_\beta}$.
\end{enumerate}

We say that such a construction is \dfnemph{successful} iff $\chi=\beta+1<\delta$
and $\R_\beta=\Q_\beta$ (thus, the construction has produced a non-dropping 
normal $\Sigma$-iterate $\Q_\beta=\N_{\alpha_\beta}^{\CC_\beta}$ of $\P$).
\end{definition}

\begin{lemma}\label{lem:P,Sigma_bounded_successful}
Adopt the hypotheses and notation of \ref{dfn:P,Sigma-bounded_hod_pair_con}. Then there is a 
unique successful
$(\P,\Sigma)$-bounded hod pair construction $\DD$ of $\NNN$.

Moreover, let $\beta<\lh(\DD)$ and $\Lambda^V_{\MMM_\beta}$ be the Q-structure guided 
$(\om,\kappa,\kappa+1)$-strategy for $\MMM_\beta$ \tu{(}so $\Lambda^V_{\MMM_\beta}$ is induced by 
the 
tail of $\Lambda^{\Sigma,\kappa,\kappa+1}_\MFsharp$ and 
$\Lambda_{\MMM_\beta}\sub\Lambda^V_{\MMM_\beta}$\tu{)}. 
Let $\Sigma^V_\beta$ be the hod $(\om,\kappa,\kappa+1)$-strategy for $\R_\beta$ induced by the tail 
of 
$\Sigma$ \tu{(}so $\Sigma_\beta\sub\Sigma^V_\beta$\tu{)}. Let $\Gamma^V_\beta$ be the
hod $(\om,\kappa,\kappa+1)$-strategy for $\R_\beta$ given by lifting to $\Lambda^V_{\MMM_\beta}$.
Then $\Sigma^V_\beta=\Gamma^V_\beta$.
\end{lemma}
\begin{proof}
This is partly proven in \cite{hod_mice}, but we cover some details
not presented there; it is in these details that the distinction between
hod strategies and full strategies is important.

It is easy to see that for each $\chi$, there is at most one construction of length $\chi$.
Trivially, if $\chi=0$, or $\chi$ is a limit and for 
all $\beta<\chi$, there is a construction of length $\beta$, then there is a construction of length 
$\chi$. So suppose there is an unsuccessful construction of length $\chi$; we need to see there is 
a construction of length $\chi+1$.

We assume $\chi=\beta+1$, as if $\chi=0$ or $\chi$ is a limit it is an easy variant.

Let $\CC$ be the maximal $L^{\Sigma_\beta}[\es](\R_\beta)$-construction of $\NNN|\delta$.
Let $\Psi$ be the above-$\R_\beta$, normal strategy for $\Q_\beta(\beta+1)$ given by continuing
$\Tt_\beta$ as a normal tree, using $\Sigma$. An easy variant of \ref{lem:back_stat}, together with
universality at $\delta$, \cite[Lemma 11.1]{DMATM}, shows that $\CC$ reaches a 
non-dropping 
$\Psi$-iterate $\R_{\beta+1}=\N_\alpha^\CC$ of $\Q_\beta(\beta+1)$, for some $\alpha<\delta$ such 
that for all $\xi<\alpha$, $\N_\xi^{\CC}$ is either a dropping such iterate, or a proper segment of 
such an iterate.
(With regard to universality, we don't need to iterate $\N_\delta^\CC$ in $\MFsharp$, so we don't 
need $\MFsharp$ to know any of its own iteration strategy.) Let $\Tt_{\beta+1}$ be the 
corresponding tree on $\P$.

So a length $\chi+1$ construction will exist given that $\Sigma_{\beta+1}$ 
agrees 
with the hod strategy $\Gamma$ for $\R_{\beta+1}$ given by lifting to $\Lambda_{\M_{\beta+1}}$
(where notation is as in \ref{dfn:P,Sigma-bounded_hod_pair_con}).
This follows from the ``moreover'' clause of the lemma at $\beta+1$, which we now prove. Let $\Uu$ 
be a limit length tree via both 
$\Gamma^V_{\beta+1}$ and $\Sigma^V_{\beta+1}$ (notation as in the statement of the lemma). Let 
$b=\Gamma^V_{\beta+1}(\Uu)$ and $c=\Sigma^V_{\beta+1}(\Uu)$. Because
$\Sigma$ and $\Lambda^{\Sigma,\kappa}_\MFsharp$ have hull condensation,
by taking a hull here we may assume everything is countable.

Sargsyan's argument showing that if $b$ does not drop then $b=c$
(using branch condensation for $\Sigma$) goes through here (cf. \cite[Lemma 2.15]{hod_mice}).
So assume that $b$ drops. Then because we are dealing with hod strategies,
$\Uu$ has the form $\Vvvec\conc\Vv$, where $\Vvvec$ does not drop,
$\Vv$ is normal and $\Vv\conc b$ drops.

Let $\gamma<\lh(\Vv)$ and $\alpha$ be such that $[0,\gamma]_\Vv$ does not drop,
and letting $\N=M^\Vv_\gamma$, such that $\alpha\leq\lambda^\N$ and for all $\tau<\gamma$, we have
\[ 
\lh(E^\Vv_\tau)<\eps\eqdef\bigcup_{\xi<\alpha}
\OR(\N(\xi))<\lh(E^\Vv_\gamma), \]
and $\Vv\rest[\gamma,\lh(\Vv))$
is based on $\N(\alpha)$ (and is above $\eps$). Let $\Omega$ be the above-$\eps$, hod 
$(\om,\kappa,\kappa+1)$-strategy for 
$\N(\alpha)$, given by normally extending $\Vv\rest\gamma+1$, continuing to use 
$\Gamma^V_{\beta+1}$. Let
\[ \Xx=\Tt_{\beta+1}\conc\Vvvec\conc(\Vv\rest\gamma+1) \]
and $\Sigma'=\oplus_{\beta<\alpha}\Sigma_{\N(\beta),\Xx}$. Let $\Upsilon$ be the above-$\eps$, hod 
$(\om,\kappa,\kappa+1)$-strategy for 
$\N(\alpha)$, given by normally extending $\Vv\rest\gamma+1$, continuing to use $\Sigma$.
So $\Upsilon$ is a $\Sigma'$-strategy as $(\P,\Sigma)$ is a hod 
pair. But $\Omega$ is also a $\Sigma'$-strategy, because $\Sigma$ has factor-hull condensation by 
\ref{lem:factor-hull_condensation}.

Let $\widetilde{\Vv}$ be the tree on $\N(\alpha)$ which is equivalent to 
$\Vv\rest[\gamma,\lh(\Vv))$ (the latter tree is on $\N$). Let $\tilde{b},\tilde{c}$ be the branches 
determined by $b,c$. By the previous paragraph we can use 
$\Omega$ and $\Upsilon$ to compare the phalanxes $\Phi(\widetilde{\Vv}\conc\tilde{b})$ and 
$\Phi(\widetilde{\Vv}\conc\tilde{c})$. This leads to contradiction almost as in the proof of 
\ref{lem:factor-hull_condensation}. The only slight difference is in showing that $\Q'$ has an 
iteration strategy in $\Gamma$ when $\tilde{c}$ does not drop, where $\Q'$ is as 
in the proof of \ref{lem:factor-hull_condensation}, so consider this. As before, $\delta'$ is a 
cutpoint of $\Q'$. We have a normal tree $\Yy$ via $\Sigma$ of successor length, such that $\Q'\ins 
M^\Yy_\infty$. If $b^\Yy$ drops then we can argue as in \ref{lem:factor-hull_condensation}, so 
suppose $b^\Yy$ does not drop.
Then $\Q'\pins M^\Yy_\infty$. If $\delta'$ is a cutpoint of $M^\Yy_\infty$ then we can use 
\ref{dfn:Gamma-fullness*}(\ref{item:Gamma-fullness*_non-dropping}), so suppose otherwise.
Let $E\in\es^{M^\Yy_\infty}$ be the extender of least index overlapping $\delta'$.
So $\OR(\Q')<\lh(E)$. Consider the tree $\Zz$ on $M^\Yy_\infty$, using only $E$.
So $\Q'\pins M^\Zz_1$, and note that \ref{dfn:Gamma-fullness*} applies to the stack $(\Yy,\Zz)$ and 
$\Q',\delta'$.
\end{proof}

The next lemma, and much of its proof, are similar to Sargsyan's \cite[Lemma 
3.35]{hod_mice}.

\begin{lemma}
\label{GenericInt}
Let $\kappa$ be an uncountable cardinal. Let $(\P,\Sigma)$ be such that $\P$ is
countable and reasonable and 
either
\begin{enumerate}[label=\tu{(}\roman*\tu{)}]
\item\label{item:P_std_premouse} $\P$ is an $n$-sound premouse and $\Sigma$ is the unique 
$(n,\kappa)$-strategy $\Sigma'$ for $\P$ such that if $\cof(\kappa)>\om$ then $\Sigma'$ extends to 
an $(n,\kappa+1)$-strategy, or \item\label{item:P,Sigma_hod_pair} $(\P,\Sigma)$ is a
hod pair, within scope, and $\Sigma$ is a hod $(\om,\kappa,\kappa+1)$-strategy.
\end{enumerate}
\tu{(}So in case \ref{item:P,Sigma_hod_pair} $\kappa$ is regular.\tu{)} Let 
$a\in\HC$ be such that 
$\P\in\J_1(\hat{a})$ and $\hat{a}$ is swo'd. Suppose that $\M_1^{\Sigma,\#}(\hat{a})$ 
exists and is $\Sigma$-$\kappa$-naturally iterable. Then $(\Sigma,a)$ is nice.\end{lemma}

\begin{proof}
$\Sigma$ has hull condensation, by the uniqueness of 
$\Sigma$ in case \ref{item:P_std_premouse}, and because $(\P,\Sigma)$ is within scope in case 
\ref{item:P,Sigma_hod_pair}.
\footnote{In case (i), we use the fact here that $\Sigma$ is only an $(n,\kappa)$-strategy. 
If $\kappa$ is singular then it seems difficult to deal with trees of length $(\kappa+1)$.}
It remains to see that $t_{(\Sigma,a)}$
determines itself on generic extensions.

We describe a process 
by which $\NNN[g]$ can compute
$\Sigma\rest\NNN[g]$ whenever $\NNN$ is a non-dropping $\Lambda^{X,\kappa}_\MFsharp$-iterate
of $\MFsharp=\M_1^{\Sigma,\#}(\hat{a})$ and $g$ is set-generic over $\NNN$. The result will then 
be a straightforward
corollary. So fix $\NNN$ and let
$\delta=\delta^\NNN$. Let $\Xx$ be the tree on $\MFsharp$ whose last model is $\NNN$.

Consider case (i).
If $\cof(\kappa)=\om$ let $\tau=\kappa$; otherwise let 
$\tau=\kappa+1$.
Let $\Lambda_\MFsharp$ be the $\Sigma$-$(0,{<\om},\tau)$-maximal strategy for $\MFsharp$ given by 
\ref{lem:Sigma_N_condensation}. So $\NNN$ is a $\Lambda_\MFsharp$-iterate. 
Let $\Lambda_\NNN$ be the $\Sigma$-$(0,{<\om},\tau)$-maximal strategy for 
$\NNN$ which is the tail of $\Lambda_\MFsharp$.
Let $\CC=\left<\N_\alpha\right>_{\alpha\leq\delta}$ be the
maximal $L[\es]$-construction of $\NNN|\delta$, where background extenders are
required to be in $\es^\NNN$.
Note that the hypotheses of \ref{lem:back_stat} hold
in $\NNN$ with respect to $\P,\gamma=\delta,\Sigma\rest\NNN,\CC$.

Now there is
$\alpha<\delta$ such that \ref{lem:back_stat}(ii) attains. For in
$\NNN$, $\delta$ is Woodin, and $\P$ is super-small, so we can apply the universality of 
$\N_\delta$ 
(see \cite[Lemma 11.1]{DMATM}). Note
that $\alpha<\mu$ where $\mu$ is the least strong of $\NNN$. Let $\gamma$ be a cutpoint of $\NNN$
such that $\alpha<\gamma<\mu$, and let $\theta=(\gamma^{++})^\NNN$.
Then via copying/resurrection, $\N_\alpha$, and therefore also $\P$, are normally iterable in $V$ 
via
lifting to nowhere-dropping normal trees on $\NNN$, via $\Lambda_\NNN$, based on $\NNN|\theta$. 
Let $\Sigma_\P$ be the resulting strategy for $\P$.
By the uniqueness of $\Sigma$ we have $\Sigma_\P=\Sigma$.
Note that $\theta\in\rg(i^\Xx)$.

Now consider case (ii). So $\kappa$ is regular. Let $\Lambda_\MFsharp$ be the 
$\Sigma$-$(0,\kappa,\kappa+1)$-maximal strategy for $\MFsharp$ given by 
\ref{lem:Sigma_N_condensation}. Let $\Lambda_\NNN$ be the $\Sigma$-$(0,\kappa,\kappa+1)$-maximal 
strategy for $\NNN$ which is the tail of $\Lambda_\MFsharp$. Let $\DD$ be the
$(\P,\Sigma)$-bounded hod pair construction of $\NNN$. By \ref{lem:P,Sigma_bounded_successful}, we 
have $\alpha<\delta$ and a normal
tree $\Tt$ via $\Sigma$ with last model $\R$ such that
$b^\Tt$ does not drop, $\DD$ has length $\beta+1$ and $\R=\R_\beta^\DD$, and 
$\Lambda_\beta^V=\Sigma_{\R,\Tt}$ (where $\Lambda_\beta^V$ is as in 
\ref{lem:P,Sigma_bounded_successful}; so this is just the strategy for $\R$ which lifts to 
$\Lambda^V_{\MMM_\beta^\DD}$). By
hull condensation,
$\Sigma$ has pullback consistency, so $\Sigma=\Sigma_\P$, where $\Sigma_\P$ is the pullback of 
$\Lambda^V_\beta$.
Note that $\OR(\MMM_\beta^\DD)<\mu$ where $\mu$ is the least strong of $\NNN$.
Let $\gamma$ be a cutpoint of $\NNN$ such that $\OR(\MMM_\beta^\DD)<\gamma<\mu$
and let $\theta=(\gamma^{++})^\NNN$.
And $\Sigma_\P$ is again computed by lifting to nowhere-dropping trees on $\NNN$, based on 
$\NNN|\theta$ (this time stacks of normal such trees). Again $\theta\in\rg(i^\Xx)$.

We now continue with both cases. It suffices to see that 
$\Lambda_\NNN\rest X$ is sufficiently definable over
$\NNN[g]$, where $X$ is the class of trees $\Tt\in\NNN[g]$ such that $\Tt$ is based
on $\NNN|\theta$ and is nowhere-dropping.
Iterating $\NNN$ for $\NNN|\theta$-based trees just requires computing the correct
Q-structures, which requires sufficient ordinals and knowledge of $\Sigma$. But
we don't yet know that $\Sigma``\NNN[g]\sub\NNN[g]$. We will computed the Q-structures by reducing
such trees $\Tt$ to trees in $\NNN$.

Let $\PP,\Ttdot\in\NNN|\crit(F^\NNN)$ be a poset and a $\PP$-name 
such that $\PP$ forces that $\Ttdot$ is a nowhere dropping, $\NNN|\theta$-based tree 
on $\NNN$, of
limit length, via the strategy to be described; it will follow that 
$\Ttdot^g$ is a correct tree on $\NNN$ for any $\NNN$-generic $g\sub\PP$.

\begin{claim}
Let $g$ be $\PP$-generic over $\NNN$. Let $Q=Q(\Ttdot^g)$. Then $Q\in\NNN[g]$.

In fact, let $\lambda$ be the
maximum of $\delta$, $(\lh(\Ttdot^g)^{++})^{\NNN[g]}$, and $(\car(\PP)^{++})^\NNN$.
Then there is a short tree $\Vv\in\NNN|\lambda$, $\Vv$ on $\NNN$, according to
$\Lambda_\NNN$, of successor length, such that for some $\alpha<\crit(F^\NNN)$, if
$G$ is $\Coll(\om,\lambda)$ generic over $\NNN[g]$, then in $\NNN[g][G]$, there is
an spm $Q$ which is a Q-structure for $\M(\Ttdot^g)$, and a
$\Sigma_1$-elementary embedding $\pi:Q\to M^\Vv_\infty|\alpha$. So $Q$ is unique with 
these properties and
$Q(\Ttdot^g)=Q\in\NNN[g]$.
\end{claim}

\begin{proof}
Suppose not and assume that $\PP$ forces the failure.
In $\NNN$, we first form a Boolean valued comparison of $M(\Ttdot)$ with $\NNN$,
forming a $\PP$-name for a tree $\Uudot$ on $M(\Ttdot)$ and a tree $\Vv$ on
$\NNN$. Note that $\NNN$ correctly computes
Q-structures as far as they exist during this comparison. Consider a
limit stage $(\Vv,\Uudot)\rest\lambda$ of the comparison. If a condition $q$
forces that $\Uudot\rest\lambda$ is eventually only padding then below $q$,
nothing need be done for $\Uudot$ at stage $\lambda$. Now suppose $q$ forces
otherwise. 
Suppose $p\leq q$ forces that here is a cofinal branch $b$ of $\Uudot$ such that
$Q(M(\Vv\rest\lambda))\ins
M^{\Uudot}_b$. Then below $p$, we set $[0,\lambda]_\Uudot=b$. If $p\leq q$ 
forces
otherwise, then below $p$, we declare that $\Uudot$ is \dfnemph{uncontinuable}, 
and terminate the comparison. (In the latter case 
$p$ forces that $\Uudot$ has limit length; we deal with this later.) For each 
stage
$\alpha$ of the comparison, let $\lh_\alpha$ be the index of any extender
(forced by some $p$ to be) used at that stage. For limit $\lambda$, let
$M((\Vv,\Uudot)\rest\lambda)$
be the lined up part of that stage, of height $\sup_{\alpha<\lambda}\lh_\alpha$.

\begin{subclaim} We have:
\begin{enumerate}[label=\tu{(}\alph*\tu{)}]
 \item $\Vv$ is based on $\NNN|\theta$;
 \item if $\alpha$ is such
that $[0,\alpha]_\Vv$ does not drop and $\PP$ forces that
$M^{\Uudot}_\alpha|\theta'=M^\Vv_\alpha|\theta'$, where
$\theta'=i^\Vv_{0,\alpha}(\theta)$, then the comparison terminates at stage
$\alpha$, and in fact, $\PP$ forces that $M^{\Uudot}_\alpha\ins
M^\Vv_\alpha|\theta'$;
 \item at every limit stage $\lambda$, a Q-structure for
$M((\Vv,\Uudot)\rest\lambda)$ exists;
 \item the comparison terminates \tu{(}i.e. there is $\alpha$ such that $\PP$ 
forces that either $\Uudot$ is uncontinuable, or $M^\Vv_\alpha\ins 
M^\Uudot_\alpha$, or $M^\Uudot_\alpha\ins M^\Vv_\alpha$\tu{)};
 \item there is $p\in\PP$ forcing
that if $\Uudot$ has a final model, then 
$M^\Uudot_\infty\pins M^\Vv_\infty$.
\end{enumerate}
\end{subclaim}
\begin{proof}
Part (b) implies (a) and (c). Suppose (b) fails. Let $\alpha$ be the least
failure, and let $p$ be a condition forcing this failure. Let $g\sub\PP$ be
generic with $p\in g$. Let $\Tt'$ be the tree on $\NNN$ which uses the same
extenders as does $\Tt=\Ttdot^g$, followed by $\Lambda_\NNN(\Tt)$, and let $W_0=M^{\Tt'}_\infty$. 
So $b^{\Tt'}$ is non-dropping (as $\Tt$ was nowhere dropping). Let $\Uu'$ be
the tree on $W_0$ using the same extenders as $\Uu^g$. Let $W=M^{\Uu'}_\alpha$.
So $\theta'<\OR(W)$. We can compare $(M^\Vv_\alpha,W)$, producing trees
$(\Tt_1,\Tt_2)$. The comparison begins above $\theta'$, a cardinal of
$M^\Vv_\alpha$. Note that by choice of $\theta$, all extenders used in the comparison have 
critical point $>\theta'$. Suppose $b^{\Uu'}$ drops. Then $\rho_{n+1}^\W<\theta'$,
where $n=\deg^{\Uu'}(\alpha)$. Also then,
$b^{\Tt_1}$ drops, whereas $b^{\Tt_2}$ does not, and $\Tt_1,\Tt_2$
have the same last model. But the last model $Z$ of $\Tt_1$ has
$\rho_\om(Z)\geq\theta'$, contradiction. So $b^{\Uu'}$ does not drop, and so 
neither
do $b^{\Tt_1},b^{\Tt_2}$, and
$j=k$ where $j=i^{(\Xx,\Vv,\Tt_1)}$ and $k=i^{(\Xx,\Tt',\Uu',\Tt_2)}$. But
$j(\theta)=\theta'$ and $k(\theta)>\theta'$, contradiction. This gives (b).

The usual proof that boolean-valued comparisons terminate gives (d).

So if (e) fails, then $b^\Vv$ drops, so $M^\Vv_\infty$ is unsound, and $\PP$ forces 
that $M^\Uudot_\infty=M^\Vv_\infty$. But then again the usual methods yield a 
contradiction.
\end{proof}

Now let $p$ be as in part (e), and let $g\sub\PP$ be $\NNN$-generic, with 
$p\in g$. Let $\Tt=\Ttdot^g$ and $\Uu=\Uudot^g$. Let $Q=Q(M(\Tt))$. Let 
$W_0,\Uu'$ be as before, and let $\Uu_Q$ be 
the $0$-maximal tree on $Q$ given by
$\Uu$ (with the same extenders and branches).

Suppose that $\Uu$ 
has a last model $R$. So we have $R\pins M^\Vv_\infty$ and $b^{\Uu}$ does not drop, and 
so neither do $b^{\Uu'}$ or $b^{\Uu_Q}$. Let $\pi:M^{\Uu_Q}_\infty\to 
i^{\Uu'}(Q)$ be the factor map. Then $\pi$ is a weak $0$-embedding. So by 
\ref{StrategyCondensation}, $M^{\Uu_Q}_\infty$ is a $\Sigma$-premouse. Also, 
$i^{\Uu_Q}:Q\to M^{\Uu_Q}_\infty$ is
continuous at $\delta=\delta(\Ttdot^g)$, and $M^{\Uu_Q}_\infty$ has no $E$-active 
levels above $i^{\Uu_Q}(\delta)$ and $i^{\Uu_Q}(\delta)$ is Woodin in $M^{\Uu_Q}_\infty$. It 
follows that 
$M^{\Uu_Q}_\infty\ins M^\Vv_\infty$. Also, $i^{\Uu_Q}$ is $\Sigma_1$-elementary. 
So $Q$, $\Vv$, $M^{\Uu_Q}_\infty$ and $i^{\Uu_Q}$ witness the truth of the claim, 
a contradiction.\footnote{Ostensibly $M^{\Uu_Q}_\infty$ might be a strict 
segment of the Q-structure for $M^\Vv_\infty|i^{\Uu_Q}(\delta)$, but this is not 
relevant. If one chooses $n<\om$ appropriately, and takes $\Uu_Q$ to be 
$n$-maximal instead of $0$-maximal, then one can arrange that $M^{\Uu_Q}_\infty$ 
is the Q-structure.}

Suppose now that $\Uu$ is uncontinuable, so has limit length. Let 
$b$ be the $\Uu$-cofinal branch determined by $\Lambda_\NNN$. Note that $b$ does not drop, and 
$M(\Uu)=M^\Uu_\infty$.
But this leads to the same contradiction as in the previous paragraph.
\end{proof}

This completes the proof that $\NNN[g]$ computes $\Sigma\rest\NNN[g]$. Now let
$\Phi$ be the formula ``There is no largest cardinal, there is a Woodin
cardinal $\delta$, in case (i) the $L[\es]$-construction reaches a non-dropping $\Sigma$-iterate 
of $\P$, and in case (ii) the $(\P,\Sigma)$-bounded hod pair
construction is successful at some stage $<\delta$, and every partial order $\PP$ forces
that the process described above always succeeds''. Let $\Psi$ be
the formula defining $\Sigma\rest\NNN[g]$ through the above process. Note that if
$\NNN'\ins\NNN$ and $\NNN'\models\Phi$ and $g$ is set generic over $\NNN'$, then
$\NNN'[g]$ is indeed closed under $\Sigma$, and $\Sigma\rest\NNN'[g]$ is defined
over $\NNN'[g]$ by $\Psi$. So $(\Phi,\Psi)$ generically determines $t_{(\Sigma,a)}$,
as required. (We don't actually need that the Woodin of $\NNN'$ is a cardinal of
$\NNN$.)
\end{proof}

\begin{notation}\label{ntn:fixed_notation}
Let $(\Omega,A_0)$ be nice, $t_0=t_{\Omega,A_0}$ and $\kappa_0=\kappa_t$. Let 
$\MFsharp=\M_1^{\Omega,\#}(A_0)$ and $\Lambda_\MFsharp=\Lambda^{\Omega,\kappa_0}_\MFsharp$. Let 
$(\Phi_0,\Psi_0)$ be a pair that generically determines $(\Omega,A_0)$. Let $a_0\in\RR$ code 
$A_0$ in a canonical way.\footnote{If $a_0$ can be chosen such that $\MFsharp$ codes $a_0$
then we do so, and $a_0$ is redundant.} These 
objects are fixed for 
the remainder of the paper.
\end{notation}

\begin{definition}
An hpm $\N$ is \dfnemph{$\MFsharp$-like}\footnote{The ``$\MFsharp$'' in ``$\MFsharp$-like''
is just a symbol; it does not refer to the fixed structure $\MFsharp$.} iff $\N$ 
is non-$1$-small, all proper segments of $\N$ are $1$-small, and 
$\exists\gamma\in(\delta^\N,\l(\N))$ such that $\N|\gamma\sats\Phi_0$.
\end{definition}

\begin{remark}\label{rem:generic_tree}
G-organization will use an initial segment of the tree for making a structure \emph{generically 
generic}, due to Sargsyan \cite{hod_mice}. We recall this notion 
and define some related 
notation and terminology now.

Let $\N,\P$ be transitive structures, where $\P$ is $\MFsharp$-like. Let $\QQ=\Coll(\om,\N)$. Let 
$\dot{x}_\N$ be the canonical $\QQ$-name for the real coding $\N$ determined by a $\QQ$-generic 
filter. Let $\Tt$ be a normal iteration tree on $\P$. We say that $\Tt$ is \dfnemph{making $\N$
generically generic} iff:
\begin{enumerate}[label=--]
\item $\Tt\rest\OR(\N)+1$ is a linear iteration at the least measurable of $\P$.
\item Suppose $\lh(\Tt)\geq\OR(\N)+2$ and let $\alpha+1\in(\OR(\N),\lh(\Tt))$.
Let $\delta=\delta(M^\Tt_\alpha)$ and let $\BB=\BB(M^\Tt_\alpha)$. Then 
$E^\mathcal{T}_\alpha$ is the extender $E\in\es_+(M^\Tt_\alpha)$ with
least index such that some
$p\in\QQ$ forces ``There is a $\BB$-axiom induced
by $E$ which fails for $\dot{x}_\N$''.
\end{enumerate}

Given a putative strategy $\Sigma$ for $\P$,
let $\Tt^{*\Sigma}_\N$ denote the longest putative tree $\Tt$ via $\Sigma$ which is making $\N$ 
generically generic. Clearly if $\Sigma$ is a normal $\kappa$-strategy for a large enough $\kappa$
then $\Tt\eqdef\Tt_\N^{*\Sigma}$ has 
successor length and $\QQ$ forces that $\dot{x}_\N$ is generic for 
$\BB(M_\infty^\Tt)$.

Let $\Tt^*_\N$ denote $\Tt^{*\Lambda_\MFsharp}_\N$.
\end{remark}

Sargsyan noticed (see \cite[Definition 3.37]{hod_mice}) that one can feed 
$\Omega$ into a strategy mouse $\N$ indirectly, by feeding in the branches for 
something like $\Tt^*_\M$, for various $\M\ins\N$. The key notion of
\emph{$\g$-organized $\Omega$-premice}, to come,
uses this idea, and the main point of it is due to Sargsyan.
We will only actually use a certain initial segment $\Tt^\Sigma_\N$ 
of $\Tt^{*\Sigma}_\N$:

\begin{definition}\label{dfn:P^Phi} Let $\P$ be $\MFsharp$-like.
Then $\P_{\Phi_0}$ denotes the least $\P'\pins\P$
such that for some cardinal $\delta'$ of $\P$, $\P'\sats\Phi_0+$``$\delta'$ is Woodin''.
Note that $\P_{\Phi_0}$ is a strong cutpoint of $\P$.
Given a transitive structure $\N$ and a putative strategy $\Sigma$ for $\P$,
$\Tt^\Sigma_\N$ denotes the initial segment of $\Tt^{*\Sigma}_\N$ based on $\P_{\Phi_0}$.
Let $\Tt_\N$ denote $\Tt^{\Lambda_\MFsharp}_\N$.
\end{definition}

To ensure the absoluteness of iterations making structures generically generic, 
we will require our models to add branches to iteration trees sufficiently slowly:

\begin{definition}
Let $\M$ be an aspm and $\eta<\OR(\M)$. Let $\Tt\in\M$ be a putative tree via $\Sigma^\M$.
Then $\Ss^\M_\Tt$ denotes the least $\Ss\ins\M$ such that $\Tt$ is via $\Sigma^\Ss$.
We say that $\Tt$ is \dfnemph{$\M$-reckonable above $\eta$} iff for every limit 
$\alpha\in[\eta,\lh(\Tt))$ we have the following.
Let $\zeta=\sup_{n<\om}\wfp(\OR(M^\Tt_{\alpha+n}))$.
Then:
\begin{enumerate}[label=--]
 \item if $\alpha+\om<\lh(\Tt)$ then 
$\OR(\Ss^\M_{\Tt\rest\alpha+1})+\zeta\leq\OR(\Ss^\M_{\Tt\rest\alpha+\om})$,
\item if $\lh(\Tt)\leq\alpha+\om$ then $\OR(\Ss^\M_{\Tt\rest\alpha+1})+\zeta\leq\OR(\M)$, and
\item if $\lh(\Tt)<\alpha+\om$ and $\M\sats$``$M^\Tt_\infty$ is 
wellfounded'' then $M^\Tt_\infty$ is 
wellfounded (equivalently, 
$M^\Tt_\infty\sats$``$\OR(\Ss^\M_{\Tt\rest\alpha+1})+\OR(M^\Tt_\infty)\leq\OR(V)$'').\qedhere
\end{enumerate}
\end{definition}

\begin{remark}\label{rem:absoluteness_generically_generic}
Let $\M$ be an aspm such that $\cp^\M$ is $\MFsharp$-like.
Let $\N\pins\M$ satisfy $\ZF$.
Let $\Tt\in\M$ be a putative tree via $\Sigma^\M$ (on $\cp^\M$), based on 
$\cp^\M_{\Phi_0}$, such that $\Tt$ is $\M$-reckonable above 
$\OR(\N)$. 
Then $\Tt$ is making $\N$ generically generic (in $V$)
iff $\M\sats$``$\Tt$ is making $\N$ 
generically generic''.
Moreover, let $\Uu'=\Tt^{\Sigma^\M}_\N$ (as computed in $V$)
and $\Uu=\Uu'\rest\lambda$
where $\lambda$ is largest such that
$\Uu\rest\alpha+1$ is $\M$-reckonable above $\OR(\N)$ for all $\alpha<\lambda$. 
Given $\alpha+1<\lh(\Uu')$ let $e_\alpha=E^{\Uu'}_\alpha$,
and given $\alpha+1=\lh(\Uu')$, if $M^\Uu_\alpha$ is illfounded then 
let $e_\alpha=0$,
and otherwise let $e_\alpha=1$.
Then the map $\alpha\mapsto(\Uu\rest\alpha+1,e_\alpha)$,
with domain $\lambda$, is $\rPi_2^\M(\Ll^-,\{\N\})$, 
uniformly in $\M,\N$.\footnote{There is a natural $\Sigma_1$ formula which attempts to define this 
function, which computes the correct values on the domain of the function, but might give a larger 
domain.} Further, suppose that $\Tt\eqdef\Tt^{\Sigma^\M}_{\N}$ exists, is in $\M$,
and is $\M$-reckonable above $\OR(\N)$.
Then $\{(\Tt,(M^\Tt_\infty)_{\Phi_0})\}$ is $\Sigma_1^\M(\Ll^-,\{\N\})$, uniformly in $\M,\N$.

These facts use the local definability of the
$\Coll(\om,\N)$ forcing relation. Given
$p\in\Coll(\om,\N)$, $n<\om$, a limit ordinal $\alpha<\lambda$ and
$E\in\es(M^\Tt_{\alpha+n})$ such that
$\nu(E)$ is inaccessible in $M^\Tt_{\alpha+n}$, the question of whether
$p\forces$``$E$
induces an extender algebra axiom not satisfied by $\dot{x}_\N$'' is 
computed over $\J^{\hpm}_{\nu(E)}(\Ss^\M_{\Tt\rest\alpha+1})$. (Such an 
axiom has the form
\[
\bigvee_{\gamma<\crit(E)}\varphi_\gamma\iff\bigvee_{\gamma<\nu(E)}\varphi_\gamma
,
\]
where for each $\gamma<\nu(E)$, $\varphi_\gamma\in M^\Tt_{\alpha+n}|\nu(E)$, so the
forcing relation below $p$ regarding the truth of $\varphi_\gamma$ is computed
over some proper segment of $\J^{\hpm}_{\nu(E)}(\Ss^\M_{\Tt\rest\alpha+1})$.)
\end{remark}

\begin{definition}
 Let $\R$ be an aspm such that $\cp^\R$ is $\MFsharp$-like.
 Let $\psi\in\Ll$.
 The \dfnemph{$(\g,\psi)$-hierarchy} of $\M$ is the pair
 $(\left<\M_\alpha\right>_{\alpha\leq\gamma},\left<\N_{\alpha+1}\right>_{\alpha<\gamma'})$ with
 $\gamma,\gamma'\in\Ord$ both as large as possible such that $\gamma\leq\gamma'\leq\gamma+1$ and:
 \begin{enumerate}
 \item $\N_{\alpha+1}\ins\R$ for each $\alpha<\gamma'$ and $\M_\alpha\ins\R$ for each 
$\alpha\leq\gamma$.
 \item $\M_0=\M|1$ and $\OR(\M_\lambda)=\lim_{\alpha<\lambda}\OR(\M_\alpha)$ for limit $\lambda$.
\item \label{item_2} For $\alpha<\gamma'$, $\N_{\alpha+1}$ is the least $\N\ins\R$ such 
that $\M_\alpha\pins\N$ and $\N\sats\ZF$.
\item For $\alpha<\gamma$, $\M_{\alpha+1}$ is the least $\M\ins\R$ such that $\N_{\alpha+1}\pins\M$
and for some $\Ss$ with $\N_{\alpha+1}\ins\Ss\ins\M$ we have either:
\begin{enumerate}[label=--]
\item $\Ss\sats\neg\psi$, or
\item $\Tt'\eqdef\Tt^{\Sigma^\Ss}_{\N_{\alpha+1}}$ exists, is in $\M$ and is $\M$-reckonable
above $\OR(\N_{\alpha+1})$.
\end{enumerate}
\end{enumerate}

For $\N\ins\M$,
we say that $\N$ is a \dfnemph{$(\g,\psi)$-tree activation 
level} of $\M$ iff $\N = \N_{\alpha+1}$ for some $\alpha$. 
We say that $\M$ is \dfnemph{$(\g,\psi)$-whole} iff $\M=\M_\gamma$,
and say that $\M$ is \dfnemph{$(\g,\psi)$-closed} iff $\M$ is $(\g,\psi)$-whole and $\gamma$ is a 
limit.

We abbreviate $(\g,\mathrm{true})$ with $\g$ (for example in
the \emph{$\g$-hierarchy} of $\R$, etc). We abbreviate $(\g,\text{``}\Theta\text{ 
exists''})$ with both $(\g,\Theta)$ and $\Theta$-$\g$.
\end{definition}

\begin{remark}
 Let $\M$ be an aspm such that $\cp^\M$ is $\MFsharp$-like, and $\psi\in\Ll$.
 Let the $(\g,\psi)$-hierarchy of $\M$ be 
$(\left<\M_\alpha\right>_{\alpha\leq\gamma},\left<\N_{\alpha+1}\right>_{\alpha<\gamma'})$.
Then $\left<\M_\alpha\right>_{\alpha\leq\gamma}\rest\M$ is $\Sigma_1^\M(\Ll^-)$,
and $\left<\N_{\alpha+1}\right>_{\alpha<\gamma'}\rest\M$ is $\Delta_2^\M(\Ll^-)$, uniformly in 
$\M$; this follows easily from \ref{rem:absoluteness_generically_generic}.\footnote{
If $\gamma'=\alpha+1$ and $\l(\M)=\gamma'+1$ then 
$\left<\N_{\alpha+1}\right>_{\alpha<\gamma'}$ is not $\Sigma_1^\M(\Ll^-)$ because 
$\N_{\gamma'}\sats\ZF$.}
Similarly, there is $\varrho_\psi\in\Ll$ such that $\M\sats\varrho_\psi$ iff
$\M$ is $(\g,\psi)$-whole, uniformly in $\M$.
\end{remark}

\newcommand{\aspm}{\mathrm{aspm}}
\begin{definition}
Let ``$V$ is an aspm'' be the natural formula $\psi\in\Ll$ such that for any transitive 
$\Ll$-structure $\M$, $\M\sats\psi$ iff $\M$ is an aspm.
\end{definition}

\begin{definition}[$\varphi_{(\g,\psi)}$]\label{dfn:varphi_g}
For $\psi\in\Ll$, $\varphi_{(\g,\psi)}$ denotes the $\Ll$-formula of one free variable 
$\Tt$ asserting (when interpreted over transitive $\Ll$-structures) ``$V$ is an aspm, $\cp$ is an 
$\MFsharp$-like hpm, the $(\g,\psi)$-hierachy of $V$ has the form
\[ (\left<\M_\alpha\right>_{\alpha\leq\gamma},\left<\N_{\alpha+1}\right>_{\alpha<\gamma+1}) \]
with $\N\eqdef\N_{\gamma+1}\pins V$, $\Tt$ is a limit length iteration tree via $\Sigma^V$ (on 
$\cp$),
based on $\cp_{\Phi_0}$,
making $\N$ generically generic,
$\Tt$ is $V$-reckonable above $\OR(\N)$,
and $\Sigma^V(\Tt)$ is undefined.''

We have $\varphi_\g=\varphi_{(\g,\mathrm{true})}$; let
$\varphi_\mathrm{G}=\varphi_{(\g,\Theta)}$.
\end{definition}

The notion \emph{$\g$-organized $\Omega$-premouse} below is a variant of 
Sargsyan's reorganized hybrid strategy premouse, \cite[Definition 3.37]{hod_mice}:

\begin{definition}\label{dfn:g_organized_F_premouse}
$\gOmega=(\Lambda_\MFsharp,\varphi_\g)$
and $\GOmega=(\Lambda_\MFsharp,\varphi_\rmG)$.
For example, a
$\gOmega$-premouse is a 
$(\Lambda_\MFsharp,\varphi_\g)$-premouse and
$\Lp^{\gOmega}(x)=\Lp^{(\Lambda_\MFsharp,\varphi_\g)}(x)$, etc.
A \dfnemph{$\g$-organized $\Omega$-premouse} is a $\gOmega$-premouse.
\end{definition}

So a $\g$-organized $\Omega$-pm is over $A$ for some $A\in\witri{V}$ where $\MFsharp\in\J_1(A)$.

\begin{lemma}\label{lem:gF_props}
The class of $\g$-organized $\Omega$-pms $\M$ such that $\hmPsi^\M=\emptyset$ is very condensing.
For any $\g$-organized $\Omega$-pm $\M$ not of type 3,
and any $\pi:\R\to\M$ a weak $0$-embedding, $\R$ is a $\g$-organized 
$\Omega$-pm.
\end{lemma}
\begin{proof}
These facts
follow from \ref{cor:Sigma-pm_very_condensing} and \ref{StrategyCondensation} 
respectively.\end{proof}

As in \cite[Lemma 3.38]{hod_mice}, the first consequence of $\g$-organization is the 
following. Because $t_{\Omega,A_0}$ determines itself on generic extensions,
$\g$-closure ensures closure under $\Omega$:

\begin{lemma}
\label{ClosedUnderJ}
Let $\M$ be a $\g$-closed $\g$-organized $\Omega$-pm. Then
$\M$ is closed under $\Omega$. In fact, for any set generic extension $\M[g]$ 
of $\M$, with $g\in V$\footnote{Without the assumption that $g\in V$, it seems that the 
domain of $\Omega\rest\M[g]$ might not be definable over $\M[g]$.}, $\M[g]$ is closed under
$\Omega$ and $\Omega\rest\M[g]$ is $\Ll^-$-definable over $\M[g]$, uniformly in 
$\M,g$.\end{lemma}
\begin{proof}[Proof sketch.]
We show that $\M$ is closed under $\Omega$; the generalization to generic 
extensions of $\M$ and the definability of $\Omega$ is 
similar. We assume that $\Omega$ is an operator; the strategy case is similar.

Let 
$z\in\univ{\M}\inter\dom(\Omega)$; we want to see that 
$\Omega(z)\in\univ{\M}$. Let $t=\Th_\om^{\Omega(z)}(z)$; it suffices to see that $t\in\M$. Let 
$\N,\N'\pins\M$ be tree activation levels of $\M$ with $z\in\N\pins\N'$. Then 
$\Tt\eqdef\Tt_\N\in\N'$. Let
$\MFsharp^*=\MFsharp^\Tt_{\Phi_0}$ and $\QQ=\Coll(\om,\N)$. Then in $\N'$, $\QQ$ forces that 
$\dot{x}_\N$ is extender 
algebra generic over $\MFsharp^*$. So by \ref{lem:determines_extend}, for
$w\in z^{<\om}$
and any formula $\varphi$, $\varphi(w)\in t$ iff in $\N'$,
$\QQ$ 
forces that $\check{\MFsharp^*}[\dot{x}_\N]\sats$``There is $y$ such that $\Psi_0(\check{z},y)$ 
and $y\sats\varphi(\check{w})$''.
\end{proof}

The analysis of scales in $\Lp^{^\g\Omega}(\RR)$ runs into some problems (see footnotes
\ref{ftn:motivate_Theta-g} and \ref{ftn:motivate_Theta-g_2}). So we will analyze scales in a 
slightly 
different hierarchy, which we now describe.

\begin{definition}\label{dfn:self-scaled}
Fix a natural coding of elements of $\HC$ by reals. Let $\Upsilon\sub\HC$. Given a set 
$\Upsilon\sub\HC$, 
$\Upsilon^\HCc$ denotes the set of codes for elements of $\Upsilon$ in this coding.\footnote{Note 
that for any $\J$-structure $\M$ such that $\HC^\M\in\M$, the decoding function (for 
the above 
codes), restricted to $\RR^\M$, is definable over $\HC^\M$, so 
$(\Upsilon\inter\HC^\M)^\HCc=\Upsilon^\HCc\inter\M$.} We say that 
$\Upsilon$ is \dfnemph{self-scaled} iff there are scales on $\Upsilon^\HCc$ and $\RR\cut 
\Upsilon^\HCc$ which are analytical in $\Upsilon^\HCc$ (i.e. $\Sigma^1_n(\Upsilon^\HCc)$ for some 
$n<\om$).
\end{definition}

\begin{definition}\label{dfn:Theta-g-organized}
An aspm $\M$ is \dfnemph{suitably based} iff
$\cp^\M\in\HC^\M$ is $\MFsharp$-like,
$\hmb^\M=\hat{x}$
where $x=(\HC^\M,\Upsilon)$ for some $\Upsilon\sub\HC^\M$ such that $\M\sats$``$\Upsilon$ is 
self-scaled'', and $\hmPsi^\M=\emptyset$.
Abusing terminology, we say that $\M$ is \dfnemph{over $\Upsilon$} and write
$\Upsilon^\M=\Upsilon$. Let $\M$ be a suitably based aspm over $\Upsilon$. Let 
$\vec{\leq}^\M,\vec{\leq'}^{\M}$ denote what are, in $\M$,
the least 
analytical-in-$\Upsilon^\HCc$ scales on $\Upsilon^\HCc,\RR\cut\Upsilon^\HCc$.
If there is some $\N\ins\M$ which is admissible, then working in $\M$ (or $\N$) let $U^\M,U'^{\M}$ 
denote the trees of these scales, respectively.

A \dfnemph{$\Theta$-$\g$-spm} is a suitably based
$\varphi_{\mathrm{G}}$-indexed spm.

A \dfnemph{$\Theta$-g-organized $\Omega$-premouse} is a $\Theta$-g-spm
which is a $(\Lambda_\MFsharp,\varphi_{\mathrm{G}})$-pm.
\end{definition}

In our application to core model induction, we will be most interested in the cases that either 
$\Upsilon^\M=\emptyset$ or $\Upsilon^\M=\Omega\rest\HC^\M$.

\begin{definition} Let ``$V$ is a $\Theta$-g-spm'' be the natural 
formula
$\psi\in\Ll$ such that for all transitive $\Ll$-structures $\M$,
$\M\sats\psi$ iff $\M$ is a $\Theta$-g-spm.
\end{definition}

\begin{definition}
 Let $\M$ be an aspm and let $\P\pins\J^{\hpm}(\P)\ins\M$ with $\P$ a strong cutpoint of $\M$.
 Then $\M\down\P$ denotes the aspm $\M'$ defined by induction on $\M$ as follows:
$\univ{\M'}=\univ{\M}$, $\hmb^{\M'}=\hat{\P}$, $\hmp^{\M'}=\hmp^\M$,
 $\Psi^{\M'}=\Sigma^\P$, $\hmP^{\M'}=\hmP^\M$, $\hmE^{\M'}=\hmE^\M$,
$\l(\P)+\l(\M')=\l(\M)$ and $\N\down\P\pins\M'$ for all $\N$ such that $\P\pins\N\pins\M$ (this 
determines $\hmPvec^{\M'}$).
\end{definition}

\begin{lemma}\label{lem:charac_Tg-org}
Let $\M$ be an hpm. Then the following are equivalent:
\tu{(}i\tu{)} $\M$ is a
$\Theta$-$\g$-organized $\Omega$-pm; \tu{(}ii\tu{)} 
$\M\sats$``$V$ is a $\Theta$-g-spm'' and $\hmp^\M=\MFsharp$ 
and $\Sigma_{\varphi_{\rm{G}}}^\M\sub\Lambda_\MFsharp$; \tu{(}iii\tu{)}
$\M$ is a suitably based aspm and $\hmp^\M=\MFsharp$ and for all $\N\ins\M$:
\begin{enumerate}[label=--]
 \item if $\P\pins\J^{\hpm}(\P)\ins\N$ and $\P\sats\ZF^-$ and every $\R$ such that 
$\P\ins\R\pins\N$ has $\Theta^\R=\OR(\P)$ \tu{(}possibly $\R=\P$; therefore $\P$ is a strong 
cutpoint of $\N$\tu{)} then
 $\N\down\P$ is a g-organized $\Omega$-pm, and
 \item if there are arbitrarily large $\R\pins\N$ satisfying ``$\Theta$ does not exist'' then $\N$ 
is passive.
\end{enumerate}
\end{lemma}

\begin{lemma}\label{lem:Th-g_very_con}
The class of $\Theta$-$\g$-organized 
$\Omega$-premice is very condensing.
\end{lemma}
\begin{proof}
By \ref{cor:Sigma-pm_very_condensing}.
\end{proof}

\begin{corollary}\label{cor:Tg_pullback} Let $\M$ be an $n$-sound $\Theta$-$\g$-organized 
$\Omega$-premouse and let $\pi:\N\to\M$ be a weak $n$-embedding. If $\M$ is $n$-maximally iterable 
then so is $\N$.
\end{corollary}

\begin{remark}\label{rem:bad_non-wo_spm_definition}
It seems that one might try 
to define strategy premice over non-wellordered 
sets $A$ by feeding in branches $b_x$ for multiple trees $\Tt_x$ 
simultaneously, thus avoiding the need to select a single tree $\Tt$. However, 
we do not see how to arrange this in such a manner that the branch predicate 
$B$ is always amenable. For example, suppose $A=\RR$, and $\N|\eta$ is given, and we have 
identified, for each $x\in\RR$, a tree $\Tt_x\in\N|\eta$, and now we want to 
feed in $b_x=\Sigma(\Tt_x)$, simultaneously. Let's say we have arranged that 
$\lambda=\lh(\Tt_x)$ is independent of $x$. Then we can easily knit together 
the predicates used to define $\BBB(\N|\eta,\Tt_x,b_x)$, as $x$ ranges over 
$\RR$. Let $\M$ be the resulting structure and let $B=B^\M$. For $B$ to be 
amenable, for each $\alpha<\lambda$, we must have that the function 
$B_\alpha$ is in $\M$, where $B_\alpha(x)= b_x\inter\alpha$. But it seems 
that even $B_2$ could contain non-trivial information, and maybe 
$B_2\notin\M$; note that essentially, $B_2\sub\RR$. Maybe one could first add the sets $B_\alpha$ 
(amenably). But even if one achieved this, it seems that the first problem described in 
\ref{rem:bad_spm_definition} 
would be an obstacle to proving that the resulting hierarchy has nice 
condensation.
\end{remark}
\section{$\Hh^\M$, the local $\hod_{a_0}^\M$}\label{sec:HOD_J_analysis}
\begin{lemma}\label{Sigma1Elem}
Let $\M$ be a $\Theta$-g-organized $\Omega$-pm such that $\M\sats``\Theta$ exists''.
Let $\theta = \Theta^\M$. Let $n_0\leq\om$ be such that $\M$ is $n_0$-sound and $\rho_{n_0}^\M \geq 
\theta$.
Let $\gamma_0=\l(\M)$. Assume that for all $(\xi, k) <_\lex
(\gamma_0,n_0)$, $\M|\xi$ is countably $\GOmega$-$(k,\om_1+1)$-iterable. Assume $\DC_{\RR^\M}$.
Then \tu{(}i\tu{)}
\[ \M|\theta\elem_{\Sigma_1(\Ll^-)}\M\]
and \tu{(}ii\tu{)} for any $(\xi,k)<_\lex(\gamma_0,n_0)$ with $\theta\leq\xi$, and any
$a\in\M|\theta$,
\[ \cHull_{n+1}^{\M|\xi}(\RR^\M\un\{a\})\pins\M|\theta.\]
\end{lemma}
\begin{proof}
(i) from (ii): Let $\varphi\in\Ll^-$ be $\Sigma_1$
and $a\in\M|\theta$.
Suppose
$\M\models\varphi(a)$. We
must show that $\M|\theta\models\varphi(a)$. Let $\xi<\gamma_0$ be least such
that $\M|(\xi+1)\models\varphi(a)$. We need to see that $\xi<\theta$. Assume $\theta\leq\xi$. Fix 
$n<\omega$ and an $\rSigma_{n+1}$
formula $\psi\in\Ll$ such that $\M|\xi\models\psi(a)$, and for any 
hpm
$\N$ and $a'\in\N$, if $\N\models\psi(a')$ then
$\J^\hpm(\N)\models\varphi(a')$. Let
\[ \Hh=\cHull_{n+1}^{\M|\xi}(\RR^\M\un\{a\}).\]
Then $a\in\Hh$ and
$\J_1(\Hh)\models\varphi(a)$. But by (ii),
$\Hh\pins\M|\theta$, a contradiction.

(ii):
For 
$\eta<\theta$, let
$\Hh_\eta=\cHull_{n+1}^{\M|\xi}(\RR^\M\un\eta)$,
and
$\pi_\eta:\Hh_\eta\to\M|\xi$ be the
uncollapse. Note that $\crit(\pi_\eta)$ exists iff $\Hh_\eta\sats$``$\Theta$ exists'',
and $\crit(\pi_\eta)=\Theta^{\Hh_\eta}$ when they exist. Let $\theta_\eta=\Theta^{\Hh_\eta}$ (where 
$\Theta^{\Hh_\eta}=\OR(\Hh_\eta)$ if $\Hh_\eta\sats$``$\Theta$ does not exist''). Then 
$\Hh_\eta\in\M|\theta$ and
$\theta_\eta<\theta$, since $\rho_{n+1}^{\M|\xi}\neq\om$. We say $\eta$
is
a \dfnemph{generator} iff $\eta=\theta_\eta$. The generators are
club in $\theta$. Let $\Hh'_\eta$ be the least $\Hh\pins\M|\theta$
such that $\eta\leq\OR(\Hh)$ and $\rho_\om^\Hh=\om$. Now
$\cHull_{n+1}^{\M|\xi}(\RR^\M\un\{a\})=\Hh_\eta$ for some generator $\eta$. So
the following claim finishes the proof:

\begin{claim} Let $\eta<\theta$ be a generator. Then:
\begin{enumerate}[label=--]
 \item $\Hh_\eta\ins
\Hh'_\eta\pins\M|\theta$.
 \item If $\eta$ is the least generator then
$\rho_{n+1}^{\Hh_\eta}=\om$ and $p_{n+1}^{\Hh_\eta}=\emptyset$.
 \item If $\zeta<\eta$ is the largest generator
$<\eta$, then $\rho_{n+1}^{\Hh_\eta}=\om$ and $p_{n+1}^{\Hh_\eta}=\{\zeta\}$.
\item If $\eta$ is
a limit of generators then $\rho_{n+1}^{\Hh_\eta}=\eta$ and
$p_{n+1}^{\Hh_\eta}=\emptyset$.
\end{enumerate}
\end{claim}

\begin{proof} The proof is by induction on $\eta$.

Suppose $\eta$ is the least generator. Clearly 
$\hmb^{\Hh_\eta}=\hmb^\M$ and $\om_1^\M<\eta$ and
$\Hh_\eta=\cHull_{n+1}^{\Hh_\eta}(\RR^\M)$, which gives that $\rho_{n+1}^{\Hh_\eta}=\om$ and 
$p_{n+1}^{\Hh_\eta}=\emptyset$ and
$\Hh_\eta$ is a fully sound $\Theta$-g-organized $\Omega$-pm. So by $\DC_{\RR^\M}$, countable 
iterability and \ref{cor:Tg_pullback}, we have
$\Hh_\eta\pins\M|\theta$, and $\Hh_\eta=\Hh'_\eta$ since $\eta=\Theta^{\Hh_\eta}$.

Now suppose $\zeta$ is the largest generator $<\eta$. Then
\[ \eta\sub
Y\eqdef\Hull_{n+1}^{\M|\xi}(\RR^\M\un\{\zeta\}),\]
so $\rho_{n+1}^{\Hh_\eta}=\om$
and $p_{n+1}^{\Hh_\eta}\leq\{\zeta\}$. But $\Hh'_\zeta\in Y$, so $\Hh'_\zeta\sub Y$
and $\Hh_\zeta\in Y$. Therefore $p_{n+1}^{\Hh_\eta}=\{\zeta\}$ and $\Hh_\eta$ is
$(n+1)$-solid, and $(n+1)$-sound, so fully sound. The rest is as in the previous
case; again we get $\Hh'_\eta=\Hh_\eta$.

Suppose $\eta$ is a limit of generators. The
$\rSigma_{n+1}$ facts about $\Hh_\eta$ follow readily by induction. Since
$\rho_{n+1}^{\Hh_\eta}=\eta=\Theta^{\Hh_\eta}$ and $\Hh_\eta$ is $(n+1)$-sound, and
$\Hh_\eta$ cannot have extenders overlapping $\eta$, comparison
gives $\Hh_\eta\ins\Hh'_\eta$, as
required.
\end{proof}\let\qed\relax
\end{proof}

\begin{definition}\label{dfn:T^M}
Let $\M$ be a $\Theta$-g-organized $\Omega$-pm satisfying ``$\Theta$ exists'' and 
$\theta=\Theta^\M$.
Let
\[ \tilde{T}^\M\eqdef\Th_{\Sigma_1(\Ll^-)}^{\M|\theta}(\theta\un\{a_0\}). \]
Let $W^\M=\J_\theta[\tilde{T}^\M]$ and
$T^\M=(W^\M,\tilde{T}^\M)$.
We say that a set of ordinals 
$A$ is $\OD_{a_0}^\M$ iff $A\in\M$ and there is 
$\xi<\l(\M)$ such that $A$ is definable from $a_0$ and ordinal parameters over 
$\M|\xi$.\footnote{Note that this provides much more expressive power than
$\OD^{\univ{\M}}_{a_0}$.}
\end{definition}

\begin{remark}
With $\M$ as above, note that $\MFsharp,U^\M,U^{'\M}\in W^\M$ (for $\MFsharp$,
this uses the parameter $a_0$) and 
$\Sigma^{\M|\theta}\in\J(W^\M)$.
Let $\pow_{<\theta}$ denote the bounded subsets of $\theta$.
By \ref{Sigma1Elem}, if the 
hypotheses of \ref{Sigma1Elem} hold, then
\[ \tilde{T}^\M=\Th_{\Sigma_1(\Ll^-)}^\M(\theta\un\{a_0\}) \]
and
$\pow_{<\theta}\cap\OD_{a_0}^\M=
\pow_{<\theta}\cap W^\M$.
\end{remark}

\begin{definition}
A $\Theta$-g-organized $\Omega$-pm is \dfnemph{relevant} iff $\M\sats$``$\Theta$ 
exists'' and $\ex\N\pins\M[\Theta^\M<\OR(\N)$ and $\N\sats\ZF]$.
\end{definition}

\begin{definition}\label{dfn:H^M}
Adopt the hypotheses of \ref{Sigma1Elem}, and suppose $\M$ is relevant.
We define a g-organized $\Omega$-pm $\Hh\eqdef\Hh^\M$ over $\widehat{T^\M}$, with 
$\OR(\Hh)=\OR(\M)$, much as in 
\cite{ScalesK(R)}.
(We show in \ref{local def} that $\Hh$ is indeed a g-organized $\Omega$-pm. It is natural to 
consider $\Hh^\M$ as a locally defined $\H_{a_0}^\M$.)

Let $\theta=\Theta^\M$. Set $\hmb^\Hh=\widehat{T^\M}$, $\hmp^\Hh=\MFsharp$ and 
$\hmPsi^\Hh=\Sigma^{\M|\theta}$.
For $\alpha\geq 1$ define the predicates of $\Hh|\alpha$
by restricting those of $\M|\theta+\alpha$, setting
(i) $\hmP^{\Hh|\alpha}=\hmP^{\M|\theta+\alpha}$ and 
(ii) $\hmE^{\Hh|\alpha}=\hmE^{\M|\theta+\alpha}\cap \Hh|\alpha$.
\end{definition}

Continue with the notation above.
Note that $\PP\in\Hh|2$, where $\PP$ is the Vopenka 
algebra defined over 
$\M|\theta$ as in
\cite{ScalesK(R)}.
Let $\zeta>\theta$ be least such that $\M|\zeta\sats\ZF$.
Note that $\M|\alpha$ is passive for all $\alpha\leq\zeta$, because $\M|\theta$ 
is $(\g,\Theta)$-whole. For $\alpha\geq\zeta$ we have $\theta+\alpha=\alpha$,
and to see that $\Hh$ is indeed a g-organized $\Omega$-pm
we will need to consider how $\M|\alpha=\M|(\theta+\alpha)$ relates 
to $\Hh|\alpha$.
We will observe that for
$\alpha\geq\zeta$, $\Hh|\alpha$ is a g-organized $\Omega$-pm, $\M|\alpha$ is a 
symmetric submodel of a generic extension of $\Hh|\alpha$ (via $\PP$), 
$\Sigma^{\Hh|\alpha}=\Sigma^{\M|\alpha}$, that $\M|\alpha$ is a $(\g,\Theta)$-activation level of 
$\M$ iff $\Hh|\alpha$ is a $\g$-activation level 
of $\Hh$, and $\Tt_{\M|\alpha}\rest\gamma=\Tt_{\Hh|\alpha}\rest\gamma$ for enough
$\gamma$ that condition (i) above will be appropriate.
We will also need to see that the fine structures of $\Hh|\alpha$ and $\M|\alpha$ correspond
appropriately. The fine structural correspondence is mostly as in
\cite{ScalesK(R)}, so we omit most of the details, but give a summary.

\begin{definition}
Adopt the hypotheses of \ref{dfn:H^M} and the notation above. For $\alpha\geq\zeta$ and 
$\Ii=\Hh|\alpha$ we define the $\Ll$-structure
\[ 
\Hh_\alpha(\RR^\M)=\Ii(\RR^\M)=(\J_{\alpha}^{\hmPvec^{\Hh}}(\widehat{T^\M}\un\HC^\M),
\hmPvec^{\Ii},\widehat{T^\M},\hmE^{\Ii},\hmP^{\Ii};\MFsharp,\Sigma^{\M|\theta}
).\qedhere \]
\end{definition}

Truth in
$\Ii(\RR^\M)$ can be reduced to truth in
$\Ii$ via forcing with $\PP$. And
$\Ii(\RR^\M)$ determines $\M|\alpha$: if 
$\M|\theta\in\Hh_\zeta(\RR^\M)$ then
$\es_+^{\Ii}$ determines $\es_+^{\M|\alpha}\rest[\theta,\alpha]$ by the local
definability of the forcing; because $\MFsharp,U,U'\in\Hh|1$ and by induction applied to 
relevant initial segments of $\M|\theta$, we do have
$\M|\theta\in\Hh_\zeta(\RR^\M)$.
The main facts, which generalize \cite[3.9]{ScalesK(R)}, are summarized as follows:

\begin{lemma}
\label{local def} Under the hypotheses of \ref{dfn:H^M} and with $\zeta$ as above, we have:
\begin{enumerate}[label=\tu{(}\arabic*\tu{)}]
\item\label{local def i:0} For relevant $\N\ins\M$, 
$\N||\OR(\N)$ is $\Sigma_1(\Ll^-)$ over 
$\Hh^{\N}(\RR^\M)$, and $\N$ is $\Sigma_1(\Ll)$ over 
$\Hh^{\N}(\RR^\M)$, uniformly in $\N$.
\item $\Hh$ is an $n_0$-sound g-organized $\Omega$-pm \tu{(}over $\widehat{T^\M}$, with 
$\hmPsi^\Hh=\Sigma^{\M|\theta}$\tu{)}, $\theta$ is a cardinal of $\Hh$, and $\zeta$ is least such 
that $\Hh|\zeta\sats\ZF$.
\item\label{local def i:2} For all $(\beta,k)\leq_\lex(\l(\M),n_0)$ with
$\zeta\leq\beta$, we have
$\rho_k(\Hh|\beta) = \rho_k(\M|\beta)$ and $p_k(\Hh|\beta) =
p_k(\M|\beta)\backslash\{\theta\}$.
\item\label{item:fsr_forcing_ext} For all $\beta\in[\zeta,\l(\M)]$, for 
any $p\in\PP$, $\Hh_\beta(\RR^\M)$ is a symmetric inner model of a $\PP$-forcing 
extension of $\Hh|\beta$.
\item\label{item:fsr_forcing_ext_determines_M} For all $\beta\in[\zeta,\l(\M)]$, $\M|\beta$ is 
determined by $\Hh_\beta(\RR^\M)$ as 
described above.
\item\label{local def i:3} Let $\beta\in[\theta,\l(\M)]$. Then $\M|\beta$
is $(\g,\Theta)$-whole iff either $\beta=\theta$, or $\beta>\zeta$ and 
$\Hh|\beta$ is $\g$-whole. Similarly, $\M|\beta$ is a $(\g,\Theta)$-activation level of $\M$ iff 
$\Hh|\beta$ is a $\g$-activation level of $\Hh$.
\end{enumerate}
\end{lemma}
\begin{proof}[Proof sketch]
For most of the details, see the proof of \cite[3.9]{ScalesK(R)}. We just
give enough of a sketch to describe the new features.

As usual, \ref{local def i:0} will follow from the proof, and by 
induction, we may assume that \ref{local def i:0} holds for $\N\ins\M|\theta$.
This implies $\M|\theta\in\Hh_\zeta(\RR^\M)$, unless there is no relevant $\xi<\theta$ (a fact 
regarding which $T^\M$ informs us). In the latter case, 
$\M|\theta=\J^\hpm_\theta(\hmb^\M;\MFsharp,\emptyset)$. But 
$U^\M\in W^\M$, so 
$\Upsilon^\M,\hmb^\M\in\Hh_\zeta(\RR^\M)$, which suffices.

Let $\eta\in[\zeta,\l(\M)]$. We say that $\M|\eta,\Hh|\eta$ are \dfnemph{fine structurally related} 
iff
\ref{local def i:2}, \ref{item:fsr_forcing_ext} and 
\ref{item:fsr_forcing_ext_determines_M} hold 
for $\beta\leq\eta$.
We say that $\M|\eta$, $\Hh|\eta$ are \dfnemph{$\g$-related} iff \ref{local def i:3} holds for 
$\beta\leq\eta$. We say that $\M|\eta,\Hh|\eta$ are \dfnemph{related} iff they are both fine 
structurally related and $\g$-related.

\begin{claim}
For
$\eta\in[\zeta,\l(\M)]$, $\Hh|\eta$ is a g-organized 
$\Omega$-pm over $\widehat{T^\M}$, and the models
$\M|\eta,\Hh|\eta$ are related, and uniformly so in $\eta$.
\end{claim}
\begin{proof}
By induction on $\eta$; the uniformity follows from the proof.
Let $\M_0\pins\M$ be $(\g,\Theta)$-whole, with $\theta\leq\beta_0\eqdef\l(\M_0)$,
and suppose that if $\theta<\beta_0$ then the claim holds for all $\eta\in[\zeta,\beta_0]$. (By 
$(\g,\Theta)$-wholeness, either $\beta_0=\theta$ or $\zeta<\beta_0$.)
For simplicity, suppose that ($*$) there is a $(\g,\Theta)$-whole
$\M_1\ins\M$ such that $\M_0\pins\M_1$. Let $\M_1$ be least such and $\beta_1=\l(\M_1)$.
We will prove the claim for $\eta\in(\beta_0,\beta_1]$.
The fact that $\M_1$ and $\Hh^{\M_1}$ are fine structurally related is proved as in 
\cite[3.9]{ScalesK(R)} (this is actually easier than in \cite{ScalesK(R)}, as we have 
$\hmP^{\Hh|\eta}=\hmP^{\M|\eta}$ and $\hmE^{\Hh|\eta}=\emptyset=\hmE^{\M|\eta}$ for all 
$\eta\in(\beta_0,\beta_1]$).
It remains to see that they are $\g$-related. For this we need to see that
\begin{enumerate}[label=--]
 \item $\alpha$ is least such that $\alpha>\beta_0$ and $\Hh|\alpha\sats\ZF$, and
 \item 
$\Tt\eqdef\Tt_{\Hh|\alpha}=\Uu\eqdef\Uu_{\M|\alpha}$.
\end{enumerate}
The former is straightforward, using forcing as in \cite[3.9]{ScalesK(R)}.
So $\Hh|\alpha$ is the next $\g$-activation level of $\Hh$, beyond $\Hh|\beta_0$ if 
$\beta_0>\theta$, or at all if $\beta_0=\theta$.
We now prove 
by induction on $\gamma$ that $\Tt\rest\gamma+1=\Uu\rest\gamma+1$ for all 
$\gamma+1\leq\eps=\max(\lh(\Tt),\lh(\Uu))$.
But then $\Tt=\Uu$ as required.

We have $\Tt\rest\alpha+1=\Uu\rest\alpha+1$ (this part is linear iteration). So let 
$\gamma\geq\alpha$ and suppose that $\Tt\rest\gamma+1=\Uu\rest\gamma+1$ and 
$\gamma+1<\eps$; we just need to see that $E^\Tt_\gamma=E^\Uu_\gamma$ (and in particular, both are 
defined). Let $\xi$ be the largest limit ordinal such that $\xi\leq\gamma$. Let
$\Ss=\Ss^\M_{\Uu\rest\xi+1}$. Let $\delta=\l(\Ss)+\OR(M^\Uu_\gamma)$.
So $\delta\leq\l(\M_1)$.

Suppose that $E^\Tt_\gamma\neq\emptyset$.
Let 
$p\in\Coll(\om,\Hh|\alpha)$ be such that $p$ forces, 
over\footnote{This forcing is absolute, but the point is that the 
relevant forcing relation is in $\Hh|\delta$.} 
$\Hh|\delta$, that $E^\Tt_\gamma$ induces an 
axiom which fails for $\dot{x}_{\Hh|\alpha}$. Now in $\M|\delta$,
$\QQ\eqdef\Coll(\om,\M|\alpha)$ factors naturally as $\QQ_0\cross\QQ$ where 
$\QQ_0=\Coll(\om,\Hh|\alpha)$. Let $\dot{G}_0,\dot{G}_1$ be the resulting $\QQ$-names for the 
factor generics (so under the factoring just mentioned, $\dot{G}_0\cross\dot{G}_1$ corresponds to 
$\dot{G}$, the standard $\QQ$-name for the $\QQ$-generic).
Let 
$\dot{x}_{0,\M|\alpha}$ and 
$\dot{x}_{1,\M|\alpha}$ be the $\QQ$-names for the generic reals determined by $\dot{G}_0$ and 
$\dot{G}_1$. Let $p'\in\Coll(\om,\M|\alpha)$ force
that $p\in\dot{G}_0$. we have that $p'$ forces that 
$E^\Tt_\gamma$ induces an axiom which fails for $\dot{x}_{0,\M|\alpha}$. But 
assuming we have used the natural definitions, $\dot{x}_{0,\M|\alpha}$ is arithmetic in 
$\dot{x}_{\M|\alpha}$, and so it is easy to see that $p'$ forces that 
$E^\Tt_\gamma$ induces an axiom which fails for $\dot{x}_{\M|\alpha}$, as 
required.

The case that $E^\Uu_\gamma\neq\emptyset$ is similar, but we need to use 
the fact that $\M|\delta$ can be realized as a symmetric 
submodel of a $\PP$-generic extension of $\Hh|\delta$. (It 
doesn't suffice that this holds for $\M|\alpha$ and $\Hh|\alpha$, since the 
forcing relation which demonstrates the fact that $E^\Uu_\gamma$ induces a bad 
axiom need not be in $\M|\alpha$.) We omit further detail.

If ($*$) fails then it is almost the same. However, suppose there is a 
$(\g,\Theta)$-activation level $\ins\M$ beyond $\M_0$.
Then it can be that $\Tt\neq\Uu$, where $\Tt,\Uu$ are as before.
However, the preceding argument still shows that $\Tt\rest\gamma+1=\Uu\rest\gamma+1$ for 
enough ordinals $\gamma$ that the proof goes through.

The remaining details (in particular the fact that $\hmE^{\Hh_\alpha(\RR^\M)}$
determines $\hmE^{\M|\alpha}$) are as in \cite{ScalesK(R)}. This 
completes the sketch of the proof of the lemma.
\end{proof}
\renewcommand{\qed}{}
\end{proof}
\indent The next theorem relates the iterability of $\Hh$ and $\M$. The proof of 
\ref{IterabilityLifts}
uses \ref{local def} and is just like that in
\cite[3.18]{ScalesK(R)}.
\begin{theorem}
\label{IterabilityLifts}
Assume the hypotheses of \ref{dfn:H^M}. Let $\gamma\in\Ord$. Then
$\Hh^\M$ is \tu{(}countably\tu{)} $(n_0,\gamma)$-iterable iff 
$\M$ is \tu{(}countably\tu{)} above-$\Theta^\M$ $(n_0,\gamma)$-iterable.
\end{theorem}

\begin{remark}\label{rem:S-construction}
Constructions having the flavor of \ref{dfn:H^M}, as well as their inverses, are referred to as 
\dfnemph{S-constructions}. In the sequel, we will also need S-construction, performed mostly as in 
\ref{dfn:H^M}, for example, in the 
following 
context. Let $\M$ be a g-organized $\Omega$-pm. Let $\N\pins\M$ be a $\g$-whole 
strong cutpoint of $\M$. Let 
$g\sub\Coll(\om,\N)$ be $\M$-generic. Then $\M[g]$ can be reorganized as a g-organized 
$\Omega$-pm $\M[g]^*$ over $\hat{x}$ where $x=(\N,g)$, with $\hmPsi^{\M[g]^*}=\Sigma^\N$. Moreover, 
the fine structure and iterability of 
$\M[g]^*$ corresponds to the fine structure and iterability of $\M$ above $\eta$, in a manner 
similar to \ref{local def} and \ref{IterabilityLifts}. We leave the precise formulation and proofs 
of these facts to the reader.

Assume $\DC_\RR$ and suppose that $\kappa_0\geq\Theta$ (see \ref{ntn:fixed_notation}). Using 
similar arguments, we also get that
$\M\eqdef\Lp^{\gOmega}(\RR)$ and $\N\eqdef\Lp^{\GOmega}(\RR)$
and $\P\eqdef\Lp^{\GOmega}(\HC,\Omega\rest\HC)$ have the same $\pow(\RR)$ (we have
$\Omega\rest\HC\in\M\inter\N$ by \ref{ClosedUnderJ}
and \ref{lem:charac_Tg-org}). Moreover, if $\Q=\Lp^\Omega(\RR)$ is well-defined and 
$\Omega$ 
has a property along the lines of \emph{relativizes well} (see \cite[Definition 1.3.21(?)]{cmi}) 
then 
the same holds of $\Q$. In fact, $\M$, $\N$, $\P$ (and $\Q$) have literally the same extender 
sequences and for 
all $\alpha$ such that $\M|\alpha$ is $\hmE$-active, there is a straightforward translation between 
$\M|\alpha$, $\N|\alpha$, $\P|\alpha$ (and $\Q|\alpha$). (To see that $\Q|\alpha$
computes the others, note that the $\hmP$-predicates 
of the others are determined by Q-structures for trees $\Tt$, where the 
Q-structures are in 
$L^\Omega(M(\Tt))$.)
\end{remark}

\section{Scales}\label{sec:scales}
We now begin the main project of the paper: the analysis of scales in $\Theta$-g-organized 
$\Omega$-premice.\footnote{Let $\M$ be an 
hpm.
When we say that
$\M\sats$``$\rSigma_n$ has the scale property'',
recall that $\rSigma_n$ uses the language $\Ll^+$
and $\rSigma_n$ formulas are interpreted over $\core_0(\M)$,
so the statement literally means that 
$\core_0(\M)\sats$``$\rSigma_n$ has the scale property''.
Moreover, since it is a statement satisfied by $\core_0(\M)$,
it is interpreted with respect to sequences of reals 
$\left<x_n\right>_{n<\om}\in\M$.
Likewise when we say that $\M\sats$``$\rSigma_n(\RR)$ has the scale property'', and
here any parameters in $\RR^\M$ are allowed; for $\bfrSigma_n$, any parameters in $\core_0(\M)$ are 
allowed, as usual.}
In our application to the core model induction, the analysis proceeds from 
optimal determinacy hypotheses; cf.
\cite{Scalesweakgap}.\footnote{Let 
$\Sigma$ be the unique iteration strategy for $\M_1^\sharp$. Suppose 
$\Lp^{^\gTheta\Sigma}(\RR)\vDash \AD^+ + \MC$. Then in fact 
$\Lp^{^\gTheta\Sigma}(\RR)\cap \pow(\RR) = 
\Lp(\RR)\cap \pow(\RR)$. This is because in 
$L(\Lp^{^\gTheta\Sigma}(\RR))$, $L(\pow(\RR)) \vDash \AD^+ \ 
+ \ \Theta = \theta_0 + \MC$ and hence by \cite{sargsyan2014Rmice}, in 
$L(\Lp^{^\gTheta\Sigma}(\RR))$, $\pow(\RR)\subseteq 
\Lp(\RR)$.  Therefore, even though the hierarchies 
$\Lp(\RR)$ and $\Lp^{^\gTheta\Sigma}(\RR)$ are different, as far 
as sets of reals are concerned, we don't lose any information by analyzing the 
scales pattern in $\Lp^{^\gTheta\Sigma}(\RR)$ instead of that in 
$\Lp(\RR)$.}

\subsection{Scales on $\Sigma_1^\M$ sets for passive $\M$}

\begin{theorem}
\label{passiveGameScale}
Let $\M$ be a passive $\Theta$-g-organized 
$\Omega$-pm satisfying $\AD$. Assume $\DC_{\RR^\M}$. Suppose that every proper segment of $\M$ 
is countably $\GOmega$-$(\om,\om_1+1)$-iterable. Then $\M 
\sats$``$\rSigma_1(a_0)$ has the scale property''.
\end{theorem}

\begin{proof}
By $\DC_{\RR^\M}$, \ref{lem:Th-g_very_con} and \ref{cor:Tg_pullback} we may assume that $\M$ is 
countable.
For simplicity we assume that $\l(\M)$ is a limit ordinal; for the contrary 
case make the usual modifications using the $\Ss$-hierarchy as in 
\ref{thm:scale_branch_active_immediately_after_new_Sigma_1} below.
For this proof we abbreviate $\RR^\M$ with $\RR$, and likewise interpret $\HC$ and terms 
like \emph{real, analytical}, etc, over $\M$.

Let $\Phi\in\Ll^-$ be $\Sigma_1$. For $x \in 
\RR$, let
$A(x) \Leftrightarrow \M \vDash \Phi(x)$. We will define a 
$\Sigma^\M_1(a_0)$-scale on $A$.
For $x \in \RR$ and $1\leq\beta<\l(\M)$ let
$A^\beta(x) \Leftrightarrow \M|\beta \vDash \Phi(x)$.
Then $A=\bigcup_{\beta<\l(\M)}A^\beta$.
We will construct a closed game representation $x \mapsto
G^{\beta}_x$ for $A^\beta$, with $G^\beta_x$ continuously associated to 
$x$. For $u$ a partial play of $G^\beta_x$, let $G^\beta_{x,u}$ be the game in which the players 
continue the play of $G^\beta_x$ from $u$. Let
$A^\beta_k$ be the set of pairs $(x,u)$ such that $\M\sats$``$u$ is a partial play of 
$G^\beta_x$ of length $k$, and player $\plI$ has a winning quasi-strategy in 
$G^\beta_{x,u}$''.\footnote{Note 
that we require the winning quasi-strategy to be in $\M$, unlike in Steel's arguments.}
We will eventually show that $A^\beta_k \in \M$ and the map
$(\beta,k) \mapsto A^\beta_k$ is $\Sigma_1^\M(a_0)$.\footnote{Because we required that the winning 
quasi-strategy be in $\M$, we already know that $A^\beta_k$ is definable over $\M$, but this does 
not particularly help us prove that $A^\beta_k\in\M$.}

The foregoing yields a $\Sigma_1^\M(a_0)$ scale 
essentially as in \cite{Scales}. However, let us mention two small differences.

Firstly, there will be certain moves in $G^\beta_x$ which are, in the sense of $\M$, 
prewellordering equivalence classes of reals. In the scale computation, these are simply treated in 
the same manner that ordinals are treated in \cite{Scales}.
Let us set up some notation for these moves. Let $\Upsilon=\Upsilon^\M$ and
$\vec{\leq}=\left<\leq_n\right>_{n<\om}\eqdef\vec{\leq}^\M$ (see \ref{dfn:Theta-g-organized}).
For $n<\om$ let $e_n$ be the set of $\leq_n$-equivalence classes of reals. Let 
$e=\bigcup_{n<\om}e_n$. Let $W$ be the tree of the scale in the codes in $e$; so
$W$ is a tree on $\om\cross e$ and $p[W]=\Upsilon^\HCc$ (in $\M$).
Let $W'$ be likewise for $\RR\cut\Upsilon$. (If $\M$ has an admissible initial segment then we 
could just use $U^\M,U'^\M$ instead of $W,W'$.)
Then $\{(\Upsilon,\Upsilon^\HCc,W,W')\}$ is $\Delta_1^\M$.

Now secondly, because the payoff 
is closed for player $\plI$, if $\M\sats$``$\Sigma$ is a winning quasi-strategy for player 
$\plI$ for $G^\beta_{x,u}$'',
then $V$ satisfies the same. However, $\M$ might 
not have a winning quasi-strategy for player $\plI$ for $G^\beta_{x,u}$, although $V$ does. But 
this does not cause problems 
for the computation of a $\Sigma_1^\M(a_0)$ scale. For the fact that each $A^\beta_k\in\M$
ensures that $A^\beta_{k}$ is related to $A^\beta_{k+1}$ in essentially the usual manner.
That is, $A^\beta_k$ is either of the form $\ex^\RR A^\beta_{k+1}$, or $\ex\alpha<\beta[ 
A^\beta_{k+1}]$, or $\ex n\leq k\ex X\in e_n[A^\beta_{k+1}]$, or $\all^\RR A^\beta_{k+1}$.
Because the relevant computations propagating norms are made inside $\M$ -- where, in 
particular, $\leq_n$ and $\leq'_n$ are wellfounded -- this is enough for the scale computation.

Before defining $G^\beta_x$ we give an outline. Player 
$\plII$ will play reals. Player $\plI$ will (attempt to) build a countable, iterable, 
passive, $\Theta$-g-organized $\Omega$-pm $\P$ over $\Upsilon\cap\P$, 
containing all 
reals played by player $\plII$, such that $\P\sats\Phi(x)$, but for all 
$\gamma<\l(\P)$, $\P|\gamma\sats\neg\Phi(x)$. To ensure that player $\plI$ indeed 
plays an iterable $\Theta$-g-organized $\Omega$-pm, he must 
simultaneously 
build a (cofinal) very weak $0$-embedding $\pi:\P\to\M|\gamma$ for some $\gamma\leq\beta$
\footnote{One could have instead used an approach more like that used in 
\cite{ScalesK(R)}.}. To ensure that $\P$ is over $\Upsilon\inter\P$, he must also build 
various branches through $W$ and $W'$. (Here we will be interested in the case that
those branches appear in generic extensions of $\M$, which will ensure that really prove that a 
given real of $\M$ is in the set it is claimed to be in.)

We now proceed to the details. Player $\plI$ will describe his model using the language
\[ \Ll^* 
=_{\textrm{def}}
\Ll \cup \{\dot{x_i} \ | \ i
<\omega\}\cup\{\dot{\Upsilon}\}.
\]
Here $\dot{x_i}$ and $\dot{\Upsilon}$ are constants; 
$\dot{x_i}$ will denote the $i^\nth$ real played in the game. Fix
recursive maps
\begin{equation*}
m,n: \{\sigma\ |\ \sigma\textrm{ is an }\Ll^*\textrm{-formula}\} \rightarrow 
\{2n \ | \ 2 \leq n < \omega\}
\end{equation*}
which are one-to-one, have disjoint recursive ranges, and are such that whenever
$\dot{x_i}$ occurs in $\sigma$, then $i <  \min(m(\sigma),n(\sigma))$.

Fix a $\Sigma_1(\Ll^-)$ formula $\sigma_0(v_0,v_1,v_2)$ that defines over each $\M|\gamma$,
a map
\[h_\gamma:\gamma^{<\omega}\times\RR
\onto \M|\gamma.\]
Let $T$ be the
following $\mathcal{L}^*$ theory:
\begin{IEEEeqnarray*}{ll}
(1) &\text{Extensionality}\\
(2) &\text{``}V\text{ is a }\Theta\text{-}\g\text{-spm''}\\
(3)_i &\ \dot{x_i} \in \RR\\
(4)&\ \Phi(\dot{x}_2)\
\&\ \neg\ex\N\pins V\left[\N\sats\Phi(\dot{x}_2)\right]\\
(5)&\ \forall u,v,y,z\left[\sigma_0(u,v,y)\wedge
\sigma_0(u,v,z) \implies y = z\right]\\
(6)_{\varphi} &\ [\exists v \varphi(v)] \implies \exists v\exists
F\in \l(V)^{<\omega}\left[\varphi(v) \wedge
\sigma_0(F,\dot{x}_{m(\varphi)},v)\right]\\
(7)_\varphi &\ \exists v\left[\varphi(v) \wedge v\in\RR\right]
\implies\varphi(\dot{x}_{n(\varphi)})\\
(8)&\ \hmbdot=\widehat{x}\text{ where }x=(\HC,\dot{\Upsilon})\\
(9)&\ \hmpdot\text{ is an hpm over the transitive set coded by }\dot{x}_0
\end{IEEEeqnarray*}

A run of the game $G^\beta_x$ has $\omega$ rounds. In round $n$,
player $\plI$ first plays $i_n,x_{2n},\eta_n,\Lambda_{n}$ where $i_n \in \{0,1\}$, $x_{2n} \in
\RR$, $\eta_0\leq\beta$ and $\eta_{n+1}<\OR(\M|\eta_0)$, and $\Lambda_{n}\in(W\cup W')^n$; player 
$\plII$ plays then $x_{2n+1} \in
\RR$.

The payoff for player $\plI$ is mostly analogous to that in
\cite{ScalesK(R)}. Conditions \ref{item:branches_cohere_inf} and \ref{item:branches_correct} are 
new, and they
ensure that for each $i<\om$, if player $\plI$ asserts, for example,
that ``$\dot{x}_i\in\dot{\Upsilon}^\HCc$'' then $\left<\Lambda_{n,i}\right>_{n\in(i,\om)}$ is an 
infinite 
branch through $W$ witnessing that $x_i\in \Upsilon^\HCc$.

If $u = \langle(i_k,x_{2k},\eta_k,x_{2k+1}) \ | \ k < n\rangle$ is
a partial play of $G^\beta_x$, let
\begin{equation*}
T^*(u)=\{(\neg)^i\sigma\ |\ 
\sigma\textrm{ is an }\Ll^*\textrm{-sentence}\wedge n(\sigma)<n\wedge 
i=i_{n(\sigma)} \},
\end{equation*}
where $(\neg)^0\sigma=\sigma$ and $(\neg)^1\sigma=\neg\sigma$.
If $p$ is a full run of $G^\beta_x$, let
$T^*(p)$ be the union of all $T^*(p\rest n)$, for $n<\om$.
We write ``$\unique v\varphi(v)$'' for ``the unique $v$ such
that $\varphi(v)$''. For $\sigma=((a_0,b_0),\ldots,(a_{n-1},b_{n-1}))$ let 
$p_0[\sigma]=(a_0,\ldots,a_{n-1})$ and $p_1[\sigma]=(b_0,\ldots,b_{n-1})$.

A run $p = \langle(i_k,x_{2k},\eta_k,\Lambda_k,x_{2k+1}) \ | \ k <
\omega\rangle$ of $G^\beta_x$ is a win for player $\plI$ iff
\begin{enumerate}[label=(\alph*)]
\item\label{item:consistent_extension} $T^*(p)$ is a consistent extension of $T$,
\item\label{item:correct_reals} $x_0=a_0$ and $x_2 = x$,
\item for all $i,m,n<\omega$, $\text{``}\dot{x}_i(n) = m\text{''}\in T^*(p)$
iff $x_i(n) = m$,
\item if $\varphi$ and $\psi$ are $\mathcal{L}^*$-formulae with one free 
variable
and
\[ \text{``}\unique v\varphi(v) \in \Ord\ \&\ \unique v\psi(v) \in 
\Ord\text{''}\in
T^*(p),\]
then ``$\unique v\varphi(v) \leq \unique v\psi(v)$''$\in T^*(p)$ iff
$\eta_{n(\varphi)} \leq \eta_{n(\psi)}$,
\item\label{item:embedding} if $\sigma_0,\ldots,\sigma_{n-1}$ are 
$\Ll^*$-formulas with one free variable and
\[ \text{``}\unique v\sigma_k(v)\in\Ord\text{''}\in 
T^*(p) 
\]
for all $k<n$, then for any $\rSigma_1$-formula $\theta(v_0,\dots,v_{n-1},v)$,
\[ 
\theta(\unique v\sigma_0(v),\dots,\unique v\sigma_{n-1}(v),\dot{x}_0)\in
T^*(p) \]
if and only if
\[
\M|\eta_0\sats\theta(\eta_{n(\sigma_0)},\dots,\eta_{n(\sigma_{n-1})},a_0),
\]
\item\label{item:branches_cohere_inf} for all $i<m\leq n<\om$, 
$\Lambda_{m,i}\ins\Lambda_{n,i}$ and $p_0[\Lambda_{n,i}]=x_i\rest n$,
\item\label{item:branches_correct} for all $i<m<\om$, if 
$\text{``}\dot{x}_i\in\dot{\Upsilon}^\HCc\text{''}\in T^*(p)$ then $\Lambda_{m,i}\in W$, and 
otherwise 
$\Lambda_{m,i}\in W'$.
\end{enumerate}

Because of the payoff conditions, we could have added a sentence like ``$\hmpdot$ is 
$\MFsharp$-like'' to 
$T$ (or any other sentences satisfied by all initial segments of 
$\M$), without any significant effect.

We next define the notion of \textit{honesty} and show that the only winning
strategy for player $\plI$ is to be honest.
A partial play
\[ u =
\langle(i_k,x_{2k},\eta_k,\Lambda_k,x_{2k+1}) \ | \ k < n\rangle \]
is
$(\beta,x)$-\dfnemph{honest} iff $\M|\beta\sats\Phi(x)$ and if $n>0$ then letting $\eta$ be least 
such that $\M|\eta\sats\Phi(x)$, we have:
\begin{enumerate}[label=(\roman*)]
\item\label{item:honesty_correct_reals} $x_0=a_0$ and if $n>1$ then $x_2=x$.
\item\label{item:honesty_correct_theory} Let $I_u$ be any interpretation of $\Ll^*$ in which 
$\dot{x}_i^{I_u}=x_i$ for 
$0 < i < 2n$ and $\dot{\Upsilon}^{I_u}=\Upsilon$. Then $(\M|\eta,I_u)\sats T^*(u)$.
\item\label{item:honesty_embedding_existence} Let $\left<\sigma_i\right>_{i<m}$ enumerate all 
formulas 
$\sigma\in\Ll^*$
of one free variable such that $n(\sigma) < n$ and
$(\M|\eta,I_u) \sats\text{``}\unique v\sigma(v) \in\Ord$\text{''}.
For $k<m$, let $\delta_k\in\M|\eta$ be such that
$(\M|\eta,I_u)\sats\sigma_k(\delta_k)$. Then in $\M$,
$\Coll(\om,\RR)$ forces the existence of a partial embedding
\[ \pi:\M|\eta\partialto\M|\eta_0\]
such that $\OR(\M|\eta)\un\{a_0\}\sub\dom(\pi)$,
$\pi(\delta_k)=\eta_{n(\sigma_k)}$ for each $k<m$, 
and $\pi$ is $\rSigma_1$-elementary on its domain.
\item\label{item:honesty_Upsilon} For each $i<m<n$, $\Lambda_{m,i}\ins\Lambda_{n-1,i}$ and 
$x_i\rest 
m=p_0[\Lambda_{m,i}]$, and if $x_i\in \Upsilon^\HCc$ (if $x_i\notin \Upsilon^\HCc$) then 
there is $f\in\M\inter[W_{x_i}]$ ($f\in\M\inter[W'_{x_i}]$) such that $f\rest 
m=p_1[\Lambda_{m,i}]$''.
\end{enumerate}

Let $Q^\beta_k(x,u)$ iff $u$ is a
$(\beta,x)$-honest position of length $k$.

The following two claims complete our proof of Theorem \ref{passiveGameScale}. 
Their proofs are similar to those of \cite[Claims 
4.2, 4.3]{ScalesK(R)}.

\begin{claim}\label{clm:honesty_computable}  
$Q^\beta_k \in \M$ for all $\beta,k$, and the map
$(\beta,k)\mapsto Q^\beta_k$ is $\Sigma_1^\M(a_0)$.\end{claim}
\begin{proof}[Proof Sketch]
For condition \ref{item:honesty_Upsilon}, observe that there is $k<\om$ such that every infinite 
branch $b\in\M$ through $W$ or $W'$,
is in fact in $\Ss_k((\HC,\Upsilon,\MFsharp))$. For let $b=(x,f)\in\M$ be a branch through, say, 
$W$.
Because $W$ is the tree of $\vec{\leq}$, by $\AC_{\om,\RR}$ in $\M$ (where $\AD$ holds),
there is $\left<x_n\right>_{n<\om}\in\M$ such that for each $n$, $x_m\rest n=x\rest n$ and 
$x_m\leq_i x_n\leq_i x_m$ for each $i<n\leq m$. But $\left<x_n\right>_{n<\om}$ determines $b$, and 
gives the observation.

Regarding the other conditions, the proof is mostly like that of \cite[Claim 4.2]{ScalesK(R)}, 
but
we modify some details and give a complete proof of some points only hinted at in 
\cite{ScalesK(R)}. Let $\gamma=\OR(\M|\beta)$,
$A=\Th_1^{\M|\beta}(\gamma\un\{a_0\})$ and $A'=\gamma\un\{A\}$. Let $\lambda\in\Ord$ be 
least such that 
$\J_\lambda(A')$ is admissible. The ``embedding game'' 
$\mathcal{G}$ (see 
\cite[Claim 4.2]{ScalesK(R)}) is definable from $A$ and is fully analysed in 
$\J_\alpha(A')$ for some $\alpha<\lambda$. Now we claim that for 
each 
$\alpha<\lambda$, \[ 
t_\alpha=\Th_1^{\J_\alpha(A')}(A')\in\M. \]
This suffices. For if $N$ is any structure 
with $A'\sub N$ and satisfying ``$V=L[A']$, I 
see a full analysis of 
$\mathcal{G}$ but no proper segment of me does'', then $N$ is wellfounded and 
so $N=\J_\alpha(A')$ for some $\alpha$ (since otherwise the 
wellfounded part of $N$ is 
admissible, contradicting the minimality of $N$). Therefore $\M$ can identify 
the theory of the unique such $N$, allowing the rest of the proof of 
\cite[Claim 4.2]{ScalesK(R)} to go through.

So we show that $t_\alpha\in\M$.
Let $\leq$ be a prewellorder of $\RR^\M$ of length $\geq\gamma$, with 
$\leq$ in $\M$.
Say that a structure $N$ (possibly illfounded) 
is \dfnemph{good} iff $N$ extends $A'$ and 
$N\sats$``$V=L[A']$'' and $N=\Hull_1^N(A')$ and 
$\Th_1^N(A')$ is $(\bfSigma^1_1(\leq))^\M$ (in the codes given by 
$\leq$). We claim that for every $\alpha<\lambda$, $\J_\alpha(A')$ is good (and 
therefore $t_\alpha\in\M$). All requirements are clear other than the fact that 
$t_\alpha$ is $(\bfSigma^1_1(\leq))^\M$.

Now if there is any illfounded good $N$, then the wellfounded part of $N$ is 
admissible, and therefore $\J_\alpha(A')\pins N$ for each $\alpha<\lambda$, 
which easily gives the claim. So suppose all good structures are wellfounded.

We claim that there is a largest good structure. For suppose not. Let $S$ be 
the set of all $\Sigma_1$ theories of good structures. Clearly $S\in\M$. Now 
for each $N\in S$ let $t_N=\Th_1^N(A')$. Let $t=\bigcup S$. Then 
$t\in\M$, and $t=\Th_1^{N}(A')$ for $N=\J_\xi(A')$, for some 
ordinal $\xi$. Moreover, $N=\Hull_1^N(A')$. But then by the 
coding 
lemma applied in $\M$, $N$ is good, contradiction.

So let $N$ be the largest good structure. Let $N=\J_\xi(A')$ and 
$N'=\J_{\xi+1}(A')$. We claim that $N\elem_1 N'$, and 
therefore that $N$ is admissible, completing the proof. So suppose 
otherwise. We claim that $N'$ is good, for a contradiction. Clearly
$N'=\Hull_1^{N'}(A')$, so we just need to see that 
$t'=\Th_1^{N'}(A')$ is $(\bfSigma^1_1(\leq))^\M$. By the 
coding lemma, it suffices to see that $t'\in\M$. Now $t'$ is recursively 
equivalent to $\oplus_{n<\om}T_n$ where $T_n=\Th_n^N(A')$. But 
each of these theories are in $\M$ since $T_1=t_N\in\M$. Therefore, by the 
coding lemma, each $T_n$ is $(\bfSigma_1^1(\leq))^\M$. Let $T$ be the set of 
parameters $x\in\RR$ coding (relative to $(\bfSigma^1_1(\leq))^\M$) one of the 
theories $T_n$, for some $n<\om$. Then $T\in\M$ because in fact, $T$ is 
$(\bfSigma^1_{10}(\leq))^\M$. Therefore $\oplus_{n<\om}T_n\in\M$, as required.
\end{proof}

Because $\G$ is fully analysed inside $\M$, the existence of the embedding in condition 
\ref{item:honesty_embedding_existence} of $(\beta,x)$-honesty is actually absolute between 
$\M^{\Coll(\om,\RR)}$ and $V^{\Coll(\om,\RR)}$.

\begin{claim} $A^\beta_k=Q^\beta_k$.
\end{claim}
\begin{proof}[Proof Sketch]
Let $u\in Q^\beta_k$. Then as in \cite{ScalesK(R)} there is 
$\Sigma\in\M^{\Coll(\om,\RR)}$ which is a winning quasi-strategy for player $\plI$ 
in $G^\beta_{x,u}$. For every $\alpha<\l(\M)$ and 
$n<\om$, the $\Sigma_0$ forcing relation for $\Ss_n(\M|\alpha)$ is in $\M$.
(Note here that for $x\in\HC$, the $\Sigma_0$ forcing relation restricted to elements of 
$\trancl(x)$ is essentially in $\HC$, as it is trivial on conditions $p\notin\trancl(x)$.)
Let $\tilde{\Sigma}\in\M$ and $p\in\Coll(\om,\RR)$ be such that 
in $\M$, $p\forces$``$\tilde{\Sigma}$ is a winning quasi-strategy''.
Let $\Sigma'$ be the set of all partial plays $v$ extending $u$ such that for some $q\leq p$,
in $\M$, $q\forces$``$v$ is according to $\tilde{\Sigma}$''. Then $\Sigma'\in\M$, and it is easy to 
see that $\Sigma'$ is a winning quasi-strategy, so $u\in A^\beta_k$ as required.

Now consider the converse. Let $(u,x)\in A^\beta_k$ and let $\Sigma\in\M$ be a winning 
quasi-strategy witnessing this. Let $G$ be 
$(\M,\Coll(\om,\RR^\M))$-generic (recall we have reduced to the case that $\M$ is countable).
As in the proof of \cite[Claim 4.3]{ScalesK(R)}, but working in $\M[G]$ (where we have $\Sigma$),
let let $\N$ be a model produced by playing 
$G^\beta_{x,u}$ according to $\Sigma$ and having player $\plII$ play out all reals in $\RR^\M$.   
Let $\pi':\N\partialto\M|\eta_0$ be the partial embedding, 
with $\dom(\pi)=\OR(\N)\un\{a_0\}$, provided by 
payoff condition \ref{item:embedding}; so $\pi'$ is $\Sigma_1$-elementary on its 
domain. Now $\hmp^\N$ is 
an hpm over $A_0$ (as $a_0=\dot{x_0}^\N$). Using $\pi'$ and since $\hmpdot\in\Ll$, it easily 
follows that $\hmp^\N=\MFsharp$, and that $\pi$ extends uniquely to very weak $0$-embedding 
$\pi:\N\partialto\M|\eta_0$ which is $\Sigma_1$-elementary on its domain. It follows that
$\N$ is a $\Theta$-$\g$-spm with $\RR^\N=\RR^\M$, and in fact, $\N$ 
is a $\Theta$-$\g$-organized $\Omega$-pm over some $\Upsilon'$, by \ref{lem:Th-g_very_con}. 

Actually,
$\Upsilon'=\Upsilon$. This is because
because player $\plI$ built witnessing branches 
through $W,W'$, and because if $x\in\RR^\M$ and $\M[G]\sats$``$x\in p[W]$'', for example,
then $\M\sats$``$x\in p[W]$''. The latter is because the relevant forcing relations are in 
$\M$, and so, if $p\forces$``$b\in[W_x]$'' then $\M$ can compute the left-most branch 
$b'\in[W_x]$ such that for all $n<\om$, there is some $q\leq p$ forcing ``$b\rest n=b'\rest 
n$''. Similar considerations also give condition \ref{item:honesty_Upsilon} of 
$(\beta,x)$-honesty (the relevant branches are in $\M$, not just $\M[G]$).

Now $\N\sats\Phi(x)$ but no $\N'\pins\N$ satisfies $\Phi(x)$, so 
$\l(\N)=\alpha+1$ for some $\alpha$, and $\N|\alpha$ projects to $\om$.
But $\N|\alpha$ is $\GOmega$-$(\om,\om_1+1)$-iterable,
by \ref{IterabilityLifts} and using $\pi$ as in 
\cite[Claim 4.3]{ScalesK(R)}.
The rest is as in \cite{ScalesK(R)}.
\end{proof}
This completes our sketch of the proof.
\end{proof}

\begin{remark}
In the circumstances of the preceding theorem,
if $\M$ has no admissible proper segment, then there is an alternate scale construction.
We include this also, as it yields some extra information. It is 
related to Moschovakis' construction of inductive 
scales on inductive sets.

Let $Q\sub\RR\cross\RR^{<\om}$.
We say that $Q$ is \dfnemph{open} iff $(v,\wvec\conc(x))\in Q$ for 
all $(v,\wvec)\in Q$ and $x\in\RR$. We say that $Q$ is a \dfnemph{basic payoff} iff $Q$ is open,
and definable over $(\HC,\Upsilon^\M,\MFsharp)$.

Let $Q$ be a basic payoff. For $\alpha\leq\omega\cdot\l(\M)$ let
$Q_{<\alpha}=\bigcup_{\beta<\alpha}Q_\beta$, where $Q_0=Q$ and for 
$1\leq\alpha<\omega\cdot\l(\M)$ and 
$\wvec\in\RR^{<\om}$,
\[ (v,\wvec)\in Q_\alpha\iff\mathsf{q}^\RR x[(v,\wvec\conc(x))\in Q_{<\alpha}], \]
where if $\lh(\wvec)$ is even then $\mathsf{q}^\RR=\all^\RR$, and otherwise 
$\mathsf{q}^\RR=\ex^\RR$.
Let $v\in Q'_\alpha$ iff $(v,\emptyset)\in Q_\alpha$, and likewise $Q'_{<\alpha}$.
Let $v\in\RR$. The game $\G^Q_v$ is that where players $\plI$ and $\plII$ 
alternate playing reals $x_0,x_1,\ldots$ (player $\plI$ moving first), and player $\plII$ wins 
iff there is $n<\om$ such that $(v,(x_0,\ldots,x_{n-1}))\in Q$. For $\wvec\in\RR^{<\om}$,
let $\G^Q_{v,\wvec}$ be the game like $\G^Q_v$, except that we interpret $\wvec$ as the first 
$\lh(\wvec)$ moves. Clearly if $(v,\wvec)\in Q_\alpha$
then $\plII$ has a winning quasi-strategy $\G^Q_{v,\wvec}$.

We say $P\sub\RR$ is $\IND^\M$ iff $P=Q'_{<\omega\cdot\l(\M)}$ for some
basic payoff $Q$.

\begin{claim}
$\pow(\RR)\inter\Sigma_1^{\M}=\IND^\M$. 
\end{claim}
\begin{proof}[Proof Sketch]
The fact that $\IND^\M\sub\pow(\RR)\inter\Sigma_1^\M$ is routine.
We now sketch a proof that $\Sigma_1^\M\sub\IND^\M$.
Fix a $\Sigma_1(\Ll^-)$-formula $\Phi$. We define a basic payoff $Q$, implicitly,
by directly defining the corresponding games $\G^Q_v$.
In the definition of the game, some moves are specified as integers (or formulas, etc),
but we take all moves to literally be reals. In some places, one player will play several items 
consecutively, or in a block, but for convenience, we also assume that literally the other 
player plays a dummy real between each consecutive pair of such items.\footnote{We supress these 
dummy reals from the definition of the game as we ignore their values. Their point is that they 
allow us to use the notation $\G^Q_{v,\xvec}$ even when $\xvec$ is a partial play stopping in the 
middle
of some block of items played consecutively by a single player.} At certain points, 
given $d<\om$, we will have a \dfnemph{delay of length $d$},
which is just a string of $d$ alternating moves, whose values will be ignored.\footnote{These 
moves help calibrate the length of inductive computations of winning 
quasi-strategies, as explained later.} We refer to player 
$\plII$ as ``player 
$\ex$'' and player $\plI$ as ``player $\all$''. In $\G^Q_v$, player $\ex$ attempts to prove that 
$\M\sats\Phi(v)$, roughly by describing a strictly descending sequence 
$\left<\M_{n+1}\right>_{n<N}$ of (putative) proper segments of $\M$ and making 
claims about formulas they satisfy.\footnote{In what follows, the (putative) model 
$\M_{n+1}$ is described in round $n$ and is denoted $\N$ in our discussion. If player $\ex$ 
plays according to a winning strategy of simple enough complexity then the models 
$\M_{n+1}$ exist and $\M_{n+1}\pins\M_n$, where $\M_0=\M$.} Player $\all$ keeps player $\ex$ 
honest, by playing reals for 
which player $\ex$ must furnish witnesses to his assertions.
For $\alpha\in[1,\lh(\M)]$, we will get $v\in 
Q'_{<\omega\alpha}$ iff $\M|\alpha\sats\Phi(v)$, thereby proving the claim.

$\G^Q_v$ will be broken into rounds, each of which consists of a finite sequence 
of real moves. Suppose we have a partial play $p$ consisting of $n$ complete rounds,
after which neither player has already won the game. Then $p$ will determine 
a $\Sigma_1(\Ll^-)$-formula $\varphi^n=\varphi^n(p)$ and
$\wvec_n=\wvec_n(p)\in\RR^{<\om}$, where
$\varphi^0=\Phi$ and $\wvec_0=(v)$.
Player $\ex$ will have claimed that $\M\sats\varphi^n(\wvec_n)$.

Let $\varphi\mapsto\left<\varphi_m\right>_{m<\om}$
be the natural recursive function sending $\Sigma_1(\Ll^-)$ formulas $\varphi$ to
sequences $\left<\varphi_m\right>_{m<\om}$ with $\varphi_m\in\Ll^-$,
such that for all $\alpha<\l(\M)$ and $\xvec\in\RR^{<\om}$, we have 
$\M|(\alpha+1)\sats\varphi(\xvec)$
iff there is $m<\om$ such that either $\alpha>0$ and $\M|\alpha\sats\varphi_m(\xvec)$, or
$\alpha=0$ and $(\HC,\Upsilon^\M,\MFsharp)\sats\varphi_m(\xvec)$.

Round $n$ proceeds as follows. Player $\ex$ first plays a code $(m,\psi,z)$ for a 
witness to the claim that $\M\sats\varphi^n(\wvec_n)$, where $m<\om$ and 
$\psi\in\Sigma_1(\Ll^-)\un\{\emptyset\}$ and 
$z\in\RR$,
claiming that $\N\sats\varphi^n_m(\wvec_n)$ where:
\begin{enumerate}[label=--]
\item if $\psi=\emptyset$ then $\N=(\HC,\Upsilon^\M,\MFsharp)$, and
\item if $\psi\neq\emptyset$ then there is $\N'\pins\M$
satisfying $\all^\RR x[\psi(x,z)]$, and $\N$ is the least such $\N'$.
\end{enumerate}

Next, player $\all$ can either 
\emph{dispute} or \emph{accept} the existence of $\N$, where she must accept 
if 
$\psi=\emptyset$.

Suppose $\all$ disputes. Then $\all$ plays $x\in\RR$, then $d<\om$,
which is followed by a delay of length $d$; neither player has yet won.
 Set $\wvec_{n+1}=(x,z)$ and $\varphi^{n+1}=\psi$.

Now suppose $\all$ accepts and $\psi\neq\emptyset$. Let
\[ \all X_0\ex X_1\ldots\all 
X_{2k}\ex 
X_{2k+1}[\varphi^*(\ \cdot\ ,X_0,\ldots,X_{2k+1})] \]
be the prenex normal form of $\varphi^n_m(\ \cdot\ )$,
with $\varphi^*\in\Sigma_1(\Ll^-)$ (here ``$\ \cdot\ $'' represents free variables), and then 
pass in the natural way from $(k,\varphi^*,\psi)$ to a $\Sigma_1(\Ll^-)$ formula 
$\tilde{\varphi}$ such that if $\M\sats\all^\RR x\psi(x,z)$,
then letting $\N\ins\M$ be least satisfying $\all^\RR x\psi(x,z)$, we have
\[ \N\sats\varphi^n_m(\wvec_n)\ \iff\ \N\sats\varrho(\wvec_n,z)\ \iff\ 
\M\sats\varrho(\wvec_n,z) \]
where $\varrho(\ \cdot\ )$ is the formula
\[ \all^\RR x_0\ex^\RR x_1\ldots\all^\RR x_{2k}\ex^\RR 
x_{2k+1}[\tilde{\varphi}(\ \cdot\ ,x_0,\ldots,x_{2k+1})]. \]
Then $\all$ plays $x_0\in\RR$, $\ex$ plays $x_1\in\RR$, etc, producing 
$\xvec=(x_0,\ldots,x_{2k+1})$. Then $\all$ plays $d<\om$, which is followed by a delay of 
length $d$. This completes the round; neither player has 
yet won. Set $\wvec_{n+1}=(\wvec_n,z,\xvec)$ and $\varphi^{n+1}=\tilde{\varphi}$.

Finally suppose that $\psi=\emptyset$, so $\all$ accepts. Pass in the natural way from $\varphi^n_m$
to a $\Sigma_1(\Ll^-)$ formula $\varphi^*$ such that
$(\HC,\Upsilon^\M,\MFsharp)\sats\varphi^n_m(\wvec_n)$ iff
\[ (\HC,\Upsilon^\M,\MFsharp)\sats\all^\RR x_0\ex^\RR x_1\ldots\all^\RR x_{2k}\ex^\RR 
x_{2k+1}[\varphi^*(\wvec_n,x_0,\ldots,x_{2k+1})]. \]
Then $\xvec=(x_0,\ldots,x_{2k+1})$ is played out in the obvious manner.
This finishes the game; $\ex$ 
wins iff $(\HC,\Upsilon^\M,\MFsharp)\sats\varphi^*(\wvec_n,\xvec)$.

This completes the description of round $n$.
We declare $\ex$ the winner iff he wins at some finite stage (in the situation of the previous 
paragraph). This completes the definition of $\G^Q_v$, and hence the implicit definition of $Q$.

\begin{subclaim}
 Let $p$ be a partial play of $\G^Q_v$ consisting of $n$ full rounds, after which neither 
player has yet won. Let $\varphi^n=\varphi^n(p)$ and $\wvec_n=\wvec_n(p)$.
Let $\alpha\in[1,\l(\M)]$. Then:
\begin{enumerate}[label=--]
\item $(v,p)\in Q_{<\om\alpha}$ iff $\M|\alpha\sats\varphi^n(\wvec_n)$.
\item Let $w=(m,\psi,z)$ be a valid move for player $\ex$ in $\G^Q_v$, following 
$p$, and $p'=p\conc w$.
Then $(v,p')\in Q_{<\om\alpha}$ iff either:
\begin{enumerate}[label=--]
\item $\psi=\emptyset$ and $(\HC,\Upsilon^\M,\MFsharp)\sats\varphi^n_m(\wvec_n)$, or
\item $\psi\neq\emptyset$ and there is $\N\pins\M|\alpha$ satisfying $\all^\RR 
x\psi(x,z)$, and the least such $\N$ satisfies $\varphi^n_m(\wvec_n)$.
\end{enumerate}
 \end{enumerate}
\end{subclaim}
 \begin{proof} This is a straightforward induction on $\alpha$, which we omit.\footnote{Let us just 
illustrate how delays help to calibrate the ranks of winning strategies for $\ex$.
Let $p'$ be as above, and adopt the notation there. Suppose that 
$\M|\alpha\sats\all^\RR x\psi(x,z)$, and $\M|\alpha$ is least such. Let 
$p^*=p'\conc\text{``dispute''}$. Since the putative $\N$ does not exist (from the perspective of 
$\M|\alpha$) we want to know that $p^*\notin Q_{<\om\alpha}$.
For $x\in\RR$ and $d<\om$ let $\beta_{x,d}$ be the least $\beta$ such that $p^*\conc(x,d)\in 
Q_\beta$, if such $\beta$ exists. Then $\beta_{x,d}$ does exist,
and $\beta<\om\alpha$, since $\M|\alpha\sats\psi(x,z)$.
Moreover, for each $\gamma<\alpha$, there is $x$ such that $p^*\conc(0,x)\notin 
Q_{<\om\gamma}$ (by the minimality of $\M|\alpha$). But $\sup_{x,d}\beta_{x,d}$ is a limit because 
the arbitrary $d<\om$ is followed by a delay of length $d$ (after which the next round starts), so
$\sup_{d,x}\beta_{d,x}=\om\alpha$.}
\end{proof}

Applying the first conclusion of the subclaim to the case that 
$\alpha=\l(\M)$ and
$p=\emptyset$ (so $n=0$), we have proved the claim.
\end{proof}

Because of the preceding claim, we just need to prove the next one:

\begin{claim}
$\M\sats$``Every $\IND^\M$ set has a $\Sigma_1^\M(a_0)$ scale''.
\end{claim}
\begin{proof}
This is a standard calculation, but here is a sketch. Fix a basic payoff $Q$.
We define a scale on $Q_{<\om\cdot\l(\M)}$ which is $\Sigma_1^\M(a_0)$.

Using the periodicity theorems and determinacy, over $(\HC,\Upsilon^\M,\MFsharp)$ we can define 
from the parameter $a_0$ a very good scale $\vec{\leq}^0$ on $Q$.
(Use $(a_0,\MFsharp)$ to determine the code $a_0'$ for $\MFsharp$ relative to $a_0$,
and from $a_0'$, define a scale on the set $C$ of all codes for $\MFsharp$, and on $\RR\cut 
C$. Then produce scales on Boolean combinations
of $\Upsilon^\M$, $C$ and projective sets first by reducing to the case of \emph{disjoint} 
unions of intersections of $\Upsilon^\M$, $\RR\cut\Upsilon^\M$, $C$, $\RR\cut C$ and projective 
sets.)

Now propagate $\vec{\leq}^0$ to scales on $Q^{<\beta}$, for $\beta\leq\omega\cdot\l(\M)$,
in the usual manner. (For limit $\beta$, 
$\vec{\leq}^{<\beta}=\bigcup_{\gamma<\beta}\vec{\leq}^{<\gamma}$. For each $\beta$, 
$\leq^{<\beta}_0$ is the prewellorder of the norm on $Q^{<\beta}$ given by $x\mapsto\gamma$ where 
$\gamma$ is least such that $x\in Q^\gamma$.
For successor $\beta$, the remaining norms are given by propagating $\vec{\leq}^{<\beta-1}$
using the periodicity theorems, interleaving integer norms in the usual way to yield a very good 
scale.) The propagation process is $\Sigma_1^\M$, so the scale is $\Sigma_1^\M(a_0)$.
\end{proof}
\end{remark}

We now proceed to a variant of \ref{passiveGameScale} we will need, in which $\M$ is 
$\hmP$-active but satisfies ``$\Theta$ does not exist''. (Because $\M$ is a 
$\Theta$-$\g$-spm, this can only happen if $\l(\M)=\alpha+1$ for some $\alpha$ where
$\M|\alpha\sats$``$\Theta$ exists''.)

\begin{definition}
Let $\R=\core_0(\M)$ where $\M$ is an hpm over $A$. Let $\beta<\l(\R)$ and 
$n<\om$ and $H=\Ss_{\omega\beta+n}^{\Pvec^\R}(A)$ be such that
\[ \hmp^\R,\hmPsi^\R,\hmmu^\R,\hme^\R\in H. \]
Define the $\Ll^+$-structure
\[ \R\restpred(\beta,n)=(H,\hmPvec,A;\hmE,\hmP;\hmp^\R,\hmPsi^\R,\hmmu^\R,\hme^\R), \]
where 
$\hmPvec=\hmPvec^\R\inter H$, $\hmE=\hmE^\R\inter H$ and $\hmP=\hmP^\R\inter H$.
\end{definition}

Note that
$\univ{\R|\beta}=\univ{\R\restpred(\beta,0)}$ and 
$\hmPvec^{\R|\beta}=\hmPvec^{\R\restpred(\beta,0)}$, but the $\hmE$ and $\hmP$ predicates of 
$\R|\beta$ and $\R\restpred(\beta,0)$ can differ.

\begin{theorem}
\label{thm:scale_branch_active_immediately_after_new_Sigma_1}
Let $\M$ be a countably iterable $\Theta$-g-organized 
$\Omega$-pm satisfying $\AD$. Assume $\DC_{\RR^\M}$. Suppose $\l(\M)=\beta_0+1$,
and if $\beta_0>0$ then $\M|\beta_0\not\elem_{\rSigma_1(\Ll^-,\RR)}\M$.
Then $\M\sats$``$\rSigma_1^\M(\RR)$ has the scale property''.
\end{theorem}
\begin{proof}
By \ref{passiveGameScale} we may assume that $\M$ is $\hmP$-active.
The proof is given by modifying that of \ref{passiveGameScale} as follows.
We again work with $\HC=\HC^\M$. We assume for simplicity that $\Upsilon^\M=\emptyset$;
otherwise make adaptations as in \ref{passiveGameScale}.

We have $\M^-\ins\M_0\eqdef\M|\beta_0$. Let $\OR(\M_0)=\OR(\M^-)+\lambda_0$ and 
$b_0=b^\M\inter\lambda_0$ and $b^\M=b_0\un(\lambda_0+b_1)$.
(So $b_1\sub\om$. Note that $b_0\in\M_0$.)
Let $z_0\in\RR$, $d<\om$, $k_0\in[1,\om)$, $\varrho_0\in\Ll$, $\psi_0\in\Sigma_1(\Ll^-)$, 
$\Psi_0\in\Ll$ be such that:
\begin{enumerate}[label=--]
\item $z_0\geq_T a_0$;  $a_0$ is computed by the $d^\nth$ Turing machine $\Phi^{z_0}_d$ with oracle 
$z_0$,
 \item $\M\sats\psi_0(z_0)$ and $\M_0\sats\Psi_0(z_0)\&\neg\psi_0(z_0)$, and for all
hpms $\N$, if $z\in\N$ 
and $\N\sats\Psi_0(z)$ then $\J^\hpm(\N)\sats\psi_0(z)$,
 \item $\M_0=\Hull_{k_0}^{\M_0}(\RR)$ and $\M_0\sats b_0=\unique b\varrho_0(z_0,b)$,
 \item $\varrho_0,\Psi_0$ are $\rSigma_{k_0}$ formulas.
\end{enumerate}

Let $\Phi\in\Ll$ be $\Sigma_1$. For $x\in\RR$, let
$A(x)\iff\M\sats\Phi(x)$. We will show that $\M\sats$``There is a
$\Sigma^\M_1(z_0)$-scale on $A$''.
For $x \in \RR$ and 
$k\in[k_0,\om)$ let
\[ A^k(x) \Leftrightarrow \M\restpred(\beta_0,k)\sats\Phi(x).\]
Then $A=\bigcup_{k\geq k_0}A^k$.
We will a construct closed game representation $x \mapsto
G^{k}_x$ for $A^k$, and define $A^k_l$, much as before; player $\plI$ will essentially be 
attempting to build a 
structure $\R\sats\Phi(x)$ and corresponding to $\M\restpred(\beta_0,k)$.
Literally, he will not build the full $\R$ but just a countably iterable $\Theta$-g-organized 
$\Omega$-pm $\N$
corresponding to $\M_0$, with $\N$ satisfying a formula which will ensure that an $\R$ as above
is given by extending $\N$.

We proceed to the details. Let $\L^*,m,n,\sigma_0$ be as before.
Let $(k,c)\mapsto\Phi_{k,c}$ be the natural (and recursive) function with domain 
$[k_0,\om)\cross 2^{<\om}$
such that $\Phi_{k,c}\in\Ll$ is a formula with the following property.
Let $\N$ be any $\om$-sound
$\Theta$-$\g$-spm such that $\N\sats$``$\Theta$ exists'' and
$\Tt\eqdef\Tt_{\varphi_{\mathrm{G}}}^\N$ is defined. Let 
$\lambda=\OR(\tilde{\N})$ where $\tilde{\N}$ is the largest $\varphi_{\mathrm{G}}$-whole initial 
segment of $\N$. Let $y,z\in\RR^\N$ and 
$b\in\N\inter\pow(<\OR(\N))$ and suppose $\N\sats b=\unique b'\varrho_0(z,b')$. 
Let
\[ \hmPvec=\hmPvec^\N\conc(\hmE^\N,\hmP^\N),\]
\[ b^*=(\lambda+b)\un(\OR(\N)+c), \]
\[ \hmP=(\{\Tt\}\cross b^*)\inter\Ss_k^{\hmPvec}(\N), \]
and $\R$ be the $\Ll$-structure
\[ \R=(\Ss_k^{\hmPvec}(\N),\hmPvec,A^\N;\emptyset,\hmP,\hmp^\N,\hmPsi^\N). \]
Then $\R\sats\Phi(x)$ iff $\N\sats\Phi_{k,c}(x,z)$.

For $k\geq k_0$ and $c\in 2^{<\om}$ let $T'_{k,c}$ be the theory given by modifying the theory $T$ 
of \ref{passiveGameScale} by replacing formulas (4) and (9) respectively with (4') and (9') below, 
and adding (10'):
\begin{eqnarray*} (4')&\ 
\Phi_{k,c}(\dot{x}_2,\dot{x}_0)\ \&\ \Psi_0(\dot{x}_0)\ \&\ \neg\psi_0(\dot{x}_0)\\
 (9')&\ \hmpdot\text{ is an hpm over the transitive set coded by }\Phi^{\dot{x}_0}_d\\
 (10')&\ \dot{\Upsilon}=\emptyset\ \&\ V\text{ is }\om\text{-sound }\&\ V=\Hull_{k_0}^V(\RR)
\end{eqnarray*}

In round $n$ of $G^k_x$,
player $\plI$ first plays $i_n,x_{2n},\eta_n$ where $i_n \in \{0,1\}$, $x_{2n} \in
\RR$, $\eta_n<\OR(\M_0)$; player 
$\plII$ plays then $x_{2n+1} \in
\RR$. Define $T^*(u)$, etc, as before. Let $s<\om$ be such that for any transitive 
structure $\N$, $\OR(\Ss(\N))=\OR(\N)+s$. The payoff for player $\plI$ is given by modifying that 
of \ref{passiveGameScale} as follows.
Drop conditions \ref{item:branches_cohere_inf} and \ref{item:branches_correct},
replace condition \ref{item:consistent_extension} with
\begin{enumerate}
\item[(a')] $T^*(p)$ is a consistent extension of $T'_{k,c}$,
where $c=b_1\inter s\cdot k$,
\end{enumerate}
modify conditions \ref{item:correct_reals}, 
\ref{item:embedding} by replacing ``$a_0$'', ``$\rSigma_1$'' and ``$\M|\eta_0$'' respectively with 
``$z_0$'', ``$\rSigma_{k_0+5}$'' and ``$\M_0$'',
and retain the remaining conditions unmodified.

We say that a partial play $u$ of $G^k_x$ is $(k,x)$-\dfnemph{honest} iff
$\M\restpred(\beta_0,k)\sats\Phi(x)$
and if $n>0$ then the modifications of properties
\ref{item:honesty_correct_reals}--\ref{item:honesty_embedding_existence} of 
$(\beta,x)$-honesty of \ref{passiveGameScale}
hold,
given by replacing ``$a_0$'', ``$\M|\eta$'', ``$\M|\eta_0$'', ``$\Upsilon$'' and ``$\rSigma_1$''
respectively with ``$z_0$'', ``$\M_0$'', ``$\M_0$'', ``$\emptyset$'' and ``$\rSigma_{k_0+5}$''.
Let $Q^k_l(x,u)$ iff $u$ is a $(k,x)$-honest position of length $l$.

\begin{claim}
$Q^k_l \in \M$ and the map
$(k,l)\mapsto Q^k_l$ is $\Sigma_1^\M(z_0)$.\end{claim}
\begin{proof}
As before, using that $b_1$ is $\Sigma_1^\M(\{\beta_0\})$
to compute $c=b_1\inter s\cdot k$.
\end{proof}

\begin{claim} $A^k_l=Q^k_l$.
\end{claim}
\begin{proof}
$Q^k_l\sub A^k_l$ as before. For the converse, let $\N$ and 
$\pi:\N\partialto\M_0$ 
be produced as before.
We get $\N=\M_0$ because $\N$ is sound, $\rho_\om^\N=\om$ and 
$\N\sats\Psi(z_0)\&\neg\psi(z_0)$, and $\N$ is
sufficiently iterable above $\Theta^\N$ as $\N=\Hull_{k_0}^\N(\RR)$ and if $\N$ is relevant 
then $\pi$ induces a near 
$k_0$-embedding $\Hh^\N\to\Hh^{\M_0}$. Therefore $b_0$ is the unique $b'\in\N$ such that 
$\N\sats\varrho_0(b',z_0)$. Since $\N\sats\Phi_{k,c}(x,z_0)$ where $c=b_1\inter k\cdot s$,
it follows that
$\M\restpred(\beta_0,k)\sats\Phi(x)$.\end{proof}

This completes the proof.
\end{proof}

\subsection{$\Sigma_1$ gaps}

\begin{definition}
Let $\elem_\RR^-$ abbreviate $\elem_{\rSigma_1(\Ll^-,\RR)}$.
Let $\M$ be an hpm with $\HC^\M\in\M|1$. Let 
$\alpha \leq
\beta \leq\l(\M)$. The interval $[\alpha,\beta]$ is a 
\dfnemph{$\Sigma_1$ gap} of $\M$ iff:
\begin{enumerate}[label=--]
 \item $\M|\alpha\elem_\RR^-\M|\beta$,
 \item $\all\alpha'\in[1,\alpha)$, $\M|\alpha'\not\elem_\RR^-\M|\alpha$, and 
$\all\beta'\in(\beta,\l(\M)]$, $\M|\beta\not\elem_\RR^-\M|\beta'$,
 \item if $\beta=\l(\M)$ then $\M'\eqdef\J^\hpm(\M)$ is an hpm (i.e. $\M$ is $\om$-sound and 
${<\om}$-condensing),
$\HC^{\M'}=\HC^\M$ and $\M\not\elem_\RR^-\M'$.\qedhere
\end{enumerate}
\end{definition}

\begin{definition}\label{Gaps} Let $0<n<\om$. Let $\M$ be an hpm with $\HC^\M\in\M|1$ and $b\in
\core_0(\M)$. The \dfnemph{$\rSigma_n$ type realized by $b$ over $\M$} is
\[
\rSigma^{\M}_{n,b}\eqdef\{\varphi(v)\in \Ll^+ \ | \ \varphi \textrm{ is 
either }
\rSigma_n \textrm{ or } \rPi_n \textrm{ and } \core_0(\M) \vDash
\varphi(b) \}.
\]

Let $[\alpha,\beta]$ be a $\Sigma_1$ gap of $\M$. The gap is \dfnemph{admissible} iff 
$\M|\alpha$ is admissible. The gap is 
\dfnemph{strong}
iff it is admissible and letting $n<\omega$ be least such that
$\rho_n^{\M|\beta}=\om$, every $\rSigma_n$ type realized over
$\M|\beta$ is realized over $\M|\gamma$ for some $\gamma < \beta$. The gap 
is \dfnemph{weak} iff it is admissible but not strong.
\end{definition}

There are no new scales inside the $\Sigma_1$ gaps in which we are interested. The proof of the 
following
theorems are routine generalizations of the corresponding proofs in
\cite{Scales}.
\begin{theorem}[Kechris-Solovay]
\label{InsideGap}
Let $\M$ be a $\Theta$-$\g$-organized $\Omega$-pm satisfying $\AD$.
Assume $\DC_{\RR^\M}$ and that $\M$ is countably $(0,\om_1+1)$-iterable. 
Let $[\alpha,\beta]$ be a $\Sigma_1$ gap of $\M$. Then:
\begin{enumerate}
\item\label{item:no_uniformization} There is a $\Pi_1^{\M|\alpha}$ subset of 
$\RR^\M\cross\RR^\M$ not uniformized in $\M|\beta$.
\item Let $\alpha\leq \gamma < \beta$ and $1\leq n<\om$, and either let 
$\Gamma=\bfrPi_n^{\M|\gamma}$ or suppose $(\alpha,1)<_\lex(\gamma,n)$ and let 
$\Gamma=\bfrSigma_n^{\M|\gamma}$. Then $\M\sats$``$\Gamma$ does not have the scale 
property''.
\end{enumerate}
\end{theorem}
A relation witnesing \ref{InsideGap}(\ref{item:no_uniformization}) is 
$(\RR^\M)^2\cut\mathcal{C}^{\M|\alpha}$ where 
$\mathcal{C}^{\M|\alpha}(x,y)$ iff $x,y\in\RR^\M$ and there is $\gamma < \alpha$ such that $y$ 
is definable over $\M|\gamma$ from parameters in 
$\Ord\un\{x\}$. The same relation witnesses that there is no new scale 
definable over the end of
a strong gap:
\begin{theorem}[Martin]
\label{StrongGap}
Let $\M,[\alpha,\beta]$ be as in \ref{InsideGap}.
Suppose that $\beta<\l(\M)$ and $[\alpha,\beta]$ is a strong gap of $\M$. Then:
\begin{enumerate}
\item There is a $\Pi_1^{\M|\alpha}$ subset of $\RR^\M\cross\RR^\M$ not uniformized in  
$\M|\beta+1$.
\item Let $n<\om$, and either let $\Gamma=\bfrPi_n^{\M|\beta}$
or suppose $(\alpha,1)<_\lex(\beta,n)$ and let $\Gamma=\bfrSigma_n^{\M|\beta}$.
Then $\M\sats$``$\Gamma$ does not have the scale property''.
\end{enumerate}
\end{theorem}

\begin{remark}\label{rem:non_adm_gap_dealt_with}
The only case remaining in the analysis of scales in 
$\Lp^{\GOmega}(\RR,\Upsilon)$, where $\Upsilon$ is self-scaled, is at the end of a weak gap. For 
let 
$\M$ be a $\Theta$-g-organized 
$\Omega$-pm and let $[\alpha,\beta]$ be a gap of $\M$. Suppose $[\alpha,\beta]$ 
is inadmissible. Then $\alpha=\beta$ and $\M|\alpha\sats$``$\Theta$ does not 
exist''. Note then $\M\sats$``$\bfrSigma_1^{\M|\alpha}$ has the scale property'',
by \ref{passiveGameScale} and 
\ref{thm:scale_branch_active_immediately_after_new_Sigma_1}.\footnote{\label{ftn:motivate_Theta-g}
It is important here that our 
structure is $\Theta$-$\g$-organized, as opposed to $\g$-organized,
since $\g$-organized structures can satisfy ``$\Theta$ does not exist'', be of limit length, and be 
$\hmP$-active. We do not see how to generalize the proof of 
\ref{thm:scale_branch_active_immediately_after_new_Sigma_1} to deal with this 
case.} 
Combined with the argument in \cite{Scales}, this ensures that 
$\J(\M|\alpha)\sats$``Every set of reals has a scale'', assuming that 
$\RR^{\J(\M|\alpha)}=\RR^\M$ and $\J(\M|\alpha)\sats\AD$. The ends of strong gaps have just 
been dealt with, so we are left with weak gaps. We deal with weak gaps in three cases, as described 
in the 
introduction.
\end{remark}

\subsection{Scales at the end of a weak gap from strong determinacy}
The first scale construction for weak gaps proceeds from a strong determinacy assumption. It is 
most useful for weak gaps $[\alpha,\beta]$ of $\Lp^{\GOmega}(\RR,\Upsilon)$ where 
$\Omega\rest\HC\notin\Lp^{\GOmega}(\RR,\Upsilon)|\alpha$.

\begin{theorem}
\label{gapInTheMouse}
Let $\R$ be a $\Theta$-g-organized $\Omega$-pm satisfying $\AD$.
Assume $\DC_{\RR^\M}$ and that $\R$ is countably $\GOmega$-$(0,\om_1+1)$-iterable. 
Let 
$[\alpha,\beta]$ be a weak gap of $\R$ with $\beta<\l(\R)$. Let 
$n+1<\om$ be least such that $\rho_{n+1}^{\R|\beta}=\om$. Then $\R\sats$
``$\bfrSigma_{n+1}^{\R|\beta}$ has the scale 
property''.
\end{theorem}

\begin{proof}[Proof Sketch]
The proof is almost that of \cite[Theorem 4.16]{ScalesK(R)},
so we only sketch it. However, our approach is a little different from that used in 
\cite{ScalesK(R)}.\footnote{This is because the authors do not see, in the proof of \cite[Claim 
4.18]{ScalesK(R)}, and in the notation of that 
proof, why $\N=\M$, because it is not clear that $\N$ is sound. Our approach gets around this 
problem, and also simplifies the proof, because it eliminates the need for the 
``bounding integers'' $m_k$ and $n_k$ played by player $\plI$ in the game $G^i_x$ of 
\cite{ScalesK(R)}.}
For simplicity, we
assume that $\Upsilon^\R=\emptyset$ and $n=0$ and $\beta$ is a limit ordinal. (If 
$\Upsilon^\R\neq\emptyset$ 
make changes as in the proof of \ref{passiveGameScale}.) Let $\M=\R|\beta$.

Let $p =p_1^\M$ and let $w_1\in\RR^\M$ be such 
that $w_1\geq_T a_0$ and the solidity witness(es) $W$ for $p$ 
is in $\Hull_1^\M(p,w_1)$ and $\Sigma\eqdef\rSigma_{1,(p,w_1)}^{\M}$ is a
non-reflecting type. Let $\M^\restpred_\gamma$ denote
$\M\restpred(\gamma,0)$.\footnote{In \cite{ScalesK(R)}, this is denoted $\M||\gamma$.}.
We now define a sequence 
$\langle \beta_i,Y_i,\psi_i,\xi_i  \rangle_{i<\om}$ by recursion on $i$, as follows:
\begin{eqnarray*}\beta_0&=&\text{ least 
}\gamma>\hmmu^\M\text{ such that }\max(p)<\OR(\M^\restpred_\gamma),\\
Y_i&=&\Hull_\om^{\M^\restpred_{\beta_i}}(
\RR^\M\cup\{p\}),\\
\psi_i&=&\text{ least } \psi\in \Sigma\text{ 
such that 
}\M^\restpred_{\beta_i} \vDash \neg
\psi((p,w_1)),\end{eqnarray*}
and then if $\M$ is either $\hmE$-passive or $\hmE$-active type 3, let $\xi_i=0$ and
\begin{eqnarray*}
\beta_{i+1}&=&\text{ least }\gamma\text{ 
such that 
}\M^\restpred_{\gamma} \vDash \psi_i((p,w_1)),\end{eqnarray*}
and otherwise ($\M$ is $\hmE$-active type 1 or 2),
let\footnote{Recall that $E$ 
is the $\M$-amenable predicate coding the active extender of $\M$.}
\begin{eqnarray*}
\xi_i&=&\sup(Y_i\cap((\hmmu^\M)^+)^\M),\\
\beta_{i+1}&=&\text{ least }\gamma\text{ such that }\M^\restpred_{\gamma}\sats\psi_i((p,w_1))\text{ 
and }\\
&&E^\M\cap \M^\restpred_\gamma\text{ measures all sets in 
}\M|\xi_i.\end{eqnarray*}

\begin{claim}\label{claim:union} $\bigcup_{i<\omega} Y_i =\core_0(\M)$. In particular, 
$\l(\core_0(\M))=\lim_{i<\om}\beta_i$.
\end{claim}
\begin{proof}
Let $\N$ be
the transitive collapse of $\bigcup_{i<\omega} Y_i$ and let $\pi\maps\N
\rightarrow \bigcup_{i<\omega} Y_i$ be the uncollapse map. Let $\beta_\om=\sup_{i<\om}\beta_i$. 
Note that $\M^\restpred_{\beta_\om}\sats\Sigma$ and 
$H\sub\rg(\pi)$ where $H=\Hull_1^{\M}(\{p,w_1\})$, $\beta_i\in H$, $\pi$ is 
$\Sigma_1$-elementary on $\pi^{-1}``H$, and the latter is $\in$-cofinal in $\N$.\footnote{So
$\Th_1^\M(\{\beta_0,\beta_1,\ldots\})$ is recorded in $\Sigma$;
it would not have made any difference to add the parameter $\beta_i$ to $Y_{i+1}$.} In particular, 
$\pi$ is a weak $0$-embedding.
So essentially by 
\ref{lem:Th-g_very_con}, $\N$ is a $\Theta$-g-organized $\Omega$-pm, and 
clearly $\HC^\N=\HC^\M$.

Let $\pi(p^*)=p$. It is easy to see that $\N=\Hull_1^\N(\RR^\N\un\{p^*\})$. But $\pi^{-1}(W)$ is 
a generalized solidity witness for $p^*$.\footnote{This only uses the 
$\Sigma_0$-elementarity of $\pi$.
Actually $W\in H$, so $\pi$ is even $\Sigma_1$-elementary on $\pi^{-1}(W)$.
But we would in general need $\Sigma_2$-elementarity to infer already that 
$\pi^{-1}(W)$ is the \emph{standard} solidity witness for $p^*$.}
So $\N$ is $(1,p^*)$-solid. Therefore $\N$ is $1$-sound and 
$p_1^\N=p^*$. 
Since trees on $\N$ can
be lifted to trees on $\M$ via $\pi$, $\N$ is countably
${\GOmega}$-$(0,\om_1+1)$--iterable. Since $\N$ is also minimal realizing $\Sigma$, therefore 
$\N 
= 
\M$.

The fact that $\pi=\id$ now follows as usual, using the fact that $p^*=p$.
\end{proof}

Using notation mostly as in the proof of \cite[Theorem 
4.16]{ScalesK(R)}, we define the game $G^k_x$ mostly as there, with some modifications. Player 
$\plI$ describes his model 
using the language $\Ll^*=\Ll^+\cup 
\{\dot{x}_i, \dot{\beta}_i, \dot{\M}_i\}_{i<\omega}\cup \{\dot{G},\dot{p},\dot{W}\}$; 
the 
symbols in $\Ll^*\cut\Ll^+$ are constants. Let $B_0$ be defined from $\Ll^*$ as in 
\cite{ScalesK(R)}.\footnote{That is, in the manner that $B_0$ is defined from the $\Ll$ of 
\cite{ScalesK(R)}. The symbols $\Ll$ and $\Ll^*$ have had their roles 
interchanged from \cite{ScalesK(R)}.} Let $S_0$ be the set of sentences $\varphi\in B_0$ such 
that $i\in\{1,2\}$ whenever $\dot{x}_i$ appears in $\varphi$, and 
$(\core_0(\M),I)\sats\varphi$ where $I$ is the assignment
\[ (\dot{x_1},\dot{x_2},\dot{G},\dot{p},\dot{W},\left<\dot{\beta}_i,\dot{\M}_i\right>_{i<\om})^I=
(w_1,w_2,p,p,W,\left<\beta_i,\M^\restpred_{\beta_i}\right>_{i<\om}).\]
A run of $G^k_x$ has the form
\begin{IEEEeqnarray*}{cccccccc}
\plI  & \ \ \ & T_0,s_0,\eta_0 & & T_1,s_1,\eta_1 & & \cdots\\
\plII &\ \ \ & &\ s_1\ & &\ s_3\ & & \cdots
\end{IEEEeqnarray*}
where $T_i,s_i$ are as in \cite{ScalesK(R)} and $\eta_i\in\OR(\M)$ .
The winning conditions for 
player $\plI$ are the winning conditions 
(1)--(6)\footnote{We have no need for the integer moves $m_k$, nor any 
version of 
condition (8) used in 
\cite{ScalesK(R)}.} of \cite{ScalesK(R)} verbatim (other than a small notational difference), 
and \dfnemph{$(k,x)$-honesty} is as in
\cite{ScalesK(R)} except that we drop condition (iv) from there. Define $A^k_l$ (strategic) and 
$Q^k_l$ (honest) in the obvious manner (the analogue of $A^k_l$ was denoted $P^k_l$ in 
\cite{ScalesK(R)}).
\begin{claim}
$A^k_l=Q^k_l$.
\end{claim}
\begin{proof}[Proof Sketch]
Consider the proof that every strategic position is honest. We use notation mostly as in 
the proof of \cite[Claim 4.19]{ScalesK(R)}, with a couple of changes. 
Let $\N$ be the reduct of $\Aa$ to an $\Ll^+$-structure. Let $\N_i$ be (the $\Ll^+$-structure) 
$\dot{\M}^{\mathcal{A}}_i$. Because $\Aa\sats S_0$, 
$\N_i=\N^\restpred_{\beta_i^*}$ and $\N$ is the ``union'' of the $\N_i$. Let $p^*=\dot{p}^\Aa=G^*$. 
As in the proof of \cite[Claim 4.19]{ScalesK(R)} we get that $\N$ is a countably 
$\GOmega$-$(0,\om_1+1)$-iterable 
$\Theta$-g-organized $\Omega$-pm which is minimal for realizing $\Sigma$. Clearly 
$\Upsilon^\N=\emptyset=\Upsilon^\M$. Also, 
$\N$ is sound with $\rho_1^\N=\RR^\N$ and $p_1^\N=p^*$. For let $H=\Hull_1^\N(\RR^\N\un p^*)$. Then 
because $\Aa\sats S_0$, we have:
\begin{enumerate}[label=--]
 \item $\N_i\in H$ for each $i$ (it follows that $H=\univ{\N}$),
\item $W^*$ is a generalized solidity witness for $p^*$  (so $\N=\M$ and 
$p^*=p$),
\item $W^*=W$, $\beta_i^* = \beta_i$ and $\N_i = \M_{\beta_i}^\restpred$ for all $i$.\qedhere
\end{enumerate}
\end{proof}
\begin{claim}
$Q^k_l\in\M$ for all $k,l$, and the map $(k,l)\mapsto Q^k_l$ is $\rSigma_1^\M(p,w_1,w_2)$.
\end{claim}
\begin{proof}[Proof sketch]
The proof is the same as that of \cite[Claim 4.20]{ScalesK(R)} (except that condition (iv) of 
\cite{ScalesK(R)} is not involved, so the use of the Coding Lemma regarding this condition is 
avoided). In the computation of the definability of (v) we still use the Coding Lemma; it is here 
that we use our assumption that $\J_1(\M)\sats\AD$ (beyond that $\M\sats\AD$).
\end{proof}
The remaining details are as in \cite{ScalesK(R)}.
\end{proof}

\subsection{Scales at the end of a weak gap from optimal determinacy}
As described in \cite{Scalesweakgap}, typically in the core model induction, one does not have the 
stronger determinacy hypothesis at the stage required to apply \ref{gapInTheMouse}. So we need
generalizations of \cite[Theorem 4.17]{ScalesK(R)} and \cite[Theorem 0.1]{Scalesweakgap}, which are 
the second and third cases of our scale constructions for weak gaps, respectively.

\begin{definition}
 Let $\M$ be a $\Theta$-$\g$-organized $\Omega$-pm.
We say $\M$ is \dfnemph{subtle} iff $\M\sats$``$\Theta$ exists'' and
either $\M$ is $\hmP$-active or there is an $\M$-total $E\in\es_+^\M$.
We say 
$\M$ is \dfnemph{self-analysed} iff for every subtle $\N\ins\M$ there is 
$\P\ins\M$ such that $\N\pins\P$ and $\P$ is admissible.
We say 
$\M$ is \dfnemph{self-coded} iff for every subtle 
$\N\ins\M$ there is $\P\pins\M$ such that $\N\ins\P$ and $\rho_\om^\P=\om$.
\end{definition}

Note that if $\M\sats$``$\Theta$ does not exist'' or $\M$ has no 
active segment above $\Theta^\M$ then
$\M$ is self-coded.

\begin{theorem}\label{thm:self_anal_coded}
Let $\M$ be a $\Theta$-g-organized $\Omega$-pm satisfying $\AD$.
Assume $\DC_{\RR^\M}$ and that every proper segment of $\M$ is countably 
$\GOmega$-$(\om,\om_1+1)$-iterable.
Suppose that $\M$ ends a weak gap of $\M$,
and $\M$ is either self-analysed or self-coded.
Let $n<\om$ be least 
such that $\rho_{n+1}^\M=\om$. Then $\M\sats$``$\bfrSigma_{n+1}^\M$ has the 
scale property''.
\end{theorem}
\begin{proof}[Proof Sketch]
 The proof is similar to that of \ref{gapInTheMouse}, but we use the fact that $\M$ is either 
self-analysed or self-coded to reduce the reliance on determinacy.\footnote{Of course determinacy 
is still required in the, supressed, norm propagation part of the argument.}

Suppose first that $\M$ is passive. We assume for simplicity that $\Upsilon^\M=\emptyset$, 
$\l(\M)$ is a limit and $n=0$. We define most things, including $Y_k$ and $B_k$, 
as in the 
proof of \ref{gapInTheMouse}. Fix 
$x\in\RR$ and $i<\om$; we want to define the game $G^i_x$. Let $m:B_0\cross B_0\to\om$ and 
$n:B_0\to\om$ be recursive and injective with disjoint ranges, and such that for all 
$\varphi,\psi\in B_0$, $\varphi,\psi$ have support $m(\varphi,\psi)$ and $\varphi$ has support 
$n(\varphi)$ and if $\varphi\neq\psi$ then $m(\varphi,\varphi)<m(\varphi,\psi)$. A run of $G^i_x$ 
consists of the same types of objects as in the proof of \ref{gapInTheMouse}, except that we 
also require that $\eta_k\in Y_k$. The rules of $G^i_x$ are (1)--(5) as stated in 
\cite{ScalesK(R)}, along with rule \ref{item:order_pres} below, which requires player $\plI$ to play 
a wellfounded model, and rule \ref{item:subtle_elem} below, which requires player $\plI$ to 
build, for each subtle initial segment $\P$ of his model, a partial embedding $\P\to\R$  for 
some $\R\ins\M$, which is elementary on ordinal parameters (but these embeddings need not agree with 
one another):

\begin{enumerate}[label=(\arabic*)]
\setcounter{enumi}{5}
\item\label{item:order_pres} if $\varphi,\psi\in B_0$ each have one free 
variable and
\[ \text{``}\unique v\varphi(v) \in \Ord\ \&\ \unique v\psi(v) \in 
\Ord\text{''}\in
T^*,\]
then ``$\unique v\varphi(v) \leq \unique v\psi(v)$''$\in T^*$ iff
$\eta_{n(\varphi)} \leq \eta_{n(\psi)}$,
\item\label{item:subtle_elem} if $\psi,\sigma_0,\ldots,\sigma_{j-1}\in B_0$ each have one free 
variable and $k<\om$ and
\[ \text{``}\unique v\psi(v)<\l(\dot{\M}_k)\ \&\ \dot{\M}_k|(\unique 
v\psi(v))
\text{ is subtle}\text{''}\in T^* \]
and for all $i<j$,
\[ \text{``}\unique v\sigma_i(v)\in\OR(\dot{\M_k}|(\unique 
v\psi(v)))\text{''}\in 
T^*
\]
then $\eta_{m(\psi,\psi)}<\l(\M_k)$
and for any $\Ll$-formula $\theta(v_0,\dots,v_{j-1},u)$,
\[
\text{``}\dot{\M_k}|(\unique v\psi(v))\sats
\theta(\unique v\sigma_0(v),\dots,\unique v\sigma_{j-1}(v),\dot{x}_1)\text{''}\in
T^* \]
if and only if
\[
\M|\eta_{m(\psi,\psi)}\sats
\theta(\eta_{m(\psi,\sigma_0)},\dots,\eta_{
m(\psi,\sigma_{j-1})},w_1).
\]
\end{enumerate}

We omit most of the remaining details, including the precise 
formulation of \emph{$x$-honesty} (of a position in $G^i_x$). The analysis of 
commitments made pertaining to rule \ref{item:order_pres} are dealt with as in \cite{Scales}. 
Consider rule \ref{item:subtle_elem}. If $\M$ is self-analysed then the analogue of condition (v) 
of \emph{$x$-honest} from \cite{ScalesK(R)} can be computed in some admissible proper segment of 
$\M$ (without the Coding Lemma). Suppose $\M$ is self-coded but not self-analysed.
Then there is $\R\pins\M$ such that $\rho_\om^\R=\om$ and every subtle initial 
segment of $\M$ is a segment of $\R$. One can therefore use the Coding 
Lemma as in the 
proof of Claim \ref{clm:honesty_computable} to compute the analogue of 
condition (v) over $\R$. In rule \ref{item:subtle_elem} we have required elementarity with 
respect to $w_1$ (and ordinals) just to ensure elementarity with respect to $a_0$ (and ordinals).

This completes a sketch of the proof in the passive case. Now suppose that $\M$ is 
active. The scale construction in this case combines elements 
of \ref{thm:scale_branch_active_immediately_after_new_Sigma_1} and 
\ref{gapInTheMouse}, and we just outline what is new. Since $\M$ is not subtle, 
$\M\sats$``$\Theta$ does not 
exist'' and $\M$ is $\hmP$-active, so because $\M$ is $\Theta$-$\g$-organized,
$\l(\M)=\beta_0+1$ for some $\beta_0>0$, and $\M|\beta_0\sats$``$\Theta$ exists'' and 
$\rho_\om^{\M|\beta_0}=\om$. Therefore $n=0$. Let 
$\Tt=\Tt^\M$. Assume $\Upsilon^\M=\emptyset$, and also that $\lh(\Tt^\M)>\om$; the case 
that $\lh(\Tt)=\om$ is simpler, partly because then $\Tt$ is linear, as $\M$ is 
$\Theta$-$\g$-organized. Let $\M_0=\M|\beta_0$. Note that $\M_0,\beta_0\in\Hull_1^\M(\emptyset)$.
Let $k_0\in[1,\om)$ be such that $\M_0=\Hull_{k_0}^{\M_0}(\RR)$.
Let $\lambda_0^\M$ be the limit ordinal such that 
$\lh(\Tt^\M)=\lambda_0^\M+\om$.
Let $b_0^\M=b^\M\inter\lambda_0^\M$.
Let $p,w_1,\Sigma$ be as 
usual, except with the added requirement that $b_0^\M\in\Hull_{k_0}^{\M_0}(\{w_1\})$.
In $G^i_x$, player $\plI$ is required to build a $\Theta$-$\g$-spm $\N$ with $w_1\in\RR^\N$
and with $p^*\in\N$ such that $\rSigma_{1,(p^*,w_1)}^\N=\Sigma$ and 
$\N\restpred(\beta_0^*,i)\sats\Phi(x)$, where $\l(\N)=\beta^*_0+1$,
and letting $\N_0=\N|\beta_0^*$, is also required to build a partial embedding 
$\pi:\N_0\partialto\M_0$, with domain $\OR(\N_0)\un\{w_1\}$, such that $\pi$ is 
$\Sigma_{k_0+5}$-elementary on its domain. We leave to the reader the precise formulation of 
$G^i_x$, and of honesty.

Because player $\plI$ is only required to embed $\N_0\partialto\M_0$,
and $\rho_{k_0}^{\M_0}=\om$, the Coding Lemma argument shows that honesty is sufficiently 
simply computable.
The fact that ``strategic (for player 
$\plI$) implies honest'' is as follows. Let $\N$ and $\pi:\N_0\partialto\M_0$ be produced by a 
generic run against a winning strategy for player $\plI$, as usual. Then $\N=\M$. For we have 
$\N_0\pins\M$ as usual.
Since $\rSigma^\N_{1,(p^*,w_1)}=\Sigma$, it therefore suffices to see that 
$b^\N=\Lambda_{\MFsharp}(\Tt^\N)$.

We claim that $\pi$
induces a hull embedding $(\Tt^\N\conc b^\N)\to(\Tt^\M\conc b^\M)$, which suffices. For clearly
$\pi$ induces a hull embedding $\Tt^\N\to\Tt^\M$. Let $\lambda_0^\N,b_0^\N$ be defined over $\N$, 
analogously to $\lambda_0^\M,b_0^\M$ over $\M$.
Let $\varrho_0$ be an $\rSigma_{k_0}$ formula such that $b_0^\M=(\unique 
b\varrho_0(b,w_1))^\M$. Since $\rSigma^\N_{1,(p^*,w_1)}=\Sigma$,
then $b_0^\N=(\unique b\varrho_0(b,w_1))^\N$.
But let $c^\M=b^\M\inter[\lambda_0^\M,\lh(\Tt^\M))$ and $c^\N$ 
likewise for $\N$. Then $c^\M=c^\N$ because of how they are determined by $\Sigma$.
The claim easily follows.
\end{proof}

We now proceed to the generalization of \cite[Theorem 0.1]{Scalesweakgap}, the final scale 
construction of the paper. While it uses only the weaker 
determinacy 
assumption, it requires a mouse capturing hypothesis, as in \cite{Scalesweakgap}.

\begin{definition}
Suppose $V$ is an hpm and $\HC$ exists. Let $\Gamma$ be a pointclass of the form 
$\rSigma_1^{V|\alpha}\inter\pow(\RR)$ 
for some $\alpha<\l(V)$. In this setting, for $x\in\RR$, we write $C_\Gamma(x)$ for the set 
of all $y\in\RR$ such that for some ordinal $\gamma<\om_1$, $y$ (as a subset of $\om$) is 
$\Delta_\Gamma(\{\gamma,x\})$.
Let $x\in\HC$ be such that $x$ is transitive and $f:\om\onto x$. Then $c_f\in\RR$ 
denotes the code for $(x,\in)$ determined by $f$.
And $C_\Gamma(x)$ denotes the set of all 
$y\in\pow(x)$ such that for all $f:\om\onto x$ we have $f^{-1}(y)\in 
C_\Gamma(c_f)$.
\end{definition}

\begin{lemma}\label{lem:C_Gamma_characterize} Let $\P$ be a $\Theta$-g-organized $\Omega$-pm 
satisfying $\AD$.
Let $\Q\pins\P$ be such that $\Q$ is passive and admissible. Work in $\P$. Let $\Gamma$ be the 
pointclass $\rSigma_1^\Q\inter\pow(\RR)$. Let $x\in\HC$ with $x$ transitive and infinite. Then 
for all $y\in\HC$, the following are equivalent:
\begin{enumerate}[label=\tu{(}\arabic*\tu{)}]
 \item\label{item:y_in_C_G(x)} $y\in C_\Gamma(x)$,
 \item\label{item:y_OD^alpha(x)} there is $\R\pins\Q$ such that $y$ is definable over 
$\R$ from parameters in $\Ord\un x\un\{x\}$,
 \item\label{item:comeager} for comeager many bijections $f:\om\to x$, $f^{-1}(y)\in 
C_\Gamma(c_f)$.
\end{enumerate}
\end{lemma}
\begin{proof}
 The proof is mostly like that of \cite[Theorem 3.4(?)]{twms}; we just mention a couple of points. 
For $x\in\RR$, the equivalence of \ref{item:y_in_C_G(x)} and \ref{item:y_OD^alpha(x)} follows 
because $\Q\sats\AD+\KP$. Now consider the proof 
that \ref{item:comeager} implies \ref{item:y_OD^alpha(x)}. If $\P$ satisfies 
\ref{item:comeager}, then we may take the witnessing comeager set $C$ to be a 
countable intersection of dense sets, and then $C\in\Q$. So by 
$\KP$ there is $\R\pins\Q$ such that for every $f\in C$, $f^{-1}(y)$ is definable over $\R$ 
from parameters in $\Ord\un\{c_f\}$. As in \cite{twms}, there is then some $\alpha<\om_1^\P$ and 
$n<\om$ and injection $\sigma:n\to x$ such that for comeager many bijections $f:\om\to x$ extending 
$\sigma$, $f^{-1}(y)$ is the $\alpha^\nth$ real which is definable over $\R$ from parameters in 
$\Ord\un\{c_f\}$, in the natural ordering. Letting $\delta=\l(\R)$, this defines $y$ over 
$\Q|(\delta+2)$ from parameters in $\{\delta,x\}\un\rg(\sigma)$.
\end{proof}

\begin{definition}
Let $\P,\Q,\Gamma,x$ be as in \ref{lem:C_Gamma_characterize}.
Suppose that $\MFsharp\in\J(\hat{x})$ and $\Omega^*\in\Q$ where 
$\Omega^*=\Omega\rest\HC^\P$.
Work in $\P$. Then
$\Lp^{\Gamma,{\gOmega^*}}(x)$ denotes 
$(\Lp^{\gOmega^*}(x))^\Q$.\footnote{So $\Q\sats$``$\N$ is $\gOmega^*$-$(\om,\om_1+1)$-iterable''
for all $\N\pins\Lp^{\Gamma,\gOmega^*}(x)$. Note here that $\Q\sats$``$\pow(\om_1)$ exists'' 
because $\Q\sats\AD$.} Similarly for $\Lp_+^{\Gamma,\gOmega^*}(x)$.
We say that \dfnemph{super-small $\Gamma$-${\gOmega^*}$-mouse capturing
holds on a cone} 
iff there is $z\in\RR$ such that for all transitive $x\in\HC$, if
$\MFsharp,z\in\J(\hat{x})$ then $\Lp^{\Gamma,{\gOmega^*}}(x)$ is 
super-small and
\[ C_\Gamma(x)=\Lp^{\Gamma,{^\g\Omega^*}}(x)\inter\pow(x).\qedhere \]
\end{definition}

\begin{theorem}
\label{weak gap from MC}
Let $\M$ be a fully sound, $\Theta$-$\g$-organized $\Omega$-pm
satisfying $\AD$. Suppose $[\alpha_0,\l(\M)]$ is a weak gap of $\M$ and that $\M$ is countably 
$\gOmega$-$(n,\om_1+1)$-iterable where $n<\om$ is least such that $\rho_{n+1}^\M=\om$. 
Assume $\DC_{\RR^\M}$ and
$\RR^{\J(\M)}=\RR^\M$ and $\J(\M)\sats\DC_\RR$.\footnote{$\J(\M)$
provides a universe in which we can execute certain arguments in the proof of \cite[Theorem 
0.1]{Scalesweakgap} without introducing new reals. The authors believe that \cite[Theorem 
0.1]{Scalesweakgap} should also have adopted a hypothesis along these lines. Indeed, its proof 
seems to proceed under the implicit assumption that $\RR^\M=\RR^V$.}
Suppose that $\Omega^*\in\M|\alpha_0$ where $\Omega^*=\Omega\rest\HC^\M$.
In $\M$, let $\Gamma$ be the 
pointclass $\rSigma_1^{\M|\alpha_0}\inter\pow(\RR)$, and assume that super-small 
$\Gamma$-${\gOmega^*}$-mouse 
capturing holds on a cone. Then
$\M\sats$``$\bfrSigma_{n+1}^\M$ has the scale property''.
\end{theorem}

\begin{proof}
 We follow the proof of \cite{Scalesweakgap}, making some modifications.
 By $\DC_{\RR^\M}$ we may assume that $\M$ is countable.
 By \ref{thm:self_anal_coded} we may assume that 
$\M\sats$``$\Theta$ exists'' and there is some $\xi+1\in(\Theta^\M,\l(\M))$ such that $\M|\xi\vDash 
\ZF$. Therefore $\pow(\RR)\inter\M\in\M|\xi$ and $\M|\xi\sats\ZF+\AD$.
We work mostly inside $\J(\M)$, and so we write $\RR=\RR^\M$, $\HC=\HC^\M$, etc.
We have $\Omega^*\in\M|\alpha_0$. Let $z_0\in\RR$ be in the mouse capturing cone, with 
$z_0\geq_T (a_0,t)$ where $t$ codes $\Th_1^{\MFsharp}$ relative to $a_0$, and such that 
$\{\Omega^*\}$ is $\rSigma_1^{\M|\alpha_0}(z_0)$.
For this proof, except where context dictates otherwise, \dfnemph{premouse} 
abbreviates \emph{$\g$-organized $\Omega^*$-pm over $(\N,x)$ for some $x\geq_T z_0$ and 
transitive structure $\N$ with $\MFsharp\in\J(\N,x)$}; likewise all 
related 
terminology (such as \dfnemph{iteration tree}, \dfnemph{iterability}, $\Lp$, etc).

Because $[\alpha_0,\l(\M)]$ is a $\Sigma_1$ gap of $\M$,  for $(\N,x)$ as above we have 
(with terminology as 
just described above)
\[ \Lp^\Gamma(\N,x)=\Lp(\N,x)^{\M|\alpha_0}=\Lp(\N,x)^{\M}.\]
Likewise for $\Lp_+^\Gamma(\N,x)$.

\begin{remark}\label{rem:translation_Theta-g_above_cutpoint}
Let $\N$ be a $\g$-whole premouse and $\N\pins\Q\ins\Lp^\Gamma_+(\N)$ with $\Q$ projecting to $\N$.
Then $\Q$ translates 
to some $\Q'\pins\Lp^\Gamma(\N)$, where $\OR(\Q')=\OR(\Q)$ and $\Q'$ projects to $\om$, as follows.
There is a slight wrinkle in the translation, because we must have 
$\hmPsi^{\Q'}=\emptyset$ as $\hmPsi^{\Lp^\Gamma(\N)}=\emptyset$, whereas possibly 
$\Sigma^\N\neq\emptyset$.
We have $\univ{\Q'}=\univ{\Q}$ and $\hmb^{\Q'}=\hat{\N}$.  There is $\alpha>0$ such that 
$\l(\N)+\alpha\leq\l(\Q)$ and for $\beta>\alpha$, $\Q|(\l(\N)+\beta)$ and $\Q'|\beta$ have the 
same active predicates, and for $\beta\in[1,\alpha]$, $\R\eqdef\Q|(\l(\N)+\beta)$ and 
$\R'\eqdef\Q'|\beta$ are both $\hmE$-passive,
and if $\R$ is $\hmP$-active then $\Tt^\R$ is linear, and if $\R'$ is $\hmP$-active then 
$\Tt^{\R'}$ is linear. These are linear iterations at the 
least measurable of $\MFsharp$. Because the iterations are linear, the corresponding predicates are 
trivial, so we can trivially translate between them.\footnote{$\R$ and $\R'$ can
have different predicates, because the definition of spm requires that a particular tree can only 
have a cofinal branch added at at most one segment of the spm. We must have 
$\Sigma^{\Q'|1}=\emptyset$ by definition, but possibly $\Sigma^\N\neq\emptyset$, in which case there 
can be conflict between $\R,\R'$ over which tree should have a branch added. But it is easy to see 
that if $\Q$ is large enough then $\Q$ has a $\g$-closed segment $\R$ such that $\R'$ is also 
$\g$-closed, and beyond which no disagreements arise. (If $\N_0$ is the least $\ZF$ level of $\Q$ 
such that $\N\pins\N_0$, and if $\Tt$ is non-linear and via $\Sigma^\N$, then $\Tt$ is not 
making $\N_0$ generically generic, as its linear initial segment is too short. So $\R,\R'$ 
never disagree over non-linear trees.)} It can be that $\Q=\Lp_+^\Gamma(\N)$,
but in this case $\Q$ is not sound, whereas $\Q'$ is sound (recall $\hmb^{\Q'}=\hat{\N})$,
whereas $\hmb^\Q=\hmb^\N$).
\end{remark}

\begin{definition}\label{dfn:k-suitable}
Let $1\leq k\leq\om$. A countable
premouse $\N$ over $A$ is
\dfnemph{$k$-suitable} iff there is a strictly increasing sequence
$\left<\delta_i\right>_{i<k}$ such that
\begin{enumerate}[label=(\alph*)]
 \item For all $\delta\in(\rank(A),\OR(\N))$, we have $\N\sats$``$\delta$ is Woodin'' if and only 
if $\delta=\delta_i$ for some $i<k$.
 \item If $k=\om$ then $\OR(\N)=\sup_{i<\om}\delta_i$, and if $k<\om$ then 
$\OR(\N)=\sup_{i<\om}(\delta_{k-1}^{+i})^\N$.
\item\label{item:cutpoint} If $\N|\eta$ is a  strong 
cutpoint of $\N$ then
$\N|(\eta^+)^\N=\Lp_+^{\Gamma}(\N|\eta)$.
 \item\label{item:Qstructure} Let $\xi\in(\rank(A),\OR(\N))$ be such that $\N\sats$``$\xi$ 
is not Woodin''. Then $C_\Gamma(\N|\xi)\sats$``$\xi$ is not Woodin''.
 \end{enumerate}
We write $\delta^\N_i=\delta_i$ and $\delta_{-1}^\N=0$.
\end{definition}

Let $\N$ be $k$-suitable over $A$ and let $\xi\in(\rank(A),\OR(\N))$ be a limit ordinal 
 such that 
$\N\sats$``$\xi$ isn't Woodin''. Let $\Q\pins\N$ be the Q-structure for 
$\xi$. Let $\alpha$ be such that $\xi=\OR(\N|\alpha)$. Suppose that $\N|\alpha\pins\Q$.
Then $\alpha=\xi$ and $\N|\xi$ is $\g$-closed. In particular, $\N|\xi$ is $\g$-whole,
so $\Lp^\Gamma_+(\N|\xi)$ translates to an initial segment of $\Lp^\Gamma(\N|\xi)$.
Assume that $\N$ is reasonably iterable. If $\xi$ is a strong cutpoint of 
$\Q$, our mouse capturing hypothesis combined 
with \ref{item:Qstructure} therefore gives that $\Q\pins\Lp^{\Gamma}_+(\N|\xi)$. 
Moreover, note that if $\xi$ is a cardinal of $\N$ then $\N|\xi$ \emph{is} a 
strong cutpoint of $\Q$, since $\N$ has only finitely many Woodins. On the other hand, 
if $\xi$ is not a (strong) cutpoint of $\Q$, then one can show that $\Q\notin\Lp_+^\Gamma(\N|\xi)$,
but $\Q$ is coded over $\Lp_+^\Gamma(\N|\xi)$ (here $\Lp_+^\Gamma(\N|\xi)$ translates to a proper 
segment of $\Lp^\Gamma(\N|\xi)$).\footnote{Suppose $\xi$ is not a cutpoint of $\Q$.
Then by definition $\Q\npins\Lp^{\Gamma}_+(\N|\xi)$. Let $E\in\es_+^\Q$ be least overlapping $\xi$ 
and $\kappa=\crit(E)$. Since $\kappa$ is a limit of Woodins in $\Q$, $\kappa$ is not a cardinal of 
$\N$. Let $\P\pins\N$ be least such that $\Q\ins\P$ and $\rho_\om^\P\leq\kappa$,
and let $n<\om$ be such that $\rho_{n+1}^\P\leq\kappa<\rho_n^\P$. We claim that 
$\Lp_+^\Gamma(\N|\xi)=U$ where $U=\Ult_n(\P,E)$ (and note that $U\sats$``$\xi$ is 
Woodin'', but $\Q$ is computable from $U$, as $\Q\ins\P$ and $\P=\core_{n+1}(U)$). For $\xi$ is a 
strong cutpoint of $U$, and $U$ is $\xi$-sound but not fully sound. So it suffices to see that 
there is an above-$\kappa$, $(n,\om_1)$-iteration strategy for $\P$ in $\M|\alpha_0$. Let 
$\R\pins\N$ be least such that $\Q\ins\R$ and $\rho_\om^\R<\kappa$ (so $\P\ins\R$). Note that 
$\kappa$ is a limit of strong cutpoints of $\R$
and of Woodins of $\R$.
Let $\gamma\in(\rho_\om^\R,\kappa)$ be a strong cutpoint of $\R$,
and let $\eta$ be the least Woodin of $\R$ above $\gamma$.
Then $\eta$ is a strong cutpoint of $\R$. Since $C_\Gamma(\N|\eta)\sats$``$\eta$ is not Woodin'',
and by our mouse capturing hypothesis, therefore $\R\pins\Lp_+^\Gamma(\N|\eta)$.
In particular, there is an above-$\eta$ iteration strategy for $\R$ in $\M|\alpha_0$,
which yields the desired strategy.}

\begin{definition}[$\Gamma$-guided] Let $\P$ be
$k$-suitable and $\Tt\in\HC$ be a
normal iteration tree on $\P$. We say $\Tt$ is
\dfnemph{Q-guided} iff for each limit $\lambda<\lh(\Tt)$,
$\Q=\Q(\Tt\rest\lambda,[0,\lambda]_\Tt)$ exists and 
$\Phi(\Tt\rest\lambda)\conc(\Q,\delta(\Tt))$ is
$(\om,\om_1+1)$-iterable. We say that $\Tt$ is \dfnemph{$\Gamma$-guided} iff it is 
Q-guided, as witnessed by iteration strategies in $\M|\alpha_0$.\end{definition}

\begin{remark}\label{rem:unique_Q_guided}
Let $\P$ be $k$-suitable. For a normal tree $\Tt$ on $\P$ of limit 
length there is at most one $\Tt$-cofinal branch $b$ such that $\Tt\conc b$ is $\Q$-guided. (Let 
$b_0,b_1$ be distinct such branches; we can successfully compare the phalanxes $\Phi(\Tt\conc b_0)$ 
and $\Phi(\Tt\conc b_1)$. Standard fine structure and the fact that $\P$ has at most $\om$-many 
Woodins then leads to contradiction.) Therefore if $\Tt\conc b$ is normal, via an 
$(\om,\om_1+1)$-iteration 
strategy for $\P$, is based on $[\delta_{i-1}^\P,\delta_i^\P)$ and $\Q(\Tt,b)$ exists, then 
$\Tt\conc b$ is $\Q$-guided.
\end{remark}

\begin{definition}
\label{shortMaximalTrees}
Let $\N$ be a $\g$-whole premouse. We write $\Q^{\Gamma}_\tame(\N)$ 
for the unique $\Q\ins\Lp_+^{\Gamma}(\N)$ such that $\Q$ is a Q-structure for $\N$, if such 
exists.\footnote{The ``$\tame$'' is for \emph{tame}. While $\Q$ might not be tame, $\OR(\N)$ is 
a strong cutpoint of $\Q$.}

Let $k\leq \omega$, $\P$ be $k$-suitable and $\Tt$ a normal, limit length, $\Gamma$-guided tree on 
$\P$. We say that $\T$ is \dfnemph{short} iff $\Q^\Gamma_\tame(M(\Tt))$ exists; otherwise 
that $\T$ is
\dfnemph{maximal}.
\end{definition}

\begin{definition}\label{dfn:ss}
Let $\P$ be $k$-suitable. Let $\Tt$ be an iteration tree on $\P$. We say that $\Tt$ is 
\dfnemph{suitability strict} iff for every $\alpha<\lh(\Tt)$:
\begin{enumerate}[label=(\arabic*)]
 \item If $[0,\alpha]_\Tt$ does not drop then $M^\Tt_\alpha$ is $k$-suitable.
 \item\label{item:ss_drops} If $[0,\alpha]_\Tt$ drops and there are trees 
$\Uu,\Vv$ such that $\Tt\rest\alpha+1=\Uu\conc\Vv$, where 
$\Uu$ has last model $\R$, $b^\Uu$ does not drop, and there is $i\in[0,k)$ such that $\Vv$ is based 
on 
$[\delta_{i-1}^{\R},(\delta_i^{+\om})^{\R})$, then no $\Q\unlhd M^\Tt_\alpha$ is $(i+1)$-suitable.
 \end{enumerate}
 
 Let $\Sigma$ be a (partial) iteration strategy for $\P$. We say that $\Sigma$ is 
\dfnemph{suitability strict} iff every tree $\Tt$ via $\Sigma$ is suitability strict.
\end{definition}

\begin{definition}\label{shortTreeIterable}
Let $\P$ be $k$-suitable. We say that $\P$ is \dfnemph{short tree
iterable} iff for every normal $\Gamma$-guided tree $\Tt$ on $\P$, we have:
\begin{enumerate}[label=(\arabic*)]
 \item $\Tt$ is suitability strict.
 \item If $\T$ has limit length and is short then there is $b$ such that $\Tt\conc b$ is a 
$\Gamma$-guided tree.
  \item If $\T$ has successor length then every one-step putative normal extension 
of $\Tt$ is an iteration tree.
\end{enumerate}

Let $\P$ be short tree iterable. The \dfnemph{short tree strategy} $\Psishort_\P$ for $\P$ is the 
partial 
iteration strategy $\Psi$ for $\P$, such that $\Psi(\Tt)=b$ iff $\Tt$ is 
normal and short and $\Tt\conc b$ is $\Gamma$-guided. (By \ref{rem:unique_Q_guided} this specifies 
$\Psishort_\P$ uniquely.)
\end{definition}

\begin{lemma}\label{lem:short_tree_strat_tracks}Let $\N$ be $k$-suitable.
\begin{enumerate}[label=\tu{(}\arabic*\tu{)}]
\item\label{item:Psi_N_def} The function $\N\mapsto\Psishort_\N$, where $\N$ is short-tree 
iterable, is in $\M$; in fact, $\Psishort_\N$ is $\Gamma(\{\N,z_0\})$-definable, uniformly in 
$\N$.\footnote{But it seems that we might have 
$\Psishort_\N\notin\M|\alpha_0$.}
\item\label{item:full_strat_short_strat} Suppose there is a suitability strict normal 
$(\om,\om_1+1)$-strategy $\Sigma$ for $\N$. Then $\N$ is short tree iterable and 
$\Psishort_\N\sub\Sigma$. Moreover, for any $\Tt$ via $\Sigma$, $\Tt$ is via $\Psishort_\N$ iff for 
every 
limit $\lambda<\lh(\Tt)$, $\Q(\Tt,b)$ exists where $b=[0,\lambda)_\Tt$.
\end{enumerate}
\end{lemma}
\begin{proof}
Part \ref{item:Psi_N_def} follows from the admissibility of $\M|\alpha_0$.

Consider \ref{item:full_strat_short_strat}. Let $\Tt$ on $\N$ be normal, of limit length, via 
both $\Sigma$ and $\Psishort_\N$. Let $b=\Sigma(\Tt)$. It suffices to show that (a) if $\Q(\Tt,b)$ 
exists then $\Tt$ is short, and (b) if $\Tt$ is short then $b=\Psishort_\N(\Tt)$. (Note that if 
$\Q(\Tt,b)$ does not exist then $M^\Tt_b$ is $k$-suitable so $\Tt$ is maximal.)

Consider (a); suppose $\Q=\Q(\Tt,b)$ exists. If $b$ does not drop then $M^\Tt_b$ is suitable 
and $\delta\neq\delta_i(M^\Tt_b)$ for any $i<k$. So $C_\Gamma(M(\Tt))\sats$``$\delta$ is not 
Woodin'', so our mouse capturing hypothesis implies that $\Tt$ is short. So suppose that $b$ drops.
We can't have $C_\Gamma(M(\Tt))\sub\Q$, by suitability strictness. If $\delta$ is a cutpoint of 
$\Q$ (and so a strong cutpoint) we can then compare $\Q$ with $\Lp^\Gamma_+(M(\Tt))$; since the 
comparison is above $\delta$, we get that $\Q\ins\Lp^\Gamma_+(M(\Tt))$, so $\Tt$ is short. So 
suppose 
$\delta$ is not a cutpoint of $\Q$. Let
$E\in\es_+(\Q)$ be least such that $\kappa=\crit(E)<\delta$ and let $\Tt'$ be the normal tree 
given by $\Tt\conc\left<b,E\right>$. Then $M_\infty^{\Tt'}\sats$``$\kappa$ is a limit of Woodins'', 
so 
$b^{\Tt'}$ drops and $C_\Gamma(M(\Tt))\not\sub M_\infty^{\Tt'}$ (by suitability strictness). Also 
$M_\infty^{\Tt'}\sats$``$\delta$ is Woodin'' and $\delta$ is a cutpoint of $M_\infty^{\Tt'}$. So 
$M_\infty^{\Tt'}=\Q^\Gamma_\tame(M(\Tt))$ exists, so $\Tt$ is short.

Consider (b). Since $\Tt$ is short, $\Q=\Q(\Tt,b)$ exists. We claim that 
$\Tt\conc b$ is $\Gamma$-guided, which suffices. For it's easy to reduce to the case that 
$\delta$ is not a cutpoint of $\Q$. Let $\Tt'$ be as above,
let $\lambda=\lh(\Tt)$ and $\alpha=\pred^{\Tt'}(\lambda+1)$. 
Let $M^{*\Tt'}_{\lambda+1}=M^\Tt_\alpha|\gamma$. Then
$M^\Tt_\alpha|\gamma\sats$``$\kappa$ is a limit of cutpoints''. It follows 
that $\Tt\rest[\alpha,\lh(\Tt))$ can be considered an above-$\kappa$, normal tree on 
$M^\Tt_\alpha|\gamma$, and the iterability of the phalanx $\Phi(\Tt)\conc(\Q,\delta)$ reduces to 
the above-$\kappa$ iterability of $M^\Tt_\alpha|\gamma$, which reduces to the 
above-$\delta$ iterability of $M_\infty^{\Tt'}$ (because of the existence of 
$i^{\Tt'}_{\alpha,\lambda+1}$). But $M_\infty^{\Tt'}\ins\Lp_+^\Gamma(M(\Tt))$, so we are done.
\end{proof}

\begin{definition}
Let $A\in\pow(\RR)\inter\M$. We define the phrase $\Tt$ \dfnemph{respects $A$} as in 
\cite{Scalesweakgap}, except that we also require that $\Tt$ be suitability strict (and making any 
obvious adaptations to our setting). We define $\N$ \dfnemph{is normally $A$-iterable} as in 
\cite{Scalesweakgap}, except that we also require that $\N$ be short tree iterable. Using these 
definitions, we then define \dfnemph{(almost, locally) $A$-iterable} as in \cite{Scalesweakgap}.
\end{definition}
\begin{lemma}
\label{Comparisons}
The analogue of \tu{\cite[Lemma 1.9.1]{Scalesweakgap}} holds.
\end{lemma}
\begin{proof}This is mostly an immediate generalization. The proof in \cite{Scalesweakgap} can 
be run inside $\J(\M)$ (in fact, inside $\M$, since $\M\sats\DC_\RR$). Use suitability 
strictness to see that, for example, in the comparison of $\R|0$ with $\N|0$ (notation as in 
\cite{Scalesweakgap}), 
no tree drops on its main branch.\end{proof}

\begin{remark}
We make a further observation on the comparison above. Let $(\Tt,\Uu)$ be the $\Gamma$-guided 
portion of the comparison of, for example, $(\R|0,\N|0)$. Let $\lambda<\lh(\Tt,\Uu)$ be a limit; 
suppose $\Tt\rest\lambda$ is cofinally non-padded. So $\Q=\Q(\Tt\rest\lambda,[0,\lambda]_\Tt)$ 
exists. Then in fact, $\delta(\Tt\rest\lambda)$ is a strong cutpoint of $\Q$. For otherwise, by
the proof of \ref{lem:short_tree_strat_tracks}, $[0,\lambda]_\Tt$ drops in a manner 
which cannot be undone; i.e., for all $\alpha\geq\lambda$, $[0,\alpha]_\Tt$ drops, a contradiction. 
Similar remarks pertain to genericity iterations on $k$-suitable 
models.
\end{remark}

\begin{lemma}\label{lem:exists_A_iterable_mouse}
Let $A\in\M\inter\pow(\RR)$. Then for a cone of $s\in\RR$ 
there is an $\om$-suitable, $A$-iterable premouse over $(\MFsharp,s)$.
\end{lemma}
\begin{proof}
The following account is based on the sketch given in \cite[1.12.1]{Scalesweakgap}.\footnote{We are 
using $\g$-organized mice as our mice over reals. The authors believe that, had we used a 
hierarchy $Z$ of mice over reals more closely related to $\Theta$-$\g$-organized 
mice, then the proof in \cite[\S 7(?)]{hod_as_core_model} could be adapted to work in the 
present context. (One needs to define $Z$ such that $\Theta$-$\g$-organized mice can be realized as 
derived models of $Z$-mice, in a reasonably level-by-level manner.) Such a proof would have the 
advantage of providing some extra information. However, one would 
need to define and use the relevant Prikry forcing, so it seems to be more work overall, and our 
approach also has the advantage that it is less dependent on the precise hierarchy of mice over 
reals that is used. One might alternatively start out like \cite[\S 
7(?)]{hod_as_core_model}, but instead of using Prikry forcing, finish more like in our present 
proof.}
We give full detail here, since the 
proof is rather involved and 
the possibility of non-tame mice was not covered explicitly in \cite{Scalesweakgap}, and moreover,
comparing our proof with the remarks in \cite[Footnote 12]{Scalesweakgap}, we will not manage to 
establish the full Dodd-Jensen property for the 
iteration strategy we construct, but we will verify a version of said property 
which suffices for our purposes.

Say that a set of reals constituting a counterexample to the theorem is
\dfnemph{$\Gamma$-bad}. Suppose there is a $\Gamma$-bad set. For other pointclasses $\Gammabar$ we 
define \dfnemph{$\Gammabar$-bad} analogously.

Let $\zeta_0<\alpha_0$ and $\psi_{\Omega}\in\Sigma_1(\Ll^-)$ be such 
that $\Omega^*$ is definable over $\M|\zeta_0$ from $z_0$ and 
$\M|(\zeta_0+1)\sats\psi_\Omega(z_0)$ but $\M|\zeta_0\sats\neg\psi_\Omega(z_0)$.
Recall there is
$\xi+1\in(\theta,\l(\M))$ such that $\M|\xi\sats\ZF$.
So by \ref{Sigma1Elem} there are $\alphabar,\xibar,\betabar,\Gammabar,A$ such that:
 \begin{enumerate}[label=--]
\item $\zeta_0<\alphabar<\xibar<\betabar<\alpha_0$,
\item $ \M|\alphabar\elem_1^\RR\M|\betabar$ but $\M|\alpha'\not\elem_1^\RR\M|\alphabar$ for all
$\alpha'<\alphabar$,
\item $\Theta^{\M|\betabar}<\xibar$,
\item $\Gammabar=\rSigma_1^{\M|\alphabar}$ and $A\in\pow(\RR)^{\M|\xibar}$ and 
$\M|\xibar\sats\ZF+$``$A$ is $\Gammabar$-bad''.
\end{enumerate}
As $\M|\xibar\sats\ZF$,
$A$ really is $\Gammabar$-bad. We may assume that $\betabar$ is least such that there are 
$\alphabar,\xibar$ as 
above (relative to the fixed $\zeta_0$). Then $\betabar=\xibar+1$, $\rho_1^{\M|\betabar}=\RR$, 
$p_1^{\M|\betabar}=\{\xibar\}$ and $[\alphabar,\betabar]$ is a weak 
gap of $\M$ (the type $\rSigma_{1,(\{\xibar\},z_0)}^{\M|\betabar}$ does not reflect, using the 
choice of 
$\zeta_0,z_0$). We will show that $A$ is 
\emph{not} $\Gammabar$-bad, a contradiction.

Let $\left<A_i\right>_{i<\om}$ be a self-justifying system at the end of the gap 
$[\alphabar,\betabar]$,
with $A_0=A$. By $\AD$, in $\M|\xibar$ there is a cone of reals $s$ such that 
there is no $\om$-suitable, $A$-iterable premouse over $(\MFsharp,s)$. Let 
$z_1\geq_T z_0$ be a base for this cone such that for every $i<\om$ there is 
$\zeta<\Theta^{\M|\betabar}$ such that $A_i$ is definable over $\M|\zeta$ from $z_1$,
and a scale on $\Th_{\rPi_1}^{\M|\alphabar}(\RR)$ is definable over $\M|\betabar$ from $z_1$. We write 
$\LpGbgF$ for $\Lp^{\Gammabar}$. Recalling that $z_0$ codes $\MFsharp$, it follows that
\[ C_\Gammabar(\MFsharp,z_1)=C_\Gammabar(z_1)\psub C_\Gamma(z_1)=C_\Gamma(\MFsharp,z_1). \]
So
$\LpGbgF(\MFsharp,z_1)\pins\Lp^\Gamma(\MFsharp,z_1)$ and
both are super-small, by our mouse capturing hypothesis.
Let $\P\pins\Lp^\Gamma(\MFsharp,z_1)$ be least such that $\rho_\om^{\P}=\om$ and 
$\P\npins\LpGbgF(\MFsharp,z_1)$.
Let $\Sigma_{\P}$ be the $(\om,\om_1+1)$-strategy for $\P$. So 
$\Sigma_{\P}\in(\M|\alpha_0)\cut(\M|\alphabar)$. Let $z_2\in\RR$ code $\P$, with $z_2\geq_T 
z_1$.

We say that a pointclass $\Lambda$ is \dfnemph{lovely} iff 
$\Lambda=\rSigma_1^{\N}(z_2)\inter\pow(\RR)$ for some 
passive $\N\pins\M|\alpha_0$. Let $\left<\Gamma_i\right>_{i\in[0,9]}$ be lovely 
pointclasses such that 
$\Gammabar\sub\Delta_{\Gamma_9}$ and $(\Sigma_{\P}\rest\HC)^\HCc$ is $\Delta_{\Gamma_9}$ and for 
each 
$i\in[1,9]$, 
$\Gamma_i\sub\Delta_{\Gamma_{i-1}}$. 
Working in $\M|\xi$, let $T_0$ be the tree 
of a scale for a universal $\Gamma_0$ set. By Woodin \cite{Woodin_LCFD} applied in 
$\M|\xi$ (where $\ZF+\AD$ holds) there is 
$z_3\in\RR$ such that $z_2\leq_T z_3$ and
\[ H^*\eqdef\hod^{L_\xi[T_0,z_3]}_{T_0,z_2}\sats\text{``}\Delta_0\text{ is 
Woodin''},\]
 where $\Delta_0=\om_2^{L_\xi[T_0,z_3]}$.

Let $T_i,U_i\in H^*$ be trees projecting respectively to a universal $\Gamma_i$ set and its 
complement. Let $\Delta_i$ be least such that $V_{\Delta_i}^{H^*}$ is $\Gamma_i$-Woodin. Let 
$\lambda<\xi$ be large and such that $(V_\lambda^{H^*},\Delta_9)$ is a coarse premouse.
Let
\[ \pi_H:(H,\Delta)\to(V_\lambda^{H^*},\Delta_9)\]
be elementary, with $H\in\HC^{H^*}$, $\pi_H\in H^*$, and $z_2,T_i,U_i\in\rg(\pi)$ for each $i\leq 
9$ (let 
$U_0=\emptyset$). Let $\pi_H(T^H_i,U^H_i)=(T_i,U_i)$.
Then by arguments in \cite{twms} (using $\M|\xi$ as a background $\ZF+\AD$ model):

\begin{fact}\label{fact:strat_H} In $\M|\alpha_0$ there is a unique $(\om_1,\om_1+1)$-iteration 
strategy 
$\Lambda_H$ for $(H,\Delta)$ such that for each 
countable successor length tree $\Tt$ via $\Lambda_H$, letting $j=i^\Tt$ and $J=M_\infty^\Tt$, then
\[ p[j(T^H_8)]\sub p[T_8]\ \&\ p[j(U^H_8)]\sub p[U_8].\]
Moreover, the restriction of $\Lambda_H$ to $HC^{H^*}$ is the unique $\pi_H$-realization strategy 
in $H^*$. Further, for $i\geq 1$, $J\sats$``$j(T^H_i),j(U^H_i)$ are 
$\Coll(\om,j(\Delta))$-absolutely complementing''. Moreover,
\[ C^H\eqdef C_\Gammabar\rest V_\Delta^H\in H\ \&\ j(C^H)=C_\Gammabar\rest V_{j(\Delta)}^J; \]
\[ \Omega^H\eqdef\Omega^*\rest V_\Delta^H\in H\ \&\ j(\Omega^H)=\Omega^*\rest V^J_{j(\Delta)}.\]
\end{fact}

Let $\CC=\left<N_\alpha\right>_{\alpha\leq\Delta}$ be the maximal
$L^{\gOmega^H}[\es,(\MFsharp,z_1)]$-construction as computed in $H$ (see \ref{dfn:Fop_con}).
For every $\alpha\leq\Delta$ and 
$n<\om$, the $(n,\om_1,\om_1+1)$-strategy for $\core_n(N_\alpha)$ given by resurrection and lifting 
to $\Lambda_H$, is a $\gOmega^*$-strategy;
this is by and \ref{fact:strat_H}, \ref{lem:gF_props} and properties of the resurrection/lifting 
maps. So by \ref{lem:fine_structure}, this construction does indeed have length $\Delta+1$.

\begin{claim}
There is $\gamma<\Delta$ and $k<\om$ such that 
$\rho_{k+1}^{N_\gamma}=\om$ and $\core_\om(N_{\gamma})$ is not 
$(k,\om_1+1)$-iterable in $\M|\alphabar$.\end{claim}
\begin{proof}
It suffices to see that $\CC$ reaches 
$\P$.
We have $z_2,\P\in\HC^H$, and by the definability of $\Sigma_{\P}\rest\HC$, letting 
$\Sigma_{\P}^H=\Sigma_{\P}\rest V_\Delta^H$, we have $\Sigma_{\P}^H\in H$, and 
$\Sigma_{\P}^H$ is moved correctly by $\Lambda_H\rest\HC$. It follows that the background 
extenders used in 
$\CC$ all cohere $\Sigma_{\P}^H$, and so we can apply \ref{lem:back_stat}
(the stationarity of $\CC$ with respect to $\P$). So we just need to rule out the possibility 
that 
for some normal tree $\Tt$ on $\P$ via $\Sigma_{\P}$, with last model $\P'$, $N_\Delta\ins\P'$. 
But 
because $(\Sigma_{\P}\rest\HC)^\HCc$ and $(\Omega^*)^\HCc$ are $\Delta_{\Gamma_9}$ and $N_\Delta$ 
is definable over $(V_\Delta^H,\Omega^H)$, 
we have $\Tt\in C_{\Gamma_9}(V_\Delta^H)$. But $C_{\Gamma_9}(V_\Delta^H)\sats$``$\Delta$ is 
Woodin'', so by the universality of $N_\Delta$ (see \cite[Lemma 11.1]{DMATM}), $\Tt\notin 
C_{\Gamma_9}(V_\delta^H)$, contradiction.
\end{proof}

We will now look at the least stage where the construction produces a
fine structurally nice mouse which is not iterable in $\M|\alphabar$.
This move, and its relation to producing a mouse with $\om$ Woodins and a suitability strict 
iteration strategy, is related to, and motivated 
by, an argument shown to the first author by 
Steel, in a similar situation, though a different context.

Given a $k$-sound premouse $\N\in\HC$ and $\zeta\in\OR(\N)$, we say that $\N$ is 
\dfnemph{$(\Gammabar,k,\zeta)$-iterable} 
iff there is an above-$\zeta$, $(k,\om_1+1)$-iteration strategy for $\N$ in $\M|\alphabar$.
We say that $\N$ is $(\Gammabar,\zeta)$-iterable iff $\N$ is $(\Gammabar,m,\zeta)$-iterable,
where $m$ is defined in the next paragraph.

By the previous claim, we may let $(\gamma,m,\eta')\in\Ord^3$ be lex-least 
such that, letting $\Ss=\core_m(N_\gamma)$, $\Ss|\eta'$ is a $\g$-whole
cutpoint of 
$\Ss$ and
\[ \R'\eqdef\cHull_{m+1}^{\Ss}(\eta'\un p_{m+1}^{\Ss}) \]
is $\eta'$-sound and not $(\Gammabar,m,\eta')$-iterable. Let $\pi':\R'\to\Ss$ 
be the uncollapse. (It follows that $\pi'(p_{m+1}^{\R'}\cut\eta')=p_{m+1}^\Ss\cut\eta'$. We allow 
$\eta'<\rho_{m+1}^\Ss$, so we do need to assume $\eta'$-soundness explicitly.)
It seems that $\eta'$ could be measurable in $\R'$, which is slightly inconvenient.
So we first replace $\R'$ with a slightly larger hull $\R$,
and replace $\eta'$ with a strong cutpoint $\eta$ of $\R$.

Given a premouse $\N$ and $\eta<\OR(\N)$, we say that $\eta$ is \dfnemph{$\N$-finely measurable}
iff $\eta=\crit(E)$ for some $\N$-total measure $E$ such that either $E\in\es_+^\N$,
or $E\in\es_+^{\Ult(\N,F)}$ for some $F\in\es_+^\N$.

We claim that $\eta'<\min(\rho_m^{\R'},\rho_m^\Ss)$ and $\eta'$ is not measurable in $H$,
nor $\Ss$-finely measurable. For $\rho_m^{\R'}$ is 
the least $\rho$ such that either $\rho\notin\dom(\pi')$ or $\pi'(\rho)\geq\rho_m^{\Ss}$, by 
elementarity. We have $\eta'<\rho_m^{\R'}$
(as otherwise $\R'$ is not $(\Gammabar,m-1,\eta')$-iterable,
which implies that $\core_{m-1}(N_\gamma)$ is not $(\Gammabar,m-1,\rho_m^{N_\gamma})$-iterable,
contradicting the minimality of $m$),
so also $\eta'<\rho_m^\Ss$. Since $\eta'<\rho_m^\Ss$, if $\eta'$ is $\Ss$-finely measurable then 
$\eta'$ is measurable in $H$. But if $H\sats$``$\mu$ is a normal measure on $\eta'$'' and 
$j:H\to\Ult(H,\mu)$
is the ultrapower map, then
\[ \R'=\cHull_{m+1}^{j(\Ss)}(\eta'\un p_{m+1}^{j(\Ss)}), \]
which contradicts the minimality of $j(\eta')$ in $\Ult(H,\mu)$. (The minimality can be computed 
correctly in $H$ and its $\Lambda_H$-iterates by \ref{fact:strat_H}.)

Now let $\eta=((\eta')^+)^\Ss$. We claim that $\eta<\rho_m^\Ss$ and $\Ss|\eta$ is a 
$\g$-whole strong cutpoint of $\Ss$ and
\[ \R\eqdef\cHull_{m+1}^\Ss(\eta\un p_{m+1}^\Ss) \]
is $\eta$-sound and not $(\Gammabar,m,\eta)$-iterable.
(Therefore $\eta=((\eta')^+)^\R$ is also a strong cutpoint of $\R$.) For 
suppose $\eta=\rho_m^\Ss$.
Then $\pi'(\eta')=\eta'$ because otherwise we contradict the minimality of $m$, as above.
So $\rho_m^{\R'}=((\eta')^+)^{\R'}$ and $\eta'$ is not $\R'$-finely measurable.
But then any above-$\eta'$ tree on $\R'$ immediately drops either in model or to degree $\leq m-1$,
which contradicts the minimality of $(\gamma,m)$. In particular, $\eta=((\eta')^+)^\Ss<\OR(\Ss)$,
so 
$\Ss|\eta$ is $\g$-whole,
and since $\eta'$ is not $\Ss$-finely measurable, $\eta$ is a strong cutpoint of $\Ss$.
Clearly $\R$ is $\eta$-sound. If $\pi'(\eta')>\eta'$
then $\eta\leq\pi'(\eta')<\rho_m^\Ss$, which easily gives that $\R$ is not 
$(\Gammabar,m,\eta)$-iterable. If $\pi'(\eta')=\eta'$ then $\eta'$ is not $\R'$-finely measurable,
which implies that $\R'$ is not $(\Gammabar,m,((\eta')^+)^{\R'})$-iterable,
so $\R$ is not $(\Gammabar,m,\eta)$-iterable.

Let $\pi_0:\R\to\Ss$ be the uncollapse embedding.
Let 
$\Sigma_\R$ be the above-$\eta$,
$(m,\om_1,\om_1+1)$-strategy for $\R$ given by resurrection and lifting to $\Lambda_H$, taking 
$\pi_0$ as the base lifting map.
Let $\Tt$ be on $\R$ via $\Sigma_\R$ and $\lambda<\lh(\Tt)$, and 
let $\Uu$ be the lifted tree on $H$. Write $\CC_\lambda=i^\Uu_{0,\lambda}(\CC)$. Let 
$n=\deg^\Tt(\lambda)$. Let $(\gamma^\Tt_\lambda,\Ss^\Tt_\lambda,\pi^\Tt_\lambda)$
be the $(\gamma',\Ss',\pi')$ produced by lifting/resurrection such that
$\gamma'\leq i^\Uu_{0,\lambda}(\gamma)$ and $\Ss'=\core_n(N^{\CC_\lambda}_{\gamma'})$
and $\pi':M^\Tt_\lambda\to\Ss'$ is the lifting map.
(In particular, $\pi^\Tt_\lambda$ is a weak $n$-embedding,
and $\gamma^\Tt_\lambda=i^\Uu_{0,\lambda}(\gamma)$ iff $[0,\lambda]_\Tt$ does not drop in model.
Here if $[0,\lambda]_\Tt$ does not drop in model, the codomain of $\pi_\lambda$ is 
$i^\Uu_{0,\lambda}(\Ss)$, not $i^\Uu_{0,\lambda}(\R)$.)

Let $\Tt$ be an above-$\eta$ normal tree on $\R$, of countable limit length. Let $b$ be a 
$\Tt$-cofinal 
branch. Let $\Q_b=\Q(\Tt,b)$.
Then $k(\Tt,b)$ denotes $\om$ if $\Q_b\pins M^\Tt_b$, and denotes $\deg^\Tt(\lambda)$ otherwise.
And $\Phi_\Q(\Tt,b)$ denotes the phalanx $\Phi(\Tt)\conc(\Q_b,k)$, where $k=k(\Tt,b)$.
(In the phalanx notation, $k$ denotes the base degree corresponding to $\Q_b$.)
We say that $b$ is \dfnemph{$\Gammabar$-verified for $\Tt$} iff 
$\Phi_\Q(\Tt,b)$ is normally $(\om_1+1)$-iterable in 
$\M|\alphabar$.

\begin{claim}\label{clm:if_phalanx_not_iter}
Let $\Tt$ be normal on $\R$ via $\Sigma_\R$, of length $\lambda+1$ for some limit 
$\lambda<\om_1$.
Suppose that
$\pP\eqdef\Phi_\Q(\Tt\rest\lambda,b)$ is not normally $(\om_1+1)$-iterable in $\M|\alphabar$. 
Let $M_\lambda=M^\Tt_\lambda$, $b=b^\Tt$, $\Q=\Q(\Tt\rest\lambda,b)$,
$k=k(\Tt\rest\lambda,b)$,
$\delta=\delta(\Tt\rest\lambda)$ and $M_\Tt=M(\Tt\rest\lambda)$.
Then either:
\begin{enumerate}[label=\tu{(}\roman*\tu{)}]
 \item\label{item:delta_strong_cut_Q} $\delta$ is a strong cutpoint of $\Q=M_\lambda$, $b$ does not 
drop in model or 
degree and $\Q||(\delta^+)^\Q=\LpGbgF_+(M_\Tt)$;
or
 \item\label{item:delta_not_cut_Q} $\delta$ is not a cutpoint of $\Q$, 
and letting $E\in\es_+^\Q$ be such that $\crit(E)<\delta<\lh(E)$, with 
$\lh(E)$ minimal, and letting $\Tt^+$ be the normal tree 
$\Tt\conc\left<E\right>$, then $b^{\Tt^+}$ does not drop in model or degree, 
and $\Q||\lh(E)=\LpGbgF_+(M_\Tt)$.
\end{enumerate}
\end{claim}
\begin{proof}
Let $(\gamma_\lambda,\Ss_\lambda,\pi_\lambda)=(\gamma^\Tt_\lambda,\Ss^\Tt_\lambda,\pi^\Tt_\lambda)$.
Suppose $\delta$ is a cutpoint (hence strong cutpoint) of $\Q$. Because $\delta$ is a 
cutpoint, the difficulty 
in iterating $\pP$ gives that $\Q$ is not $(\Gammabar,k,\delta)$-iterable.
Because $\delta$ is a strong cutpoint and by standard fine structure, 
$\Q\ins\Lp^\Gamma_+(M_\Tt)$.

We leave the proof that $\Q=M_\lambda$ to the reader; assume this. We show that $b$ does not 
drop in model or degree; suppose otherwise. We have 
\[ \Q=\Hull_{k+1}^\Q(\delta\un p_{k+1}^\Q).\]
and $(\gamma_\lambda,k)<_\lex 
(i^\Uu_{0,\lambda}(\gamma),m)$ and
$p_{k+1}^{\Ss_\lambda}=\pi_{\lambda}(p_{k+1}^\Q)$ and (by
the commutativity between the copy and iteration maps after the last drop) and
\[ \rg(\pi_\lambda)\sub\R^*\eqdef\Hull_{k+1}^{\Ss_\lambda}(\pi_\lambda(\delta)\un 
p_{k+1}^{\Ss_\lambda}).\]
Let $\R'$ be the transitive collapse of $\R^*$ and let $\sigma:\Q\to\R'$ be the obvious map,
a weak $k$-embedding with
$\sigma(\delta)=\pi_\lambda(\delta)$. So $\sigma$ lifts above-$\delta$ trees on $\Q$ to 
above-$\pi_\lambda(\delta)$ trees on $\R'$. Therefore $\R'$ is not 
$(\Gammabar,k,\pi_\lambda(\delta))$-iterable. But $\R'$ is 
$\pi_\lambda(\delta)$-sound, as there are generalized $(k+1)$-solidity 
witnesses for $(\Ss_\lambda,p_{k+1}^{\Ss_\lambda})$ in $\rg(\pi_\lambda)$ (by commutativity as 
before). This contradicts the minimality of 
$(i^\Uu_{0,\lambda}(\gamma),m)$ in $M^\Uu_\lambda$.

So $b^\Tt$ does not drop. One can show $\Q||(\delta^+)^\Q\ins\LpGbgF_+(M_\Tt)$ much as above. But 
$\Q\nins\LpGbgF_+(M_\Tt)$, as $\Q$ is not $(\Gammabar,k,\delta)$-iterable. So 
$\Q||(\delta^+)^\Q=\LpGbgF_+(M_\Tt)$, as required.

Now suppose $\delta$ is not a cutpoint of $\Q$. Suppose that $b^{\Tt^+}$ drops 
in model or degree. Since $\delta$ is a strong 
cutpoint of $M_\infty^{\Tt^+}$, then as before, by choice of $(\gamma,m)$, 
$M_\infty^{\Tt^+}$ is $(\Gammabar,j,\delta)$-iterable, 
where $j=\deg^{\Tt^+}(M_\infty^{\Tt^+})$. Therefore, letting $\kappa=\crit(E)$ and 
$\xi=\lambda+1$, 
$M^{*\Tt^+}_{\xi}$ is $(\Gammabar,j,\kappa)$-iterable (we can copy trees using $i_E$). But 
$\kappa$ is a 
cutpoint of $M^{*\Tt^+}_{\xi}$. So 
$\Tt^+=(\Tt\rest\chi+1)\conc\Tt'$, where 
$\chi=\pred^\Tt(\xi)$ and $\Tt'$ is an above-$\kappa$, $j$-maximal tree on 
$M^{*\Tt^+}_\xi$. Thus, the iterability of $\pP$ can be reduced to the above-$\kappa$ iterability 
of 
$M^{*\Tt^+}_\xi$. Therefore $\pP$ is iterable in $\M|\alphabar$, 
a contradiction. So $b^{\Tt^+}$ does not drop. We then get $\Q||\lh(E)=\LpGbgF_+(M_\Tt)$ 
by the arguments just given.
\end{proof}

\begin{claim}\label{clm:partial_strategy}
Let $\Tt$ be a normal tree on $\R$, via $\Sigma_\R$, of countable limit length. Then there is 
at most one branch $\Gammabar$-verified for $\Tt$.
However, the following partial strategy
$\Psi$ is not an above-$\eta$, $(m,\om_1)$-strategy for $\R$: Given 
$\Tt$, let $\Psi(\Tt)$ be the unique branch which is $\Gammabar$-verified for $\Tt$.
\end{claim}
\begin{proof}
Uniqueness follows from the usual comparison and fine structural arguments, 
using the $\eta$-soundness of $\R$.
If existence holds then by uniqueness and because $\M|\alphabar$ 
is admissible, $\R$ is $(\Gammabar,\eta)$-iterable,
contradiction.
\end{proof}

\begin{definition}
We define the term \dfnemph{$\Gammabar$-$k$-suitable} analogously to \emph{$k$-suitable} (cf. 
\ref{dfn:k-suitable}), but 
with $\Gammabar$ replacing $\Gamma$. We likewise define \dfnemph{$\Gammabar$-$A$-iterable} and 
\dfnemph{$\Gammabar$-suitability strict}. Let $R$ be $\Gammabar$-$\om$-suitable with $z_1\in R$.
Then $\sigma_i^R$ denotes the $\Coll(\om,\delta_i^R)$-term capturing $A_i$ over $R$ (see 
\cite{twms}). Let $Q$ be a structure and $\pi:Q\to P$. We say that $\pi$ is an 
\dfnemph{$\Avec$-embedding} iff $\pi$ is $\Sigma_1$-elementary and $\sigma_i^R\in\rg(\pi)$ for all 
$i<\om$.
\end{definition}

\begin{claim}\label{clm:N_om_suitable} \tu{(}i\tu{)} $\Ss$ has infinitely many Woodins 
$\delta$ such that
$\eta<\delta<\rho_m^\Ss$. Let $\delta_\om$ be the supremum of the first $\om$-many and let 
$\N$ be the translation of $\Ss|\delta_\om$ to a $\g$-organized 
spm over $\widehat{\Ss|\eta}$ \tu{(}translated as in 
\ref{rem:translation_Theta-g_above_cutpoint}\tu{)}. Then \tu{(}ii\tu{)} $\N$ is
$\Gammabar$-$\om$-suitable.\end{claim}
\begin{proof}
We will construct a $\Gammabar$-$\om$-suitable premouse 
which is an initial segment of a $\Sigma_\R$-iterate of $\R$. This is by 
applying Claim \ref{clm:partial_strategy} and an obvious 
generalization thereof, in tandem with Claim \ref{clm:if_phalanx_not_iter}, up to $\om$ 
many times. So let $\Tt_0$ on $\R_0=\R$ be via $\Sigma_\R$ (so above 
$\delta_{-1}\eqdef\eta$), witnessing the failure of ``existence'' in 
\ref{clm:partial_strategy}, with $\Tt_0$ of minimal length. Let $\delta_0=\delta(\Tt_0)$. Let 
$b=\Sigma_\R(\Tt_0)$. So \ref{clm:if_phalanx_not_iter} applies to 
$\Phi_\Q(\Tt_0,b)$. Use notation as there, so $\Tt=\Tt_0\conc b$ 
and $\delta=\delta_0$.

Suppose first that \ref{clm:if_phalanx_not_iter}\ref{item:delta_not_cut_Q}
holds. Let 
$\kappa=\crit(E)$. Since $E$ overlaps 
$\delta$ and $b^{\Tt^+}$ does not drop in model or 
degree, $\kappa$ is a limit of Woodins of $M_\infty^{\Tt^+}$, 
and 
$\eta<\kappa<\delta<\rho_m(M_\infty^{\Tt^+})$ (recall we arranged that $\eta$ is a strong 
cutpoint of $\R$). And $M_\infty^{\Tt^+}$ 
is not $(\Gammabar,\delta)$-iterable. Now 
let 
$\delta^*_\om$ be the supremum of the first $\om$-many Woodins of $M_\infty^{\Tt^+}$ above $\eta$. 
Let 
$\zeta$ be least such that $\delta^*_\om<\lh(E^\Tt_\zeta)$. So 
$M_\infty^{\Tt^+}|\delta^*_\om=M^\Tt_\zeta|\delta^*_\om$. Note $\delta^*_\om$ is a strong 
cutpoint of $M^\Tt_\zeta$ and $\zeta\in b^{\Tt^+}$, so $[0,\zeta]_\Tt$ does not drop 
in model or degree. Therefore $M^\Tt_\zeta$ is not $(\Gammabar,\delta^*_\om)$-iterable. Now let 
$\Uu$ be the lifted tree, via $\Lambda_H$, on $H$. Then
$\eta<\pi^\Tt_\zeta(\delta^*_\om)<\rho_m(\Ss^\Tt_\zeta)$ and $\pi^\Tt_\zeta(\delta^*_\om)$ is the 
sup of 
the first 
$\om$ Woodins of 
$\Ss_\zeta$ above $\eta$, and $\Ss_\zeta$ is not 
$(\Gammabar,\pi^\Tt_\zeta(\delta^*_\om))$-iterable. By the elementarity of $i^\Uu_{0,\zeta}$, this 
gives (i).

We now verify condition \ref{item:cutpoint} of $\Gammabar$-$\om$-suitability. Let $\kappa\geq\eta$ 
be a cutpoint of $\Ss|\delta_\om$ with $\eta\leq\kappa$. Let 
$\C_\kappa$ be the $\kappa$-core of $\Ss$. We claim that ($*$) $\C_\kappa$ is not 
$(\Gammabar,\kappa)$-iterable. For
let $\xi\in b^\Tt$ be least such that $\pi_\xi(\lh(E^\Tt_\xi))>i^\Uu_{0,\xi}(\kappa)$.
Let $\kappabar$ be the least such that $\pi_\xi(\kappabar)\geq 
i^\Uu_{0,\xi}(\kappa)$. Then 
$\nu(E^\Tt_\alpha)\leq\kappabar$ for all $\alpha+1\leq_\Tt\xi$, and $\kappabar$ is a 
cutpoint of $M^\Tt_\xi$ (as $\kappa$ is a cutpoint of $\Ss$). Therefore $M^\Tt_\xi$ is not 
$(\Gammabar,\kappabar)$-iterable,
and
\[ \rg(\pi_\xi)\sub\Hull_{m+1}^{\Ss_\xi}(i^\Uu_{0,\xi}(\kappa)\un p_{m+1}^{\Ss_\xi}), \]
so $i^\Uu_{0,\xi}(\C_\kappa)$ is not $(\Gammabar,i^\Uu_{0,\xi}(\kappa))$-iterable,
giving ($*$).

Now let $\Ss|\kappa$ be a $\g$-whole strong cutpoint of $\Ss|\delta_\om$.
By the choice of $\gamma$, we have $\Ss|(\kappa^+)^\Ss\ins\LpGbgF_+(\Ss|\kappa)$.
But letting $\C_{\kappa+1}$ be the $(\kappa+1)$-core of $\Ss$, by ($*$), we have
$\LpGbgF_+(\Ss|\kappa)\pins\C_{\kappa+1}$. Condition \ref{item:cutpoint} follows.

We now verify condition \ref{item:Qstructure}. Let $\eta\leq\xi<\delta_\om$ with $\Ss\sats$``$\xi$ 
is not Woodin''; we must show that $C_\Gammabar(\Ss|\xi)\sats$``$\xi$ is not Woodin''. 
We may assume that $\Ss|\xi$ is $\g$-whole, and by \ref{item:cutpoint}, that 
$\xi$ is not a strong cutpoint of $\Ss$. Let $F\in\es^\Ss$ be least such that 
$\mu=\crit(F)\leq\xi<\lh(F)$. Note that 
$\mu$ is a limit of strong cutpoints of $\Ss|\xi$. So if $\mu=\xi$ then $\Ss|\xi$ is the 
Q-structure for $\xi$, so we are done. So suppose $\mu<\xi$. We may assume that 
$\Ss||\lh(F)\sats$``$\xi$ is Woodin'', 
because otherwise there is $\Q\pins\Ss||\lh(F)$ such that 
$\Q$ is a Q-structure for $\xi$ and $\xi$ is a strong cutpoint of 
$\Q$, and so $\Q\pins\LpGbgF_+(\Ss|\xi)$ (by choice of $\gamma$). Therefore 
$\mu$ is not a cardinal of $\Ss$. Let $\Q\pins\Ss$ be least such that $\lh(F)\leq\OR(\Q)$ 
and $\rho_\om^\Q<\mu$. Then $\Q$ collapses $\xi$. Let 
$\zeta\in[\rho_\om^\Q,\mu)$ be a $\g$-whole strong cutpoint of $\Q$. Then 
$\Q\pins\LpGbgF_+(\Ss|\zeta)$, so $\Q\in C_\Gammabar(\Ss|\zeta)$, which suffices. This completes 
the proof that $\Ss|\delta_\om$ is $\Gammabar$-$\om$-suitable in this case.

Now suppose that conclusion (a) of Claim \ref{clm:if_phalanx_not_iter} holds.
Let $\Tt_0^+=\Tt_0\conc b$ and let $\R_1=M_\infty^{\Tt_0^+}$. 
Then $b^{\Tt_0^+}$ does not drop in model or degree. And $\delta_0$ is a strong cutpoint of $\R_1$, 
$\R_1$ is $\delta_0$-sound, projects $<\delta_0$, and is not $(\Gammabar,\delta_0)$-iterable.
So the obvious modification of Claim 
\ref{clm:partial_strategy} applies to $\R_1$ above $\delta_0$. Pick $\Tt_1$ 
on $\R_1$, above $\delta_0$, like before. Again apply Claim 
\ref{clm:if_phalanx_not_iter}. If its conclusion (b) holds proceed as before, and 
otherwise let $\R_1=M_\infty^{\Tt_1^+}$ and pick $\Tt_2$ on $\R_1$, etc.

If the above process produces $\R_n$ and $\Tt_n$ for all $n<\om$, then we get (i) much as before,
and note that, letting $\delta_n$ be the $n^\nth$ Woodin of $\Ss$ above $\eta$, 
then $\Ss$ is not $(\Gammabar,\delta_n)$-iterable. Part (ii) follows much like before.
\end{proof}

\begin{claim}\label{clm:term_relation_hull}
Let $\P$ be $\Gammabar$-$\om$-suitable and let $\pi:\Q\to\P$ be an $\Avec$-embedding. Then 
\tu{(}i\tu{)} $\Q$ is 
$\Gammabar$-$\om$-suitable and for each $i<\om$, \tu{(}ii\tu{)} $\pi(\sigma^\Q_i)=\sigma^\P_i$, and 
\tu{(}iii\tu{)} $\rg(\pi)$ is cofinal in $\delta_i^\P$.
\end{claim}
\begin{proof}
Parts (i) and (ii) are by condensation of term relations for self-justifying-systems; see 
\cite{twms}.
Consider (iii). If 
$\rg(\pi)\inter\delta_i^\P$ is 
bounded in $\delta_i^\P$, then we may assume that $\crit(\pi)=\delta_i^\Q$, by taking the 
appropriate 
hull (cf. the first part of the proof of \cite[Lemma 1.16.2]{Scalesweakgap}). 
But then $\Q|\delta_i^\Q=\P|\delta_i^\Q$, and $\P|\delta_i^\Q$ is not $\Gammabar$-Woodin, but 
$\Q\sats$``$\delta_i^\Q$ is Woodin'', so $\Q$ is not $\Gammabar$-$\om$-suitable, contradiction.
\end{proof}

\begin{definition} Let $\Tt=\left<\Tt_\alpha\right>_{\alpha\leq\gamma}$ be a stack of normal 
iteration trees. We say that $\Tt$ is \dfnemph{relevant} iff for every $\alpha<\gamma$, 
$b^{\Tt_\alpha}$ does not drop. (Here we allow $\Tt_\gamma$ to be trivial, and it might drop.)
(Recall from \ref{dfn:F-iterability} that a \emph{hod} iteration strategy acts on relevant trees.)
\end{definition}

From now on we fix $\N$ as defined in Claim \ref{clm:N_om_suitable}. Let $\Sigma_\N$ be the 
hod-$(\om,\om_1,\om_1+1)$ strategy for $\N$ given by resurrection and 
lifting to $\Lambda_H$. The next claim follows from \ref{fact:strat_H}.
\begin{claim}\label{clm:Lambda_H_pres}For any successor length tree $\Uu$ on $H$ via 
$\Lambda_H$, 
$i^\Uu(N)$ is $\Gammabar$-$\om$-suitable and $i^\Uu\rest\N:\N\to i^\Uu(\N)$ is 
an $\Avec$-embedding.
\end{claim}

\begin{claim}\label{clm:Sigma_N_ss}
$\Sigma_\N$ is $\Gammabar$-suitability strict. Moreover, let $\Tt$ be via $\Sigma_\N$, of successor 
length, such that $b^\Tt$ does not drop. Then $i^\Tt$ is an $\Avec$-embedding.
 \end{claim}

\begin{proof}
Let $\Tt$ be via $\Sigma_\N$, of successor length. If $b^\Tt$ does not drop, then the 
lemma's conclusions regarding $M_\infty^\Tt$ and $i^\Tt$ follow from \ref{clm:term_relation_hull} 
and
\ref{clm:Lambda_H_pres}.

Suppose $b^\Tt$ drops and that $i<\om$ is as in \ref{dfn:ss}\ref{item:ss_drops}, but some 
$\R\ins M_\infty^\Tt$ is $\Gammabar$-$(i+1)$-suitable. For simplicity assume that $\Tt$ consists of 
just one normal tree and that $\Tt$ has minimal possible length. It follows that for every extender 
$E$ used in $\Tt$, 
$\nu(E)<\delta=\delta_i^\R$. Let $n=\deg^\Tt(b^\Tt)$. Then $\rho_{n+1}(M_\infty^\Tt)<\OR(\R)$ and 
$M_\infty^\Tt$ is $\delta$-sound. So let $\Q\ins M_\infty^\Tt$ be least such that $\R\ins \Q$ and 
$\rho_\om^\Q\leq\delta$. So
\[ \Q|(\delta^+)^\Q=\R|(\delta^+)^\R=\Lp^\Gammabar_+(\R|\delta),\]
$\Q\sats$``$\delta$ is 
Woodin'',
$\Q$ is $\delta$-sound and $\delta$ is a strong cutpoint of $\Q$. So letting $j<\om$ be such that 
$\rho_{j+1}^\Q\leq\delta<\rho_j^\Q$, $\Q$ is not 
$(\Gammabar,j,\delta)$-iterable. Let $\Uu$ be the 
$\Lambda_H$-tree on $H$ given by lifting $\Tt$.
Suppose for simplicity that $\Q=M^\Tt_\infty$. Because of the drop,
$\Ss^\Tt_\infty$ is $(\Gammabar,j,\pi^\Tt_\infty(\delta))$-iterable, so $\Q=M^\Tt_\infty$ is 
$(\Gammabar,j,\delta)$-iterable, contradiction. If $\Q\pins M^\Tt_\infty$ it is 
similar.\footnote{Suppose $M^\Tt_\infty$ is active type 3 and 
$\nu(E(M_\infty^\Tt))<\OR(\Q)<\OR(M_\infty^\Tt)$. Let $E^*\in M^\Uu_\infty$ be a background 
extender for 
$\Ss^\Tt_\infty$ and lift $\Q$ to a model in $\Ult(M^\Uu_\infty,E^*)$.}
\end{proof}

\begin{definition}
Let $\Q$ be $\Gammabar$-$\om$-suitable. Let $\Sigma$ be a hod-$(\om,\om_1,\om_1+1)$-strategy for 
$\Q$. We say that $(\Tt,\P)$ is a \dfnemph{$\Sigma$-pair} iff $\Tt$ is a 
countable tree on $\Q$ via $\Sigma$, with last model $\P$. We say that a $\Sigma$-pair $(\Tt,\P)$ 
is 
\dfnemph{non-dropping} iff $b^\Tt$ does not drop. We say that $\Sigma$ is
\dfnemph{$\Avec$-good} iff for every non-dropping $\Sigma$-pair $(\Tt,\P)$, $\P$ is 
$\Gammabar$-$\om$-suitable and $i^\Tt$ is an $\Avec$-embedding.
If $(\Tt,\P)$ is a non-dropping $\Sigma$-pair, we write $\Sigma^\Tt_\P$ for the $(\Tt,\P)$-tail of 
$\Sigma$ (that is, $\Sigma^\Tt_\P$ is the hod-$(\om,\om_1,\om_1+1)$ iteration strategy 
$\Lambda$ for $\P$ where $\Lambda(\Uu)=\Sigma(\Tt,\Uu)$).
\end{definition}

The following claim is immediate:

\begin{claim}\label{clm:propagate_ss_Ag}
 Let $\Sigma$ be a hod-$(\om,\om_1,\om_1+1)$-iteration strategy for $\Q$.
 Let $(\Tt,P)$ be a non-dropping $\Sigma$-pair. If $\Sigma$ is suitability strict then 
$\Sigma^\Tt_P$ is suitability strict. If $\Sigma$ is $\Avec$-good then $\Sigma^\Tt_P$ is 
$\Avec$-good.
\end{claim}

\begin{claim}\label{clm:unique_ss_Ag}
Let $\Q$ be $\Gammabar$-$\om$-suitable. Then there is at most one suitability strict 
$\Avec$-good hod-$(\om,\om_1,\om_1+1)$ iteration strategy for $\Q$.\end{claim}
\begin{proof}
Let $\Sigma,\Lambda$ be two such strategies, and let $\Tt$ be of limit length, via 
$\Sigma,\Lambda$, such that $b=\Sigma(\Tt)\neq\Lambda(\Tt)=c$. We may assume that $\Tt$ is 
normal. We can compare the phalanx $\Phi(\Tt)\conc b$ with the phalanx $\Phi(\Tt)\conc c$, forming 
trees $\Uu,\Vv$, using $\Sigma,\Lambda$, respectively. The comparison is successful. By suitability 
strictness, we have $M_\infty^\Uu=\P=M_\infty^\Vv$. By standard fine structure, $b^\Uu$ and 
$b^\Vv$ do not drop and $M_\infty^\Uu\sats$``$\delta(\Tt)$ is Woodin''. In particular, 
$\delta(\Tt)=\delta^\P_k$ for some $k<\om$. Because $\Sigma,\Lambda$ are $\Avec$-strategies and by 
\ref{clm:term_relation_hull}, therefore $\rg(i^\Uu)\inter\rg(i^\Vv)$ is unbounded in $\delta_k^\P$. 
But then $\rg(i^\Tt_b)\inter\rg(i^\Tt_c)$ is unbounded in $\delta_k^\P$, so $b=c$.
\end{proof}

We are now in a position to establish a version of 
Dodd-Jensen.

\begin{claim}\label{clm:DJ}
 Let $\Sigma$ be an $\Avec$-good, suitability strict strategy for $\Q$. Let $(\Tt,\P)$ be a 
non-dropping $\Sigma$-pair.
\begin{enumerate}[label=\tu{(}\arabic*\tu{)}]
\item\label{item:pullback} Let $\pi:\R\to\P$ be an $\Avec$-embedding. Then the $\pi$-pullback 
$\Lambda$ of $\Sigma^\Tt_\P$ is $\Avec$-good and suitability strict. Therefore if $\R=\Q$ then 
$\Lambda=\Sigma$.
\item\label{item:DJ} Let $\pi:\Q\to\P$ be an $\Avec$-embedding. Then for all $\alpha<\OR(\Q)$, 
$i^\Tt(\alpha)\leq\pi(\alpha)$.
\end{enumerate}
\end{claim}
\begin{proof}
The first clause of \ref{item:pullback} is proven like \ref{clm:Sigma_N_ss}. This together with 
\ref{clm:unique_ss_Ag} 
yields the second clause.
For \ref{item:DJ}, the standard Dodd-Jensen proof works; the copying does 
not break down by \ref{item:pullback}.
\end{proof}

One can now deduce that $\N$ is $\Gammabar$-$A$-iterable, because \ref{clm:propagate_ss_Ag} 
and \ref{clm:DJ} apply to $\N$ and $\Sigma_\N$, which is enough Dodd-Jensen for 
$\Sigma_\N$ to apply the proof of 
\cite[Theorem 4.6]{hodlr}. Recall that $\N$ is over $\widehat{\Ss|\eta}$. Let 
$g\sub\Coll(\om,\Ss|\eta)$ be 
$\N$-generic. Let $x\in\RR\inter(\N|1)[g]$ code $(\N|\eta,g)$. Then we can 
reorganize $\N[x]$ as a premouse $\N^*$ over $(\MFsharp,x)$, and $\N^*$ is 
$\Gammabar$-$\om$-suitable and $\Gammabar$-$A$-iterable; these facts all follow by 
S-construction (for $\g$-organized spms; cf. 
\ref{rem:S-construction}). But $x\geq_T z_1$, 
contradicting the 
choice of $z_1$. This completes the proof of \ref{lem:exists_A_iterable_mouse}.
\end{proof}

Now for simplicity assume $n=0$ and $\beta=\l(\M)$ is a limit ordinal; we allow that 
$\Upsilon^\M\neq\emptyset$.
Let $p,w_1,W,\Sigma,\left<\beta_i,Y_i,\psi_i\right>_{i<\om}$ be as in the proof of 
\ref{gapInTheMouse}. Claim \ref{claim:union} holds. Let $z = w_1$, $G=p$, 
$\Upsilon=\Upsilon^\M$, $U=U^\M$ and $U'=U'^\M$. Define the language \[ 
\Ll^*=\Ll\un\{\dot{\beta}_i,\dot{\M}_i\}_{i<\om}\un\{\dot{G},
\dot{p},\dot{W},\dot{z},\dot{\Upsilon},\dot{U},\dot{U'}\};\]
each symbol in $\Ll^*\cut\Ll$ is a constant symbol.  
Relative to these definitions, let $B_0$, $\left<B^i_0\right>_{i<\om}$ and 
$\vec{S}=\left<S_i\right>_{i<\om}$ be as in \cite{Scalesweakgap}.\footnote{As before, we use 
the symbol $\Ll^*$ where \cite{Scalesweakgap} uses $\Ll$, and vice versa.} The analogue of 
\cite[Corollary 
1.14]{Scalesweakgap} holds (the proof should be executed in $\J(\M)$, where we have 
$\left<S_i\right>_{i<\om}$, and where $\DC_\RR$ holds -- this allows us to ``intersect all the 
cones'' without introducing new reals, and also the resulting iterate $\N$ is in $\J(\M)$, hence in 
$\M$). Regarding \cite[Lemma 
1.15.1]{Scalesweakgap}, the overall proof is executed in $V$,
where $\M$ is countable, and so we may take 
$\Mbar=\M$, and we need not take any countable substructure of $V$.
The proper segments of the iteration are all in $\M$. Also see 
\cite{local_prikry} for details on the process of interleaving 
comparison with genericity iteration.\footnote{The issue is as follows. Let $\Tt$ be 
one of the trees involved in the comparison. Let $\alpha<\lh(\Tt)$; it might be that 
$[0,\alpha]_\Tt$ drops. But then the usual procedure for choosing the least extender on 
$\es_+(\M^\Tt_\alpha)$ producing a bad extender algebra axiom need not make sense, because in 
fact, the relevant extender algebra is not even in $M^\Tt_\alpha$.} Consider the
analogue of \cite[Lemma 1.16.2]{Scalesweakgap}:
\begin{lemma}
\label{1.16.2}
Let $\N$ be $\omega$-suitable and $\vec{S}$-iterable. Let 
$\pi\maps\Q \rightarrow \N$
be $\Sigma_1$-elementary with $\tau^\N_{i,j} \in \rg(\pi)$ for all
$i,j<\omega$. Then there is some $m<\omega$ such that for all $n\geq m$,
$\rg(\pi)$ is cofinal in $\delta_n^\N$.
\end{lemma}
\begin{proof}
The proof mostly follows that of \cite[1.16.2]{Scalesweakgap}. But consider the proof of its 
Claim; we adopt the same notation. Within that proof, consider the proof that $\M^*=\Mbar$. We 
prove this, as things are different. As $\M$ is countable we
have $\Mbar=\M$ and $\bar{\RR}=\RR^\M$.
Let $\Upsilon^*,U^*$, etc, be $\dot{\Upsilon}^{\M^*}$, $\dot{U}^{\M^*}$, etc.
Let $\Upsilon=\Upsilon^\M$ and $\Upsilon^-=\Upsilon^{\M^-}$, etc. We have $\rho:\M^-\to\M$ and 
$\psi^*:\Hh^*\to\Hh^-$.

First note that $\Upsilon^*=\Upsilon$, for $\rho\com\psi^*$ yields 
order-preserving maps 
$U^*\to U$ and ${U'}^*\to U'$. Therefore $\hmb^{\M^*}=\hmb^{\M}$. So essentially as in the 
proof of 
\ref{gapInTheMouse}, $\M^*$ is a $1$-sound hpm over $\hmb^{\M}$ with 
$\rho_1^{\M^*}=\om$ and $p_1^{\M^*}=p$.

By \ref{lem:gF_props}, as $\rho^*\com\psi^*\maps\Hh^*\to\Hh$ is $\Sigma_1$-elementary, we 
have that $\Hh^*$ is a $(0,\om_1+1)$-iterable $\g$-organized $\Omega$-pm over 
$T^{\M^*}$; likewise for $\Hh^{\M^*|\eta}$ for every $\eta$ such that $\M^*|\eta$ is relevant.
So $\M^*$ is a $(0,\om_1+1)$-iterable $\Theta$-$\g$-organized $\Omega$-pm over 
$\Upsilon^\M$.
So we can compare $\M^*$ with $\M$. 
Because they are both $1$-sound and minimal for realizing $\Sigma$, $\M^*=\M$.
\end{proof}

We modify the statement of \cite[Lemma 1.20.1]{Scalesweakgap} as follows: Let $\Q$ 
be 
$\om$-suitable, $j$-sound and $j$-realizable. We claim that with respect to trees above 
$\delta_{j-1}^\Q$, $\Q$ is short tree iterable, and the conclusions of \cite[Lemma 
1.20.1]{Scalesweakgap} hold, except with (a)(ii) replaced by 
``$\Q$-to-$\P$ drops'', and (b)(ii) replaced by ``$b$ drops and $\Tt\conc b$ is 
$\Gamma$-guided''. Let us argue that $\Q$ is short tree iterable above 
$\delta_{j-1}^\Q$. Assume $j=0$ for simplicity. First note that whenever $\pi:\Q\to\N$ is a 
$0$-realization, the $\pi$-pullback $(\Psishort_\N)^\pi$ of $\Psishort_\N$ 
is suitability strict. To see this argue like in the proof of \ref{clm:Sigma_N_ss}.
Then, as in the proof of \ref{lem:short_tree_strat_tracks}, it follows that 
$(\Psishort_\N)^\pi$ is precisely the short tree strategy for $\Q$. This suffices.
Now consider the uniqueness of the branch $b$ 
described in \cite[Lemma 1.20.1(b)]{Scalesweakgap}, as modified above. Given two such branches 
$b,c$, we compare 
the phalanxes $\Phi(\Tt\conc b),\Phi(\Tt\conc c)$, producing trees 
$\Uu,\Vv$. If $\Tt$ is short then note that both $\Tt\conc b$ and $\Tt\conc c$ are 
$\Gamma$-guided, so $b=c$. If $\Tt$ is maximal then $b,c$ cannot drop; rule out the possibility 
that, for example, $M_\infty^\Uu\pins M_\infty^\Vv$ and $b^\Vv$ drops, by using suitability 
strictness.

Let $\Sigma,\Q,(\Ff,\prec^*),\Q_\infty$ be defined as in \cite[\S2]{Scalesweakgap}. Then
$\Sigma,{(\Ff,{\prec^*})}\in\J(\M)$ and the analogue of \cite[Lemma 2.1.2]{Scalesweakgap} holds, 
but 
we 
mention some points. It 
seems possible that $\Q_\infty$ be illfounded because $\OR(\J(\M))=\OR(\M)+\om$. But
$\J(\M)\sats$``$\Q_\infty$ is wellfounded in the codes''. Standard arguments therefore show that 
$\Q_\infty|\delta_0^{\Q_\infty}$ is wellfounded (in fact that
$\delta_0^{\Q_\infty}\leq\Theta^\M$).\footnote{\label{ftn:motivate_Theta-g_2}Recall that at the 
start of the proof we reduced to 
the case that $\M\sats$``$\Theta$ exists''.
This reduction relied on $\M$ being
$\Theta$-$\g$-organized. This seems to be a key point at which there is a problem with the scales 
analysis for $\g$-organized mice.} The latter is enough for the scale construction to go through.
The rest of the argument is essentially as in \cite{Scalesweakgap}. This completes the 
proof.
\end{proof}

\subsection{Scales analysis within core model induction}

We finish by explaining how we use the scale existence theorems in application to the core 
model induction. Assume $\DC_\RR$.

Suppose that $\Upsilon\eqdef(\Omega\rest\HC)\cross\{z\}$ is self-scaled for some $z\in\RR$, with 
$z\geq_T a_0$. Then using the 
scales existence theorems \ref{passiveGameScale}, \ref{thm:self_anal_coded}, \ref{weak gap from 
MC} together with \ref{rem:non_adm_gap_dealt_with}, we get the scales analysis for 
$\Lp^{\GOmega}(\RR,\Upsilon)$  from optimal determinacy and super-small mouse capturing 
hypotheses (that is, through any initial segment of $\Lp^{\GOmega}(\RR,\Upsilon)$ for which these 
hypotheses hold).

We have dealt with $\Lp^{\GOmega}(\RR,\Omega\rest\HC,z)$ instead of $\Lp^{\GOmega}(\RR)$
because we seem to need extra assumptions to obtain the scales analysis from optimal assumptions in 
the 
latter. We now discuss what we need for this. In application, \emph{if} there are no 
divergent $\AD$ pointclasses, $\Omega$ will in fact be 
\emph{very} nice:

\begin{definition}\label{dfn:very_nice}
 Let $\utilde{\Gamma}$ be a boldface pointclass and $X\sub\RR$. We say that 
$\utilde{\Gamma}$ 
is an \dfnemph{$\AD$-pointclass} iff $\AD$ holds with respect to all sets in $\utilde{\Gamma}$. We 
say that $\utilde{\Gamma},X$ are \dfnemph{Wadge compatible} iff $A,X$ are Wadge 
compatible for every $A\in\utilde{\Gamma}$.

We say that $\Omega$ is \dfnemph{very nice} iff $\Upsilon\eqdef(\Omega\rest\HC)\cross\{z\}$ is 
self-scaled for some $z\in\RR$, $\J(\HC,\Upsilon)\sats\AD$, and $\Upsilon^\HCc$ is
Wadge compatible with every boldface $\AD$-pointclass.
\end{definition}

\begin{remark}\label{rem:application}
Suppose $\Omega$ is very nice and let $\Upsilon$ be as above.
We want to see that the scales analysis in 
$\Lp^{\GOmega}(\RR)$ proceeds from optimal determinacy assumptions. Let 
$\N\pins\Lp^{\GOmega}(\RR)$ be such that $\N\sats\AD$ and 
$\N$ ends a gap 
$[\alpha,\beta]$ of $\Lp^{\GOmega}(\RR)$, such that $[\alpha,\beta]$ is not strong. Suppose 
that 
if $[\alpha,\beta]$ is weak and $\Omega\rest\HC\in\N|\alpha$ then super-small mouse capturing 
for $\Gamma=\Sigma_1^{\N|\alpha}$ holds on a cone. We claim that one of the scale existence 
theorems 
\ref{passiveGameScale}, 
\ref{gapInTheMouse}, or \ref{weak gap from MC} applies.

For by \ref{rem:non_adm_gap_dealt_with} and the mouse capturing hypothesis, we may assume that the 
gap is admissible, and so weak, and that $\Omega\rest\HC\notin\N|\alpha$, so 
$\Upsilon^\HCc\notin\N|\alpha$. We 
claim that then 
$\J(\N)\sats\AD$, so \ref{gapInTheMouse} applies. If every set of reals in $\J(\N)$ is 
Wadge below $\Upsilon^\HCc$, this is because 
$\J(\HC,\Upsilon)\sats\AD$. So suppose otherwise. Let $\P\ins\N$ be least such that there 
is 
$Z\in\J(\P)$ such that $Z\not\leq_W\Upsilon^\HCc$. If $\P\pins\N$ 
then $\J_1(\P)\sats\AD$, so by the 
Wadge compatibility given by \ref{dfn:very_nice}, we have 
$\Omega\rest\HC\in\J(\P)$, so 
$\alpha\leq\l(\P)$. We claim that $\Omega\rest\HC\notin\N|\beta$. Because $\Omega$ is very 
nice 
and by 
\ref{InsideGap}, this is clear if $\Th_{\rPi_1}^{\N|\alpha}\leq_W\Upsilon^\HCc$ or 
$\Th_{\rSigma_1}^{\N|\alpha}\leq_W\Upsilon^\HCc$ (as very niceness would otherwise yield scales on 
these sets). Otherwise, by Wadge compatibility, 
$\Upsilon^\HCc<_W\Th_{\rSigma_1}^{\N|\alpha}$. But then because $\N|\alpha$ is admissible, 
$\Upsilon^\HCc\in\N|\alpha$, contradiction. So $\P=\N$. 
Since $\N$ ends a 
weak gap, there are sets $X_i\in\pow(\RR)\inter\N$ such that $\pow(\RR)\inter\J(\N)$ is exactly 
the sets which are projective in $\oplus_{i<\om}X_i$. It follows that
$\pow(\RR)\inter\J(\N)\sub\pow(\RR)\inter\J(\HC,\Upsilon)$, so 
$\J(\N)\sats\AD$ (and so $\Upsilon\in\J(\N)$).
\end{remark}

\bibliographystyle{plain}
\bibliography{Rmicebib}

\end{document}